%% file: Theta.tex
\documentclass[12pt,reqno]{amsart}
\def\draft{n}
\def\published{n}
\usepackage{
  amsmath, amssymb, multicol, mathtools, stmaryrd, xcolor,
  colortbl, mathabx, needspace, lgreek, sidecap, 
  blindtext, picins, dbnsymb, cjhebrew
}
\usepackage[export]{adjustbox} 
\usepackage{enumitem}
\usepackage{utfsym}     
\usepackage[all]{xy}
\usepackage[textwidth=6.5in,textheight=9in,headsep=0.15in,centering]{geometry}

\usepackage[pagebackref,breaklinks=true]{hyperref}\hypersetup{colorlinks,
  linkcolor={green!50!black},
  citecolor={green!50!black},
  urlcolor=blue
}

\usepackage{chngcntr}
\counterwithin{figure}{section}
\counterwithin{table}{section}

\usepackage{caption}
\if\draft y
  \usepackage[color]{showkeys}
\fi

\usepackage[stable]{footmisc} 

\input defs.tex

\begin{document} 
\newdimen\captionwidth\captionwidth=\hsize
\setcounter{secnumdepth}{4}

\title{A Fast, Strong, Topologically Meaningful, and Fun Knot Invariant}

\author{Dror~Bar-Natan}
\address{
  Department of Mathematics\\
  University of Toronto\\
  Toronto Ontario M5S 2E4\\
  Canada \ifpub{(corresponding author)}{}
}
\email{drorbn@math.toronto.edu}
\urladdr{http://www.math.toronto.edu/drorbn}

\author{Roland~van~der~Veen}
\address{
  University of Groningen, Bernoulli Institute\\
  P.O. Box 407\\
  9700 AK Groningen\\
  The Netherlands \ifpub{(corresponding author)}{}
}
\email{roland.mathematics@gmail.com}
\urladdr{http://www.rolandvdv.nl/}

\date{First edition September 22, 2025. This edition \today.}

\makeatletter
\@namedef{subjclassname@2020}{\textup{2020} Mathematics Subject Classification}
\makeatother
\subjclass[2020]{Primary 57K14, secondary 16T99}
\keywords{
  Alexander polynomial,
  loop expansion,
  solvable approximation,
  knot genus,
  fibered knots,
  ribbon knots,
  polynomial time computations,
  Feynman diagrams,
  perturbed Gaussian integration,
  Seifert formulas%
}

\thanks{This paper is available in electronic
form, along with source files and a demo {\sl Mathematica}
notebook at \url{http://drorbn.net/Theta} and at \arXiv{2509.18456}.}

\setlist[itemize]{left=0pt .. 12pt}
\begin{abstract}
\input abstract.tex
\end{abstract}

\setlist[itemize]{left=0pt .. \parindent}
\setlist[enumerate]{left=0pt .. 1.5\parindent}

\maketitle

\setcounter{tocdepth}{3}
\vspace{-10mm}
\tableofcontents

\input body.tex

\draftcut


\parpic[l]{\parbox{1.5in}{
  \includegraphics[width=1.5in]{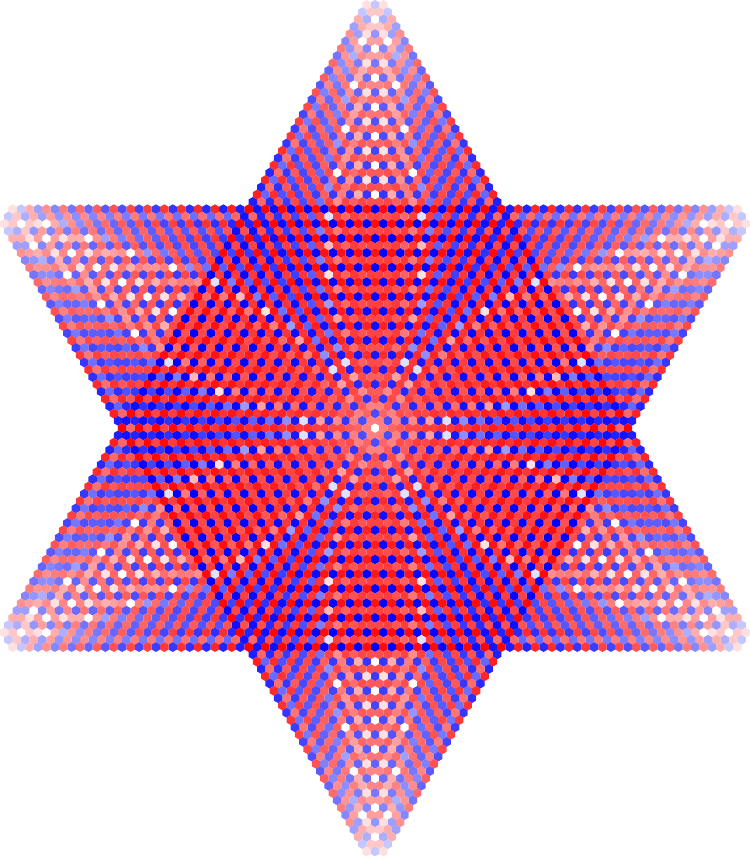}
  \vskip -8mm
  \tiny Pretzel
  \newline$(2,31,-31)$
}}

\end{document}

\endinput

%% file: defs.tex
\newcommand\ifpub[2]{\if\published y #1 \else #2 \fi}

\theoremstyle{plain}
\newtheorem{theorem}{Theorem}

\newtheorem{conjecture}[theorem]{Conjecture}

\newtheorem{dream}[theorem]{Dream}
\newtheorem{fact}[theorem]{Fact}
\newtheorem{lemma}[theorem]{Lemma}

\newtheorem{proposition}[theorem]{Proposition}

\theoremstyle{definition}

\newtheorem{question}[theorem]{Question}

\newtheorem{comment}[theorem]{Comment}

\newtheorem{discussion}[theorem]{Discussion}

\newtheorem{remark}[theorem]{Remark}

\newlength{\standardunitlength}
\newcommand{\adj}{\operatorname{adj}}

\def\car{\reflectbox{\usym{1F697}}}
\def\rac{\usym{1F697}}

\def\qed{{\linebreak[1]\null\hfill\text{$\Box$}}}
\def\endpar#1{~\hfill\fbox{\footnotesize#1}}

\newlength{\globalparindent}
\def\arXiv#1{{\href{https://arxiv.org/abs/#1}{arXiv:\linebreak[0]#1}}}

\def\weburl{http://drorbn.net/SL2PO}
\def\web#1{\href{\weburl/#1}{\begin{greek}web\end{greek}/#1}}

\def\bbQ{{\mathbb Q}}
\def\bbR{{\mathbb R}}

\def\calD{{\mathcal D}}
\def\calE{{\mathcal E}}

\def\calU{{\mathcal U}}

\def\frakb{{\mathfrak b}}

\def\frakg{{\mathfrak g}}
\def\frakh{{\mathfrak h}}

\def\tilg{{\tilde{g}}}

\def\Kh{\text{\it Kh}}

\def\Vol{\text{\it Vol}}

\def\draftcut{\if\draft y \cleardoublepage \fi}

\definecolor{lightred}{RGB}{255, 217, 217}

\def\blue{\color{blue}}

\def\red{\color{red}}

\def\imagetop#1{\vtop{\null\hbox{#1}}}

\def\eps{\epsilon}

\makeatletter
\let\DOTSI\relax 
\newcommand*{\Gint}{%
  \DOTSI
  \mathop{%
    \mathpalette\@LetterOnInt{G}%
  }%
  \mkern-\thinmuskip 
  \int
}
\newcommand*{\@LetterOnInt}[2]{%
  \sbox0{$#1\int\m@th$}%
  \sbox2{$%
    \ifx#1\displaystyle
      \textstyle
    \else
      \scriptscriptstyle
    \fi
    #2%
  \m@th$}%
  \dimen@=.4\dimexpr\ht0+\dp0\relax
  \ifdim\dimexpr\ht2+\dp2\relax>\dimen@
    \sbox2{\resizebox*{!}{\dimen@}{\unhcopy2}}%
  \fi
  \dimen@=\wd0 %
  \ifdim\wd2>\dimen@
    \dimen@=\wd2 %
  \fi
  \rlap{\hbox to \dimen@{\hfil
    $#1\vcenter{\copy2}\m@th$%
  \hfil}}%
  \ifdim\dimen@>\wd0 %
    \kern.5\dimexpr\dimen@-\wd0\relax
  \fi
}
\makeatother

\def\ip{{i^+}} \def\jp{{j^+}} \def\kp{{k^+}}
\def\ipp{{i^{+\!+}}} \def\jpp{{j^{+\!+}}} \def\kpp{{k^{+\!+}}}

\def\face{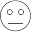}
\def\human{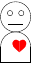}
\def\machine{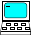}

\ifpub{\def\cellscale{0.76}}{\def\cellscale{1}}
\newsavebox\pdfbox
\newlength{\pdfheight}

\def\nbpdfPostInput{\vskip 2mm}
\def\nbpdfInput#1{{%
  \savebox{\pdfbox}{\includegraphics[scale=\cellscale, max width=\linewidth-0.27in]{#1}}%
  \settoheight{\pdfheight}{\usebox{\pdfbox}}%
  \noindent\imagetop{\ifdim\pdfheight<10mm\input{figs/face.pdf_t}\else\input{figs/human.pdf_t}\fi}\ %
  \imagetop{\usebox{\pdfbox}}%
  \nbpdfPostInput%
}}

\def\nbpdfPreOutput{~}
\def\nbpdfOutput#1{{\noindent{\imagetop{\input{figs/machine.pdf_t}}\nbpdfPreOutput\imagetop{\includegraphics[scale=\cellscale, max width=\linewidth-0.27in]{#1}}\vskip 2mm}}}

\def\nbpdfMessage#1{{}}
\def\nbpdfPrint#1{{\noindent{\imagetop{\input{figs/machine.pdf_t}}\ \imagetop{\includegraphics[scale=\cellscale]{#1}}\vskip 2mm}}}

%% file: figs/face.pdf_t
\begin{picture}(0,0)%
\includegraphics{figs/face.pdf}%
\end{picture}%
%
%
\setlength{\unitlength}{789sp}%
\begingroup\makeatletter\ifx\SetFigFont\undefined%
\gdef\SetFigFont#1#2#3#4#5{%
  \reset@font\fontsize{#1}{#2pt}%
  \fontfamily{#3}\fontseries{#4}\fontshape{#5}%
  \selectfont}%
\fi\endgroup%
\begin{picture}(1240,1240)(-19,-381)
\end{picture}%

%% file: figs/human.pdf_t
\begin{picture}(0,0)%
\includegraphics{figs/human.pdf}%
\end{picture}%
%
%
\setlength{\unitlength}{789sp}%
\begingroup\makeatletter\ifx\SetFigFont\undefined%
\gdef\SetFigFont#1#2#3#4#5{%
  \reset@font\fontsize{#1}{#2pt}%
  \fontfamily{#3}\fontseries{#4}\fontshape{#5}%
  \selectfont}%
\fi\endgroup%
\begin{picture}(1244,2442)(-21,-1583)
\end{picture}%

%% file: figs/machine.pdf_t
\begin{picture}(0,0)%
\includegraphics{figs/machine.pdf}%
\end{picture}%
%
%
\setlength{\unitlength}{789sp}%
\begingroup\makeatletter\ifx\SetFigFont\undefined%
\gdef\SetFigFont#1#2#3#4#5{%
  \reset@font\fontsize{#1}{#2pt}%
  \fontfamily{#3}\fontseries{#4}\fontshape{#5}%
  \selectfont}%
\fi\endgroup%
\begin{picture}(1244,1394)(-21,-608)
\end{picture}%

%% file: abstract.tex
In this paper we discuss a pair of polynomial knot invariants
$\Theta=(\Delta,\theta)$ which is:
\begin{itemize}
\item Theoretically and practically fast: $\Theta$ can be computed in polynomial time.
  We can compute it in full on random knots with over 300 crossings,
  and its evaluation at simple rational numbers on random knots with over
  600 crossings.
\item Strong: Its separation power is much greater than the hyperbolic volume, the
  HOMFLY-PT polynomial and Khovanov homology (taken together) on knots
  with up to 15 crossings (while being computable on much larger knots).
\item Topologically meaningful: It gives a genus bound, and there are
  reasons to hope that it would do more.
\item Fun: Scroll to Figures 1.1--1.4, 3.1, and 6.2.
\end{itemize}
$\Delta$ is merely the Alexander polynomial. $\theta$ is almost
certainly equal to an invariant that was studied extensively by
Ohtsuki \cite{Ohtsuki:TwoLoop}, continuing Rozansky, Kricker, and
Garoufalidis \cite{Rozansky:Contribution, Rozansky:Burau, Rozansky:U1RCC,
Kricker:Lines, GaroufalidisRozansky:LoopExpansion}.  Yet our formulas,
proofs, and programs are much simpler and enable its computation even
on very large knots.

%% file: body.tex
\draftcut
\section{Fun} The word ``fun'' rarely appears in the title of a math
paper, so let us start with a brief justification (more can be found in
a Quanta Magazine article,~\cite{Klarreich:QRCode}).

$\Theta$ is a pair of polynomials. The first, $\Delta$, is old news, the
Alexander polynomial~\cite{Alexander:TopologicalInvariants}. It is
a one-variable Laurent polynomial in a variable $T$. For example,
$\Delta(\lefttrefoil)=T^{-1}-1+T$. We turn such a polynomial into
a list of coefficients (for $\lefttrefoil$, it is $(1,-1,1)$),
and then to a chain of bars of varying colours: white for the zero
coefficients, and red and blue for the positive and negative coefficients
(with intensity proportional to the magnitude of the coefficients).
The result is a ``bar code'', and for the trefoil $\lefttrefoil$ it is
{\red\rule{2mm}{3mm}}{\blue\rule{2mm}{3mm}}{\red\rule{2mm}{3mm}}.

\parpic[r]{\includegraphics[width=1.2in]{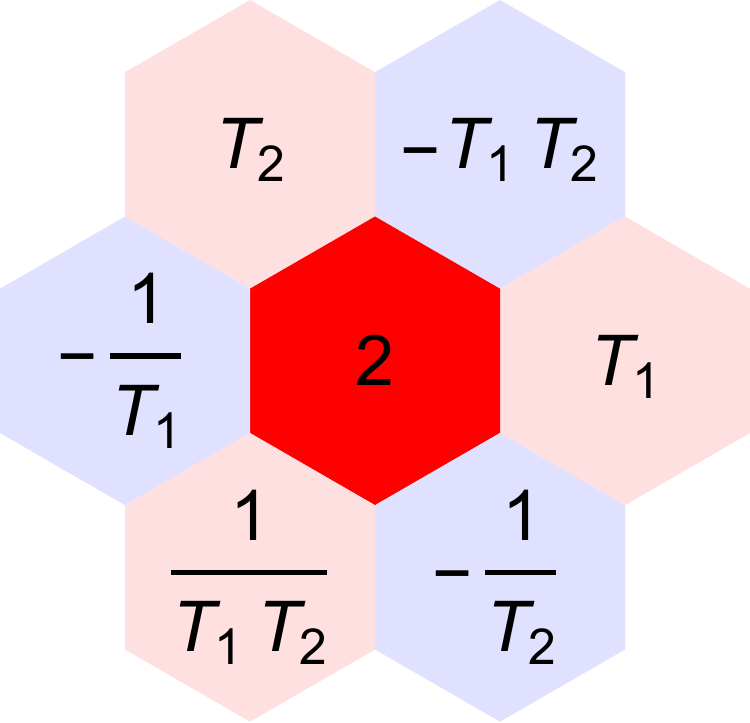}}
Similarly, $\theta$ is a 2-variable Laurent polynomial, in variables
$T_1$ and $T_2$. We can turn such a polynomial into a 2D array of
coefficients and then using the same rules, into a 2D array of colours,
namely, into a picture. To highlight a certain conjectured hexagonal
symmetry of the resulting pictures, we apply a shear transformation to
the plane before printing. So a monomial $cT_1^{n_1}T_2^{n_2}$ gets
printed at position $(n_1-n_2/2, \sqrt{3}n_2/2)$ instead of the more
straightforward $(n_1,n_2)$. On the right is the 2D picture corresponding
to the polynomial $2+T_1-T_1T_2+T_2-T_1^{-1}+T_1^{-1}T_2^{-1}-T_2^{-1}$.

\begin{figure}
\adjustbox{valign=m}{\resizebox*{!}{8.9in}{
\def\k#1#2{{\parbox{0.5in}{\centering
  \vskip 2pt
  \includegraphics[width=\linewidth,height=\linewidth]{KnotFigs/#1_#2.pdf}
  \newline\href{https://katlas.org/wiki/#1_#2}{#1\_#2}
  \vskip 2pt
}}}
\begin{tabular}{|c|c|c|c|c|c|c|c|c|c|}
  \hline \k{0}{1} & \k{3}{1} & \k{4}{1} & \k{5}{1} & \k{5}{2} & \k{6}{1} & \k{6}{2} & \k{6}{3} & \k{7}{1} & \k{7}{2} \\
  \hline \k{7}{3} & \k{7}{4} & \k{7}{5} & \k{7}{6} & \k{7}{7} & \k{8}{1} & \k{8}{2} & \k{8}{3} & \k{8}{4} & \k{8}{5} \\
  \hline \k{8}{6} & \k{8}{7} & \k{8}{8} & \k{8}{9} & \k{8}{10} & \k{8}{11} & \k{8}{12} & \k{8}{13} & \k{8}{14} & \k{8}{15} \\
  \hline \k{8}{16} & \k{8}{17} & \k{8}{18} & \k{8}{19} & \k{8}{20} & \k{8}{21} & \k{9}{1} & \k{9}{2} & \k{9}{3} & \k{9}{4} \\
  \hline \k{9}{5} & \k{9}{6} & \k{9}{7} & \k{9}{8} & \k{9}{9} & \k{9}{10} & \k{9}{11} & \k{9}{12} & \k{9}{13} & \k{9}{14} \\
  \hline \k{9}{15} & \k{9}{16} & \k{9}{17} & \k{9}{18} & \k{9}{19} & \k{9}{20} & \k{9}{21} & \k{9}{22} & \k{9}{23} & \k{9}{24} \\
  \hline \k{9}{25} & \k{9}{26} & \k{9}{27} & \k{9}{28} & \k{9}{29} & \k{9}{30} & \k{9}{31} & \k{9}{32} & \k{9}{33} & \k{9}{34} \\
  \hline \k{9}{35} & \k{9}{36} & \k{9}{37} & \k{9}{38} & \k{9}{39} & \k{9}{40} & \k{9}{41} & \k{9}{42} & \k{9}{43} & \k{9}{44} \\
  \hline \k{9}{45} & \k{9}{46} & \k{9}{47} & \k{9}{48} & \k{9}{49} & \k{10}{1} & \k{10}{2} & \k{10}{3} & \k{10}{4} & \k{10}{5} \\
  \hline \k{10}{6} & \k{10}{7} & \k{10}{8} & \k{10}{9} & \k{10}{10} & \k{10}{11} & \k{10}{12} & \k{10}{13} & \k{10}{14} & \k{10}{15} \\
  \hline \k{10}{16} & \k{10}{17} & \k{10}{18} & \k{10}{19} & \k{10}{20} & \k{10}{21} & \k{10}{22} & \k{10}{23} & \k{10}{24} & \k{10}{25} \\
  \hline \k{10}{26} & \k{10}{27} & \k{10}{28} & \k{10}{29} & \k{10}{30} & \k{10}{31} & \k{10}{32} & \k{10}{33} & \k{10}{34} & \k{10}{35} \\
  \hline \k{10}{36} & \k{10}{37} & \k{10}{38} & \k{10}{39} & \k{10}{40} & \k{10}{41} & \k{10}{42} & \k{10}{43} & \k{10}{44} & \k{10}{45} \\
  \hline \k{10}{46} & \k{10}{47} & \k{10}{48} & \k{10}{49} & \k{10}{50} & \k{10}{51} & \k{10}{52} & \k{10}{53} & \k{10}{54} & \k{10}{55} \\
  \hline \k{10}{56} & \k{10}{57} & \k{10}{58} & \k{10}{59} & \k{10}{60} & \k{10}{61} & \k{10}{62} & \k{10}{63} & \k{10}{64} & \k{10}{65} \\
  \hline \k{10}{66} & \k{10}{67} & \k{10}{68} & \k{10}{69} & \k{10}{70} & \k{10}{71} & \k{10}{72} & \k{10}{73} & \k{10}{74} & \k{10}{75} \\
  \hline \k{10}{76} & \k{10}{77} & \k{10}{78} & \k{10}{79} & \k{10}{80} & \k{10}{81} & \k{10}{82} & \k{10}{83} & \k{10}{84} & \k{10}{85} \\
  \hline \k{10}{86} & \k{10}{87} & \k{10}{88} & \k{10}{89} & \k{10}{90} & \k{10}{91} & \k{10}{92} & \k{10}{93} & \k{10}{94} & \k{10}{95} \\
  \hline \k{10}{96} & \k{10}{97} & \k{10}{98} & \k{10}{99} & \k{10}{100} & \k{10}{101} & \k{10}{102} & \k{10}{103} & \k{10}{104} & \k{10}{105} \\
  \hline \k{10}{106} & \k{10}{107} & \k{10}{108} & \k{10}{109} & \k{10}{110} & \k{10}{111} & \k{10}{112} & \k{10}{113} & \k{10}{114} & \k{10}{115} \\
  \hline \k{10}{116} & \k{10}{117} & \k{10}{118} & \k{10}{119} & \k{10}{120} & \k{10}{121} & \k{10}{122} & \k{10}{123} & \k{10}{124} & \k{10}{125} \\
  \hline \k{10}{126} & \k{10}{127} & \k{10}{128} & \k{10}{129} & \k{10}{130} & \k{10}{131} & \k{10}{132} & \k{10}{133} & \k{10}{134} & \k{10}{135} \\
  \hline \k{10}{136} & \k{10}{137} & \k{10}{138} & \k{10}{139} & \k{10}{140} & \k{10}{141} & \k{10}{142} & \k{10}{143} & \k{10}{144} & \k{10}{145} \\
  \hline \k{10}{146} & \k{10}{147} & \k{10}{148} & \k{10}{149} & \k{10}{150} & \k{10}{151} & \k{10}{152} & \k{10}{153} & \k{10}{154} & \k{10}{155} \\
  \hline \k{10}{156} & \k{10}{157} & \k{10}{158} & \k{10}{159} & \k{10}{160} & \k{10}{161} & \k{10}{162} & \k{10}{163} & \k{10}{164} & \k{10}{165} \\
  \hline
\end{tabular}
}}
$\underset{\Theta}{\rightarrow}$
\includegraphics[height=8.9in,valign=m]{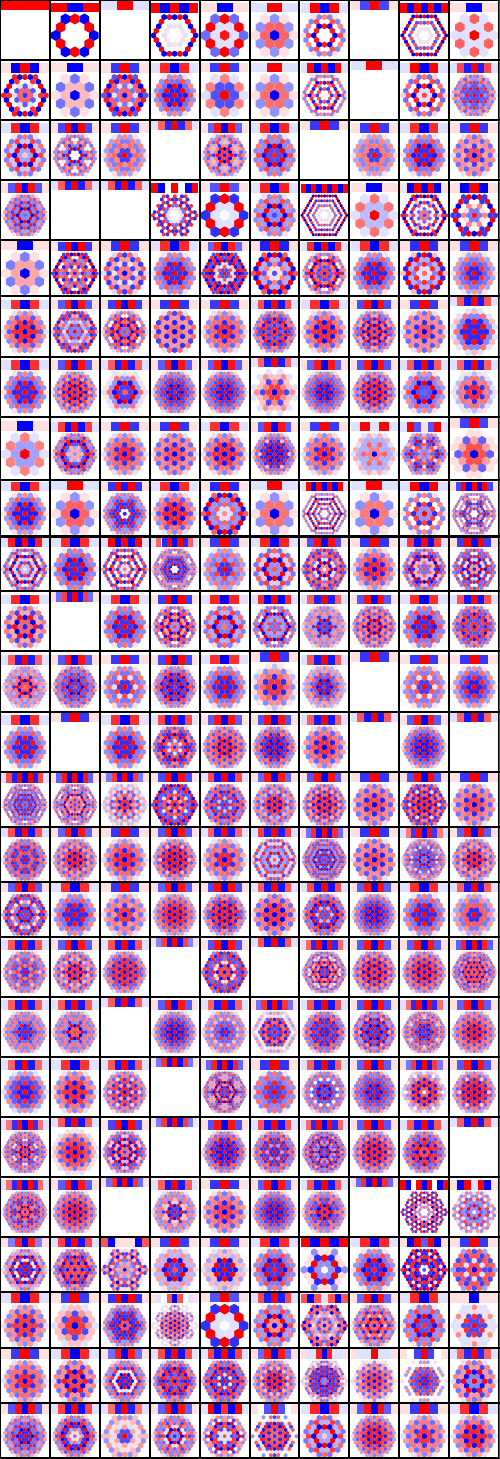}
\caption{$\Theta$ as a bar code and a QR code, for all the knots in the Rolfsen table.}
\label{fig:Rolfsen}
\end{figure}

\begin{figure}
\includegraphics[width=\linewidth]{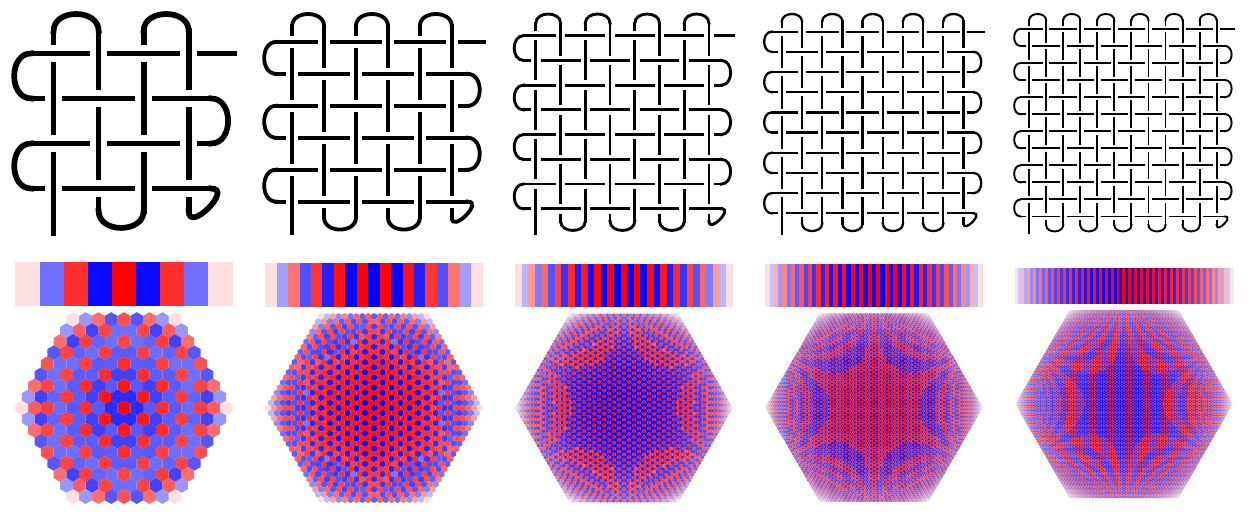}
\caption{$\Theta$ of some square weave knots, as computed by~\cite[WeaveKnots.nb]{Self}.}
\label{fig:SquareWeaves}
\end{figure}

\begin{figure}
\includegraphics[width=0.6\linewidth]{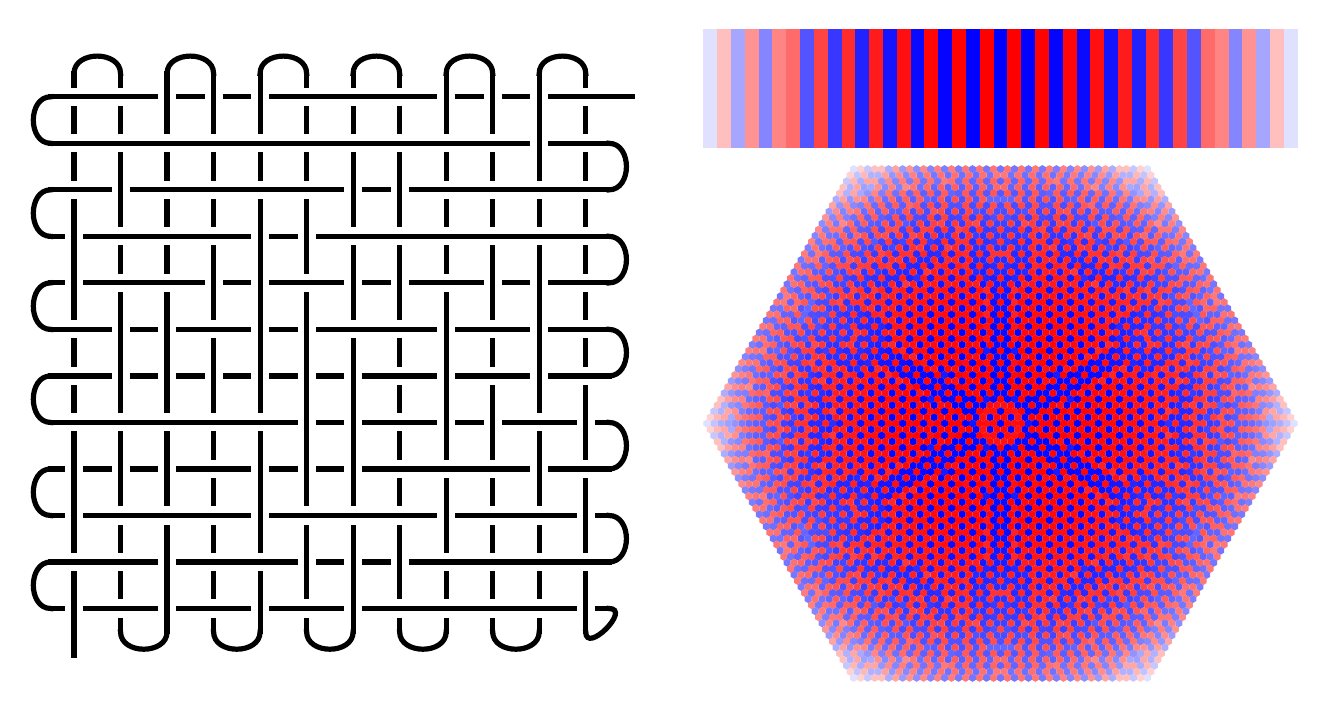}
\caption{$\Theta$ of a randomized weave knot, as computed by~\cite[WeaveKnots.nb]{Self}. Crossings
were chosen to be positive or negative with equal probabilities.}
\label{fig:RandomWeave}
\end{figure}

\begin{figure}
\includegraphics[width=\linewidth]{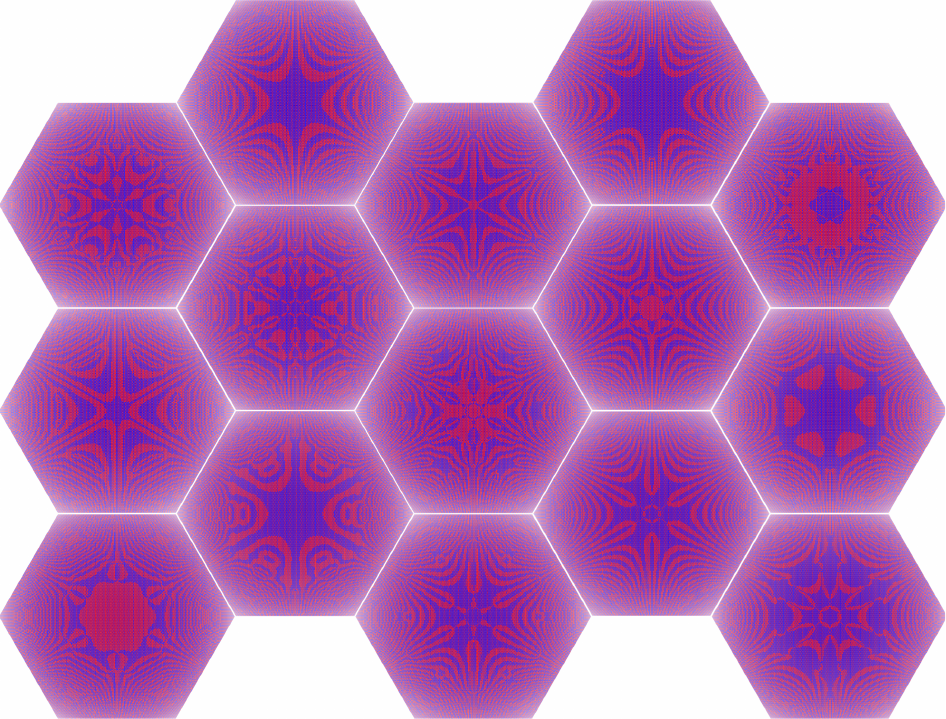}
\caption{$\theta$ (hexagonal QR code only) of the 15 largest knots that
  we have computed by September 16, 2024. They are all ``generic'' in as
  much as we know, and they all have $\geq 300$ crossings. The knots come from~\cite{DHOEBL:Random}.
  Warning: Some screens/printers may introduce spurious Moir\'e interference patterns.
}
\label{fig:300}
\end{figure}

Thus $\Theta$ becomes a pair of pictures: a bar code, and a 2D
picture that we call a ``hexagonal QR code''. For the knots in the
Rolfsen table (with the unknot prepended at the start), they are in
Figure~\ref{fig:Rolfsen}. For some alternating square weave knots, they
are in Figure~\ref{fig:SquareWeaves}, and for a random square weave, in
Figure~\ref{fig:RandomWeave}.  In addition, the hexagonal QR codes of 15
knots with $\geq 300$ crossings are in Figure~\ref{fig:300}, and $\Theta$
of a 132-crossing torus knot is in Figure~\ref{fig:T227}. Some further
computations and figures, also highlighting the parity of coefficients
rather than just their signs, are at~\cite{Lalani:ThetaSurvey}.

\parpic[r]{\parbox{3in}{\begin{center}
  \includegraphics[height=1in]{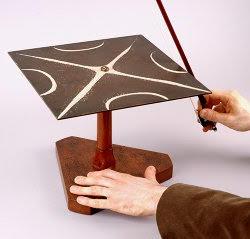}
  \includegraphics[height=1in]{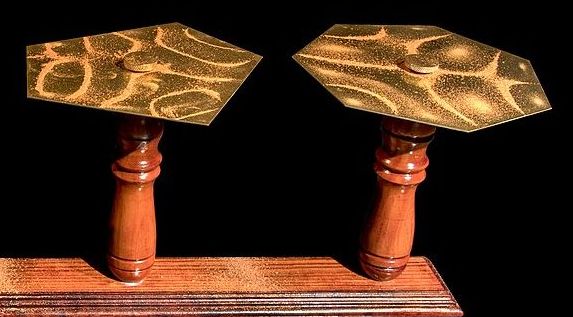}
  \newline \tiny
    left: \copyright\ 
    \href{https://www.whipplemuseum.cam.ac.uk/explore-whipple-collections/acoustics/ernst-chladni-physicist-musician-and-musical-instrument-maker}{Whipple Museum of the History of Science, University of Cambridge}; 
    right:
    \href{https://creativecommons.org/licenses/by-sa/4.0/deed.en}{CC-BY-SA 4.0}
    / \href{https://en.wikipedia.org/wiki/Chladni\%27s_law}{Wikimedia}
    / Matemateca (IME USP) / Rodrigo Tetsuo Argenton
\end{center}}}
Clearly there are patterns in these figures. There is a
hexagonal symmetry and the QR codes are nearly always hexagons
(these are independent properties). Much more can be seen in
Figure~\ref{fig:Rolfsen}. In Figure~\ref{fig:300} there seem to be
large-scale patterns perhaps reminiscent of the ``Chladni figures''
formed by powders atop vibrating plates (on right).  We can't prove any of these
things, and the last one, we can't even formulate properly. Yet they
are clearly there, too clear to be the result of chance alone.

We plan to have fun over the next few years observing and proving these patterns. We hope
that others will join us too.

\draftcut
\needspace{60mm}
\section{The Main Theorem} \label{sec:MainTheorem}

\parpic[r]{\parbox{2in}{
  \begin{center}\scalebox{1}{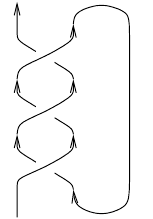}\end{center}
  \captionof{figure}{An example upright knot diagram.} \label{fig:SampleDiagram}
}}
We start with the definition of $\Theta$.
Given an oriented $n$-crossing knot $K$, we draw it in the plane as a
long knot diagram $D$ in such a way that the two strands intersecting
at each crossing are pointing up (that's always possible because we
can always rotate crossings as needed), and so that at its beginning
and at its end the knot is oriented upward. We call such a diagram an
{\em upright knot diagram.} An example of an upright knot diagram
is shown on the right.

We then label each edge of the diagram with two labels: a running index
$k$ which runs from 1 to $2n+1$, and a ``rotation number'' $\varphi_k$,
the geometric rotation number of that edge\footnote{The signed number
of times the tangent to the edge is horizontal and heading right,
with cups counted with $+1$ signs and caps with $-1$; this number
is well defined because at their ends, all edges are headed up.}. In
Figure~\ref{fig:SampleDiagram} the running index runs from $1$ to $7$,
and the rotation numbers for all edges are $0$ (and hence are omitted)
except for $\varphi_4$, which is $-1$.

Let $X$ be the set of all crossings in the diagram $D$, where we encode each crossing as a
triple (sign of the crossing, incoming over edge, incoming under edge). In our example we have
$X=\{(1,1,4),(1,5,2),(1,3,6)\}$.

We let $A$ be the $(2n+1)\times(2n+1)$ matrix of Laurent polynomials in 
a variable $T$, defined by
\[ A \coloneqq
  I - \sum_{c=(s,i,j)\in X} \left( T^sE_{i,i+1} + (1-T^s)E_{i,j+1} + E_{j,j+1} \right),
\]
where $I$ is the identity matrix and $E_{\alpha\beta}$ denotes the elementary
matrix with $1$ in row $\alpha$ and column $\beta$ and zeros elsewhere.

Alternatively, $A=I + \sum_c A_c$, where $A_c$ is a matrix of zeros except for the blocks as follows:
\begin{equation} \label{eq:A}
  \begin{array}{c}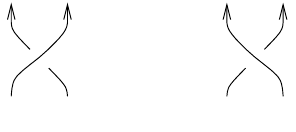\end{array}
  \qquad\longrightarrow\qquad
  \begin{array}{c|cccc}
    A_c &   \text{column }i+1  &  \text{column }j+1 \\
    \hline
    \text{row }i & -T^s  & T^s-1 \\
    \text{row }j & 0  & -1
  \end{array}
\end{equation}

We note that the determinant of $A$ is equal up to a unit to the
normalized Alexander polynomial $\Delta$ of $K$.\footnote{The informed
reader will note that $A$ is a presentation matrix for the Alexander
module of $K$, obtained by using Fox calculus on the Wirtinger
presentation of the fundamental group of the complement of $K$.} In fact,
we have that
\begin{equation} \label{eq:Delta}
  \Delta = \Delta(K) = T^{(-\varphi(D)-w(D))/2}\det(A),
\end{equation}
where $\varphi(D)\coloneqq\sum_k\varphi_k$ is the total rotation number
of $D$ and where $w(D)=\sum_cs_c$ is the writhe of $D$, namely the sum
of the signs $s_c$ of all the crossings $c$ in $D$.

We let $G=(g_{\alpha\beta})=A^{-1}$, and, thinking of it as a function
$g_{\alpha\beta}$ of a pair of edges $\alpha$ and $\beta$, we call it
the Green function of the diagram $D$. When inspired by physics (e.g.\
Fact~\ref{fact:PerturbedGaussian} and~\cite{IType}) we sometimes call
it ``the 2-point function'', and when thinking of car traffic (e.g.\
Comment~\ref{com:traffic} and~\cite{APAI, Toronto-241030}) we sometimes
call it ``the traffic function''. As an example, here are $A$ and $G$
for the knot diagram $D$ of Figure~\ref{fig:SampleDiagram}:
\[ \footnotesize \setlength{\arraycolsep}{2pt}
\left(
  \begin{array}{ccccccc}
   1 & -T & 0 & 0 & T-1 & 0 & 0 \\[2pt]
   0 & 1 & -1 & 0 & 0 & 0 & 0 \\[2pt]
   0 & 0 & 1 & -T & 0 & 0 & T-1 \\[2pt]
   0 & 0 & 0 & 1 & -1 & 0 & 0 \\[2pt]
   0 & 0 & T-1 & 0 & 1 & -T & 0 \\[2pt]
   0 & 0 & 0 & 0 & 0 & 1 & -1 \\
   0 & 0 & 0 & 0 & 0 & 0 & 1 \\
  \end{array}
\right),
  \quad \left(
    \begin{array}{ccccccc}
     1 & T & 1 & T & 1 & T & 1 \\
     0 & 1 & \frac{1}{T^2-T+1} & \frac{T}{T^2-T+1} &
       \frac{T}{T^2-T+1} & \frac{T^2}{T^2-T+1} & 1 \\
     0 & 0 & \frac{1}{T^2-T+1} & \frac{T}{T^2-T+1} &
       \frac{T}{T^2-T+1} & \frac{T^2}{T^2-T+1} & 1 \\
     0 & 0 & \frac{1-T}{T^2-T+1} & \frac{1}{T^2-T+1} &
       \frac{1}{T^2-T+1} & \frac{T}{T^2-T+1} & 1 \\
     0 & 0 & \frac{1-T}{T^2-T+1} & \frac{T-T^2}{T^2-T+1}
       & \frac{1}{T^2-T+1} & \frac{T}{T^2-T+1} & 1 \\
     0 & 0 & 0 & 0 & 0 & 1 & 1 \\
     0 & 0 & 0 & 0 & 0 & 0 & 1 \\
    \end{array}
  \right).
\]

Let $T_1$ and $T_2$ be indeterminates and let
$T_3\coloneqq T_1T_2$. Let $\Delta_\nu\coloneqq\Delta|_{T\to T_\nu}$
and $G_\nu=(g_{\nu\alpha\beta})\coloneqq G|_{T\to T_\nu}$ be $\Delta$
and $G$ subject to the substitution $T\to T_\nu$, where $\nu=1,2,3$.

Given crossings $c=(s,i,j)$, $c_0=(s_0,i_0,j_0)$, and $c_1=(s_1,i_1,j_1)$
in $X$ and an edge label $k$, let
\begin{alignat}{1}
  \label{eq:F1} F_1(c) =\ & s
    \left[1/2 - g_{3ii} +  T_2^s g_{1ii} g_{2ji} - T_2^s g_{3jj} g_{2ji}
      - (T_2^s-1) g_{3ii} g_{2ji} \right. \\
  \nonumber &\quad
    \left. + (T_3^s-1) g_{2ji} g_{3ji} - g_{1ii} g_{2jj} + 2 g_{3ii} g_{2jj}
      + g_{1ii} g_{3jj} - g_{2ii} g_{3jj} \right] \\
  \nonumber &
  + \frac{s}{T_2^s-1}
    \left[
      (T_1^s-1)T_2^s \left( g_{3jj} g_{1ji} - g_{2jj} g_{1ji} + T_2^s g_{1ji} g_{2ji} \right) \right. \\
  \nonumber &\quad
      + \left. (T_3^s-1) g_{3ji} \left( 1 - T_2^s g_{1ii} + g_{2ij} + (T_2^s-2) g_{2jj}
        - (T_1^s-1) (T_2^s+1) g_{1ji} \right)\right] \\
  \label{eq:F2} F_2(c_0,c_1) =\ &
    \frac{s_1 (T_1^{s_0}-1) (T_3^{s_1}-1) g_{1j_1i_0} g_{3j_0i_1}}{T_2^{s_1}-1}
    \left(T_2^{s_0} g_{2i_1i_0}+g_{2j_1j_0} - T_2^{s_0} g_{2j_1i_0}-g_{2i_1j_0} \right) \\
  \label{eq:F3} F_3(k) =\ & (g_{3kk}-1/2)\varphi_k
\end{alignat}

These formulas are uninspiring, yet they are easy to compute (given $G$), and they work:

\needspace{30mm}
\begin{theorem}[The Main Theorem, proof in Section~\ref{sec:Proof}] \label{thm:Main}
The following are knot invariants:
\begin{equation} \label{eq:Main}
  \theta_0(D) \coloneqq \sum_{c\in X} F_1(c) + \sum_{c_0,c_1\in X} F_2(c_0,c_1)
    + \sum_{\text{edges }k}F_3(k)
  \quad\text{and}\quad
  \theta(D) \coloneqq \Delta_1\Delta_2\Delta_3\theta_0(D).
\end{equation}
Furthermore, $\theta$ is a Laurent polynomial in $T_1$ and $T_2$, with integer coefficients.
\end{theorem}

Some comments are now in order:

\begin{comment} \label{com:normalization} The entries of $G_\nu$ are
rational functions with denominators $\Delta_\nu$, and so $\theta_0$
is valued in the ring of rational functions $\bbQ(T_1,T_2)$. The
point of $\theta$ is to clear these denominators by multiplying by
$\Delta_1\Delta_2\Delta_3$ so as to get an invariant valued in Laurent
polynomials. (There remains a potential denominator of the
form $(T_2-1)^{-1}$ coming from the explicit denominators in
Equations~\eqref{eq:F1} and~\eqref{eq:F2}. It will be shown to cancel
in Section~\ref{ssec:Residue}.)
\end{comment}

\begin{comment} \label{com:traffic}
We note following~\cite{APAI} that $g_{\alpha\beta}$ can be interpreted
as measuring ``car traffic'', assuming a stream of traffic is injected
near the start of edge $\alpha$ and a ``traffic counter'' is placed
near the end of edge $\beta$, and where cars always obey the following
traffic rules:
\begin{itemize}
\item Car travel on the edges of the knot, always in a direction consistent with the
  orientation of these edges.
\item When a car reaches a crossing on the under-strand, it travels through and continues on
  the other side.
\item When a car reaches a crossing of sign $s=\pm 1$ on the over-strand, it continues right
  through with probability $T^s$, yet with probability $1-T^s$ it falls down 
  and continues travelling on the lower strand. (It matters not that $T$ and $T^{-1}$ cannot
  be between $0$ and $1$ at the same time --- we merely use the algebraic rules of probability
  without caring about the inequalities that normally come with them).
\item When cars reach the ``end'' of the knot, the abyss that follows
  edge $2n+1$, they fall off the picture never to be seen again.
\end{itemize}

These rules can be summarized by the following pictures:
\[
  \def\DallEA{{\tiny image credits:}}
  \def\DallEB{{\tiny \href{https://labs.openai.com/}{Dall-E}}}
  \resizebox{\linewidth}{!}{\input{figs/CarRules2.pdf_t}}
\]

For further details, see~\cite{APAI}. \endpar{\ref{com:traffic}}
\end{comment}

\begin{comment} \label{com:FeynmanDiagrams}
We note without detail that there is an alternative formula for $\theta$
in terms of perturbed Gaussian integration~\cite{IType}. In that language,
and using also the traffic motifs of Discussion~\ref{com:traffic}, the
three summands in~\eqref{eq:Main} become Feynman diagrams for processes
in which cars $\car_\nu$ governed by parameter $T_\nu=T_1, T_2,$ or $T_3$ interact:

\[ \input{figs/FeynmanDiagrams.pdf_t} \]

In particular, the middle diagram which resembles the Greek letter $\Theta$ gave the invariant its
name.  \endpar{\ref{com:FeynmanDiagrams}}
\end{comment}

\begin{comment} \label{com:Computation} The computation of $G$ is a
bottleneck for the computation of $\Theta$.  It requires inverting
a $(2n+1)\times(2n+1)$ matrix whose entries are (degree 1) Laurent
polynomials in $T$. It's a daunting task yet it takes polynomial
time. Even a naive inversion using Gaussian elimination requires only
$\sim n^3$ operations in the ring $\bbQ(T)$.\footnote{Note that the
degrees of the numerators and denominators occurring in $\bbQ(T)$ are
easily bounded by $2n+1$, and their coefficients cannot grow to have more
than $n\log n$ digits, as they cannot be worse than the worst minor of
$A$. So ring operations in $\bbQ(T)$ in themselves take polynomial time.}
 So $G$ can be
computed in practice even if $n$ is in the hundreds, and everything
which then follows is not worse.

The polynomials $F_1(c)$, $F_2(c_0,c_1)$ and $F_3(k)$ are not
unique, and we are not certain that we have the cleanest possible formulas
for them. They are ugly from a human perspective, yet from a computational perspective,
having 18 terms (as is the case for $F_1(c)$) isn't really a problem;
computers don't care.

Computationally, the worst term in~\eqref{eq:Main}
is the middle one, and even it takes merely $\sim n^2$ operations in
the ring $\bbQ(T_1,T_2)$ to evaluate. \endpar{\ref{com:Computation}}
\end{comment}

\draftcut
\section{Implementation and Examples} \label{sec:Implementation}
\subsection{Implementation} A concise yet reasonably
efficient implementation is worth a thousand formulas. It completely
removes ambiguities, it tests the theories, and it allows for
experimentation. Hence our next task is to implement. The section that
follows was generated from a Mathematica~\cite{Wolfram:Mathematica}
notebook which is available at~\cite[Theta.nb]{Self}. A second implementation of $\Theta$,
using Python and SageMath (\url{https://www.sagemath.org/}) is available at
\url{https://www.rolandvdv.nl/Theta/}.

\input Implementation.tex

\draftcut
\section{Proof of the Main Theorem, Theorem~\ref{thm:Main}} \label{sec:Proof}
We divide the proof into to parts: the invariance of $\theta_0$ (and therefore of $\theta$) is in Section~\ref{ssec:Invariance},
and the polynomiality of $\theta$ is in Section~\ref{ssec:Residue}. 

\subsection{Proof of Invariance} \label{ssec:Invariance} Our proof of the
invariance of $\theta$ (Theorem~\ref{thm:Main}) is very similar,
and uses many of the same pieces, as the proof of the invariance
of $\rho_1$ in~\cite{APAI}. Thus at some places here we are briefer
than at~\cite{APAI}, and sadly, yet in the interest of saving space,
we understate here the interpretation of $g_{\alpha\beta}$ as a
``traffic function''.

Some Reidemeister moves create or lose an edge
and to avoid the need for renumbering it is beneficial to also allow
labelling the edges with non-consecutive labels. Hence we allow that, and
write $\ip$ for the successor of the label $i$ along the knot, and $\ipp$
for the successor of $\ip$ (these are $i+1$ and $i+2$ if the labelling is
by consecutive integers). Also, by convention ``$1$'' will always refer 
to the label of the first edge, and ``$2n+1$'' will always refer to the
label of the last. With this in mind, we have that
$A=I + \sum_c A_c$, with $A_c$ given by
\begin{equation} \label{eq:Ap}
  \begin{array}{c}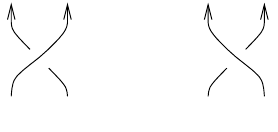\end{array}
  \qquad\longrightarrow\qquad
  \begin{array}{c|cccc}
    A_c &   \text{column }\ip  &  \text{column }\jp \\
    \hline
    \text{row }i & -T^s  & T^s-1 \\
    \text{row }j & 0  & -1
  \end{array}
\end{equation}

Like in~\cite[Lemma~3]{APAI}, the equalities $AG=I$ and $GA=I$ imply
that
for any crossing $c=(s,i,j)$ in a knot diagram $D$, the Green function
$G = (g_{\alpha\beta})$ of $D$ satisfies the following ``$g$-rules'',
with $\delta$ denoting the Kronecker delta:
\begin{equation} \label{eq:CarRules}
  g_{i\beta} = \delta_{i\beta}+T^sg_{\ip,\beta}+(1-T^s)g_{\jp,\beta},
  \qquad g_{j\beta} = \delta_{j\beta}+g_{\jp,\beta},
  \qquad g_{2n+1,\beta} = \delta_{2n+1,\beta},
\end{equation}
\begin{equation} \label{eq:CounterRules}
  g_{\alpha,\ip} = T^s g_{\alpha i} + \delta_{\alpha,\ip},
  \qquad g_{\alpha,\jp} = g_{\alpha j} + (1-T^s)g_{\alpha i} + \delta_{\alpha,\jp},
  \qquad g_{\alpha,1} = \delta_{\alpha,1}.
\end{equation}
Furthermore, the systems of equations \eqref{eq:CarRules} is equivalent
to $AG=I$ and so it fully determines $g_{\alpha\beta}$, and likewise
for the system \eqref{eq:CounterRules}, which is equivalent to $GA=I$.

Of course, the same $g$-rules also hold for $G_\nu = (g_{\nu\alpha\beta})$
for $\nu=1,2,3$, except with $T$ replaced with $T_\nu$.

We also need a variant $\tilg_{ab}$ of $g_{\alpha\beta}$, defined whenever
$a$ and $b$ are two distinct points on the edges of a knot diagram $D$,
away from the crossings. If $\alpha$ is the edge on which $a$ lies and
$\beta$ is the edge on which $b$ lies, $\tilg_{ab}$ is defined as follows:
\begin{equation} \label{eq:tilg}
  \tilg_{ab} = \begin{cases}
    g_{\alpha\beta} &
      \text{if $\alpha\neq\beta$,} \\
    g_{\alpha\beta} &
      \text{if $\alpha=\beta$ and $a<b$ relative to the orientation of the edge $\alpha=\beta$,} \\
    g_{\alpha\beta}-1 &
      \text{if $\alpha=\beta$ and $a>b$ relative to the orientation of the edge $\alpha=\beta$.} \\
  \end{cases}
\end{equation}

Of course, we can define $\tilg_{\nu ab}$ from $g_{\alpha\beta}$ in a similar way.

It is clear that $g$ and $\tilg$ contain the same information and
are easily computable from each other. The variant $\tilg$ is,
strictly speaking, not a matrix and so $g$ is a bit more suitable
for computations. Yet $\tilg$ is a bit better behaved when we try to
track, as below, the changes in $g$ and $\tilg$ under Reidemeister
moves. Reidemeister moves sometimes merge two edges into one or
break an edge into two. In such cases the points $a$ and $b$ can be
``pulled'' along with the move so as to retain their ordering along
the overall parametrization of the knot, yet mere edge labels lose
this information. From the perspective of traffic functions, $\tilg$
is somewhat more natural than $g$, as it makes sense to inject traffic
and to count traffic anywhere along an edge, provided the injection
point and the counting point are distinct.

The following discussion and lemma further exemplify the advantage
of $\tilg$ of $g$:

\parpic[r]{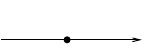}

\begin{discussion} We introduce ``null vertices'' as on the right into
knot diagrams, whose only function (as we shall see) is to cut edges
into parts that may carry different labels. When dealing with upright
knot diagrams as in Figure~\ref{fig:SampleDiagram}, we only allow null
vertices where the tangent to the knot is pointing up, so that the
rotation numbers $\varphi_k$ remain well defined on all edges.  In the
presence of null vertices the matrix $A$ becomes a bit larger (by as many
null vertices as were added to a knot diagram). The rule~\eqref{eq:Ap}
for the creation of the matrix $A$ gets an amendment for null vertices,
\[
  \input{figs/NullVertex.pdf_t}
  \qquad\longrightarrow\qquad
  \begin{array}{c|cccc}
    A_{nv} &   \text{column }k \\
    \hline
    \text{row }j & -1 \\
  \end{array},
\]
and the summation for $A$, $A=I+\sum_{c}A_c+\sum_{nv}A_{nv}$ is extended to include summands for
the null vertices. The matrix $G=A^{-1}$ and the function $g_{\alpha\beta}$ are
defined as before. The $g$-rules of~\eqref{eq:CarRules} and~\eqref{eq:CounterRules} get additions,

\noindent
\begin{minipage}{0.5\linewidth}\begin{equation} \label{eq:NullCarRules}
  g_{j\beta} = \delta_{j\beta} + g_{k\beta},
\end{equation}\end{minipage}
\begin{minipage}{0.5\linewidth}\begin{equation} \label{eq:NullCounterRules}
  \text{and}\qquad g_{\alpha k} = \delta_{\alpha k} + g_{\alpha j},
\end{equation}\end{minipage}

\picskip{0}
\noindent and it remains true that the system of equations
\eqref{eq:CarRules}$\cup$\eqref{eq:NullCarRules} (as well as
\eqref{eq:CounterRules}$\cup$\eqref{eq:NullCounterRules}) fully determines
$g_{\alpha\beta}$.  The variant $\tilg_{ab}$ is also defined as before,
except now $a$ and $b$ need to also be away from the null vertices.
\end{discussion}

\begin{lemma} \label{lem:NullVertices}
Inserting a null vertex does not change $\tilg_{ab}$ provided it is
inserted away from the points $a$ and $b$.\footnote{This statement does
not make sense for $g_{\alpha\beta}$, as inserting a null vertex changes
the dimensions of the matrix $G=(g_{\alpha\beta})$.}
\end{lemma}

\begin{proof} Let $D$ be an upright knot diagram having an edge
labelled $i$ and let $D'$ be obtained from it by adding a null vertex
within edge $i$, naming the two resulting half-edges $j$ and $k$
(in order). Let $g_{\alpha\beta}$ be the Green function for $D$, and
similarly, $g'_{\alpha\beta}$ for $D'$. We claim that
\[ \def\tif{{\text{if }}}
  g'_{\alpha\beta} = \raisebox{7pt}{$
    \raisebox{-7pt}{$
      \left\{\makebox[0pt]{\phantom{$\begin{array}{c}g_{\alpha\beta}\\g_{\alpha\beta}\\g_{\alpha\beta}\end{array}$}}\right.
    $}
    \begin{array}{cccl}
      \tif \beta=j	& \tif \beta=k	& \tif \beta\not\in\{j,k\}	& \\
      g_{ii}		& g_{ii}	& g_{i\beta}			& \tif \alpha=j \\
      g_{ii}-1		& g_{ii}	& g_{i\beta}			& \tif \alpha=k \\
      g_{\alpha i}	& g_{\alpha i}	& g_{\alpha\beta}		& \tif \alpha\not\in\{j,k\} \\
    \end{array}
  $}
\]
Indeed, all we have to do is to verify that the above-defined $g'_{\alpha\beta}$ satisfies all the
$g$-rules \eqref{eq:CarRules}$\cup$\eqref{eq:NullCarRules}, and that is easy. The lemma now follows
easily from the definition of $\tilg'$ in Equation~\eqref{eq:tilg}.\qed
\end{proof}

\begin{remark} \label{rem:F4Null}
The statement of our Main Theorem, Theorem~\ref{thm:Main}, does
not change in the presence of null vertices: There are no ``$F$''
terms for those, and their only effect on the definition of $\Theta$
in Equation~\eqref{eq:Main} is to change the edge labels that appear
within $c$, $c_1$, and $c_2$, and within the $F_3$ sum.
\end{remark}

The following theorem was not named in~\cite{APAI} yet it was stated
there as the first part of the first proof of~\cite[Theorem~1]{APAI}.

\begin{figure}
\[ 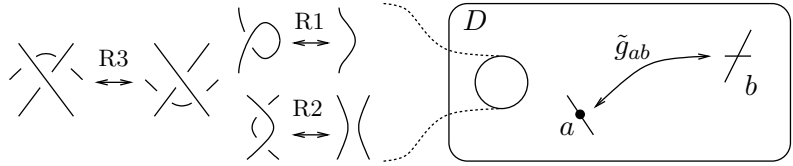 \]
\caption{
  The modified Green function $\tilg_{ab}$ is invariant under Reidemeister moves
  performed away from where it is measured.
} \label{fig:RelativeInvariance}
\end{figure}

\begin{theorem} \label{thm:RelativeInvariant}
The variant Green function $\tilg_{ab}$ is a ``relative invariant'',
meaning that once points $a$ and $b$ are fixed within a knot diagram $D$,
the value of $\tilg_{ab}$ does not change if Reidemeister moves are
performed {\em away} from the points $a$ and $b$ (an illustration appears
in Figure~\ref{fig:RelativeInvariance}). It follows that the same is also
true for $\tilg_{\nu ab}$ for $\nu=1,2,3$.
\end{theorem}

We note that $\tilg_{ab}$ is nearly the same as
$g_{\alpha\beta}$, if $a$ is on $\alpha$ and $b$ is on $\beta$. So
Theorem~\ref{thm:RelativeInvariant} also says that $g_{\alpha\beta}$
is invariant under Reidemeister moves away from $\alpha$ and $\beta$,
except for edge-renumbering issues and $\pm 1$ contributions that arise
if $\alpha$ and $\beta$ correspond to edges that get merged or broken by
the Reidemeister moves.

The proof of Theorem~\ref{thm:RelativeInvariant} is perhaps best
understood in terms of the traffic function of Discussion~\ref{com:traffic}:
One simply needs to verify that for each of the
Reidemeister moves, traffic entering the tangle diagram for the
left hand side of the move exits it in the same manner as traffic
entering the tangle diagram for the right hand side of the move, and
each of these verifications, as explained in~\cite{APAI, Oaxaca-2210,
Toronto-241030}, is very easy. Yet that proof is a bit informal, so
we opt here to give a fully formal proof along the lines of the first
halves of~\cite[Propositions 7-9]{APAI}.

\begin{figure}
\[ 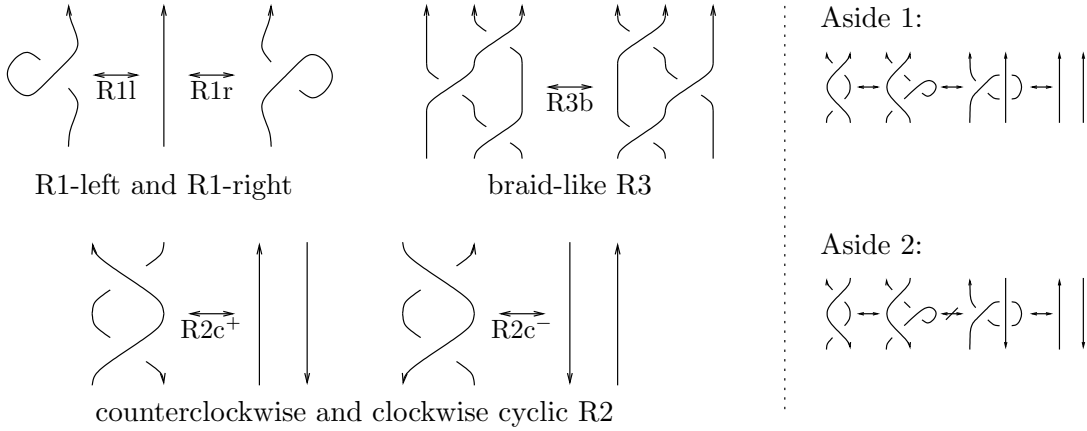 \]
\caption{
  A generating set of oriented Reidemeister moves as in~\cite[Figure~6]{Polyak:RMoves}. Aside 1: the braid-like R2b is not needed.
  Aside 2: yet R2b cannot replace R2c$^\pm$ because in the would-be proof, an unpostulated form of R3 is used (which in
  itself follows from R2c$^\pm$).
} \label{fig:RMoves}
\end{figure}

\vskip 2mm
\noindent{\em Proof of Theorem~\ref{thm:RelativeInvariant}.} We need
to know how the Green function $g_{\alpha\beta}$ changes under the
orientation-sensitive Reidemeister moves of Figure~\ref{fig:RMoves}
(note that the $g_{\alpha\beta}$ do not see the rotation numbers
and don't care if a knot diagram is upright in the sense of
Figure~\ref{fig:SampleDiagram}.

We start with R3b. Below are the two sides of the move, along with
the $g$-rules of type~\eqref{eq:CarRules} corresponding to the
crossings within, written with the assumption that $\beta$ isn't in
$\{\ip,\jp,\kp\}$, so several of the Kronecker deltas can be ignored. We
use $g$ for the Green function at the left-hand side of R3b, and $g'$
for the right-hand side:
\[
  \def\grulesA{{$\begin{array}{l}
    g_{j,\beta} = \delta_{j\beta} \!+\! Tg_{\jp,\beta} \!+\! (1\!-\!T)g_{\kp,\beta} \\
    g_{k,\beta} = \delta_{k\beta} \!+\! g_{\kp,\beta}
  \end{array}$}}
  \def\grulesB{{$\begin{array}{l}
    g_{i,\beta} = \delta_{i\beta} \!+\! Tg_{\ip,\beta} \!+\! (1\!-\!T)g_{\kpp,\beta} \\
    g_{\kp,\beta} = g_{\kpp,\beta}
  \end{array}$}}
  \def\grulesC{{$\begin{array}{l}
    g_{\ip,\beta} = Tg_{\ipp,\beta} \!+\! (1\!-\!T)g_{\jpp,\beta} \\
    g_{\jp,\beta} = g_{\jpp,\beta}
  \end{array}$}}
  \def\gprulesA{{$\begin{array}{l}
    g'_{i,\beta} = \delta_{i\beta} \!+\! Tg'_{\ip,\beta} \!+\! (1\!-\!T)g'_{\jp,\beta} \\
    g'_{j,\beta} = \delta_{j\beta} \!+\! g'_{\jp,\beta}
  \end{array}$}}
  \def\gprulesB{{$\begin{array}{l}
    g'_{\ip,\beta} = Tg'_{\ipp,\beta} \!+\! (1\!-\!T)g'_{\kp,\beta} \\
    g'_{k,\beta} = \delta_{k\beta} \!+\! g'_{\kp,\beta}
  \end{array}$}}
  \def\gprulesC{{$\begin{array}{l}
    g'_{\jp,\beta} = Tg'_{\jpp,\beta} \!+\! (1\!-\!T)g'_{\kpp,\beta} \\
    g'_{\kp,\beta} = g'_{\kpp,\beta}
  \end{array}$}}
   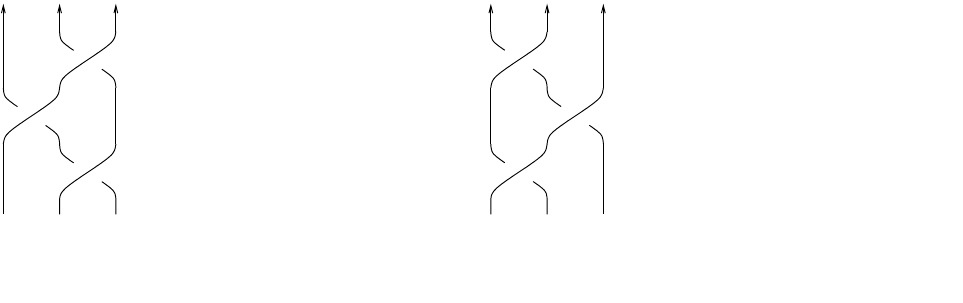
\]

Recall that along with the further $g$-rules and/or $g'$-rules
corresponding to all the non-moving knot crossings, these rules
fully determine $g_{\alpha\beta}$ and $g'_{\alpha\beta}$ for
$\beta\not\in\{\ip,\jp,\kp\}$.

A routine computation (eliminating $g_{\ip,\beta}$, $g_{\jp,\beta}$, and $g_{\kp,\beta}$) shows
that the first system of 6 equations is equivalent to the following system of 6 equations:
\needspace{30mm}
\begin{equation} \label{eq:R3LeftOuter} \begin{gathered}
  g_{i,\beta} = \delta_{i\beta} + T^2g_{\ipp,\beta} + T(1-T)g_{\jpp,\beta} + (1-T)g_{\kpp,\beta}, \\
  g_{j,\beta} = \delta_{j\beta} + Tg_{\jpp,\beta} + (1-T)g_{\kpp,\beta},
  \qquad\qquad g_{k,\beta} = \delta_{k\beta} + g_{\kpp,\beta},
\end{gathered} \end{equation}
\begin{equation} \label{eq:R3LeftInner}
  g_{\ip,\beta} = Tg_{\ipp,\beta} + (1-T)g_{\jpp,\beta},
  \qquad g_{\jp,\beta} = g_{\jpp,\beta},
  \qquad g_{\kp,\beta} = g_{\kpp,\beta}.
\end{equation}
In this system the indices $\ip$, $\jp$ and $\kp$ do not appear in \eqref{eq:R3LeftOuter} or in
the further $g$-rules corresponding to the further crossings. Hence for the purpose of
determining $g_{\alpha\beta}$ with $\alpha,\beta\not\in\{\ip,\jp,\kp\}$,
Equations~\eqref{eq:R3LeftInner} can be ignored.

Similarly eliminating $g'_{\ip,\beta}$, $g'_{\jp,\beta}$, and $g'_{\kp,\beta}$ from the second set
of equations, we find that it is equivalent to
\begin{equation} \label{eq:R3RightOuter} \begin{gathered}
  g'_{i,\beta} = \delta_{i\beta} + T^2g'_{\ipp,\beta} + T(1-T)g'_{\jpp,\beta}
    + (1-T)g'_{\kpp,\beta}, \\
  g'_{j,\beta} = \delta_{j\beta} + Tg'_{\jpp,\beta} + (1-T)g'_{\kpp,\beta},
  \qquad\qquad g'_{k,\beta} = \delta_{k\beta} + g'_{\kpp,\beta},
\end{gathered} \end{equation}
\begin{equation} \label{eq:R3RightInner}
  g'_{\ip,\beta} = Tg'_{\ipp,\beta} \!+\! (1\!-\!T)g'_{\kpp,\beta},
  \qquad g'_{\jp,\beta} = Tg'_{\jpp,\beta} \!+\! (1\!-\!T)g'_{\kpp,\beta},
  \qquad g'_{\kp,\beta} = g'_{\kpp,\beta}.
\end{equation}
Using the same logic as before, for the purpose of 
determining $g'_{\alpha\beta}$ with $\alpha,\beta\not\in\{\ip,\jp,\kp\}$, 
Equations~\eqref{eq:R3RightInner} can be ignored.

But now we compare the unignored equations, \eqref{eq:R3LeftOuter}
and \eqref{eq:R3RightOuter}, and find that they are exactly the same,
except with $g\leftrightarrow g'$, and the same is true for the
further $g$-rules and/or $g'$-rules coming from the further crossings.
Hence so long as $\alpha,\beta\not\in\{\ip,\jp,\kp\}$,
we have that $g_{\alpha\beta}=g'_{\alpha\beta}$. In the case of
the R3b move no edges merge or break up, and hence this implies that
$\tilg_{ab}=\tilg'_{ab}$ so long as $a$ and $b$ are away from the move.

Next we deal with the case of R2c$^+$.  We use the privileges afforded to
us by Lemma~\ref{lem:NullVertices} to insert 4 null vertices into the
right-hand-side of the move, and like in the case of R3b, we start
with pictures annotated with the relevant type~\eqref{eq:CarRules}
and~\eqref{eq:NullCarRules} $g$-rules, written with the assumption that
$\beta\not\in\{\ip,\jp\}$:
\[
  \def\grulesA{{$\begin{array}{l}
    g_{i,\beta} = \delta_{i,\beta} + T^{-1}g_{\ip,\beta}+(1-T^{-1})g_{\jpp,\beta} \\
    g_{\jp,\beta} = g_{\jpp,\beta}
  \end{array}$}}
  \def\grulesB{{$\begin{array}{l}
    g_{\ip,\beta} = Tg_{\ipp,\beta}+(1-T)g_{\jp,\beta} \\
    g_{j,\beta} = \delta_{j,\beta} + g_{\jp,\beta}
  \end{array}$}}
  \def\grulesC{{$\begin{array}{l}
    g'_{i,\beta} = \delta_{i,\beta} + g'_{\ip,\beta} \\
    g'_{\jp,\beta} = g'_{\jpp,\beta}
  \end{array}$}}
  \def\grulesD{{$\begin{array}{l}
    g'_{\ip,\beta} = g'_{\ipp,\beta} \\
    g'_{j,\beta} = \delta_{j,\beta} + g'_{\jp,\beta}
  \end{array}$}}
  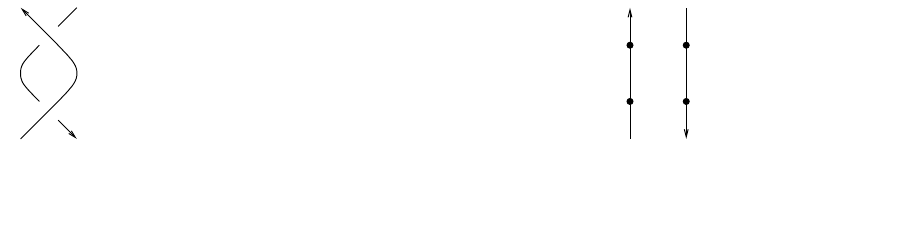
\]

As in the case of R3b, we eliminate $g_{\ip,\beta}$ and $g_{\jp,\beta}$ from the equations for the
left hand side, and find that for the purpose of determining $g_{\alpha\beta}$ with
$\beta\not\in\{\ip,\jp\}$, they are equivalent to the equations
\[
  g_{i,\beta} = \delta_{i,\beta} + g_{\ipp,\beta}
  \qquad\text{and}\qquad g_{j,\beta} = \delta_{j,\beta} + g_{\jpp,\beta}.
\]
Likewise, the right hand side is clearly equivalent to
\[
  g'_{i,\beta} = \delta_{i,\beta} + g'_{\ipp,\beta}
  \qquad\text{and}\qquad g'_{j,\beta} = \delta_{j,\beta} + g'_{\jpp,\beta},
\]
and as in the case of R3b, this establishes the invariance of $\tilg_{ab}$ under R2c moves.

For the remaining moves, R2c$^-$, R1l, and R1r, we merely display the
$g$-rules and leave it to the readers to verify that when the edges $\ip$
and/or $\jp$ are eliminated, the left hand sides become equivalent to
the right hand sides:
\[
  \def\grulesA{{$\begin{array}{l}
    g_{\ip,\beta} = T^{-1}g_{\ipp,\beta}+(1-T^{-1})g_{\jp,\beta} \\
    g_{j,\beta} = \delta_{j,\beta} + g_{\jp,\beta}
  \end{array}$}}
  \def\grulesB{{$\begin{array}{l}
    g_{i,\beta} = \delta_{i,\beta} + Tg_{\ip,\beta}+(1-T)g_{\jpp,\beta} \\
    g_{\jp,\beta} = g_{\jpp,\beta}
  \end{array}$}}
  \def\grulesC{{$\begin{array}{l}
    g'_{i,\beta} = \delta_{i,\beta} + g'_{\ip,\beta} \\
    g'_{\jp,\beta} = g'_{\jpp,\beta}
  \end{array}$}}
  \def\grulesD{{$\begin{array}{l}
    g'_{\ip,\beta} = g'_{\ipp,\beta} \\
    g'_{j,\beta} = \delta_{j,\beta} + g'_{\jp,\beta}
  \end{array}$}}
  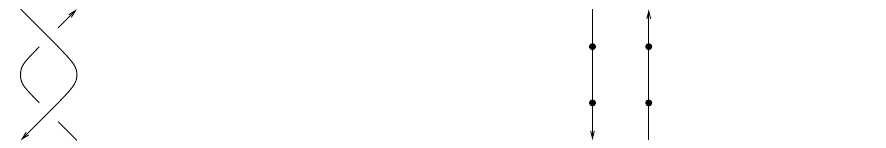
\]
\[
  \def\grulesA{{$\begin{array}{l}
    g_{\ip,\beta} = Tg_{\ipp,\beta} \\
      \hspace{12mm} + (1-T)g_{\ip,\beta} \\
    g_{i,\beta}=\delta_{i,\beta}+g_{\ip,\beta}
  \end{array}$}}
  \def\grulesB{{$\begin{array}{l}
    g'_{\ip,\beta} = g'_{\ipp,\beta} \\
    g'_{i,\beta}=\delta_{i,\beta}+g'_{\ip,\beta}
  \end{array}$}}
  \def\grulesC{{$\begin{array}{l}
    g''_{\ip,\beta} = g''_{\ipp,\beta} \\
    g''_{i,\beta}=\delta_{i,\beta} + Tg''_{\ip,\beta} \\
      \hspace{12mm} + (1-T)g''_{\ipp,\beta}
  \end{array}$}}
  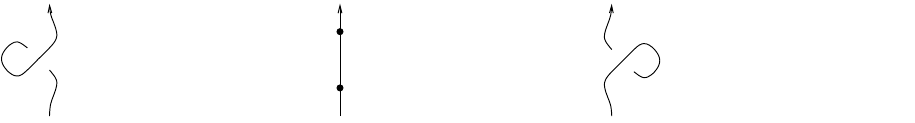
\]
\endpar{\ref{thm:RelativeInvariant}}

\begin{figure}
\[ 
  \def\r#1{$\varphi\!=\!#1$}
  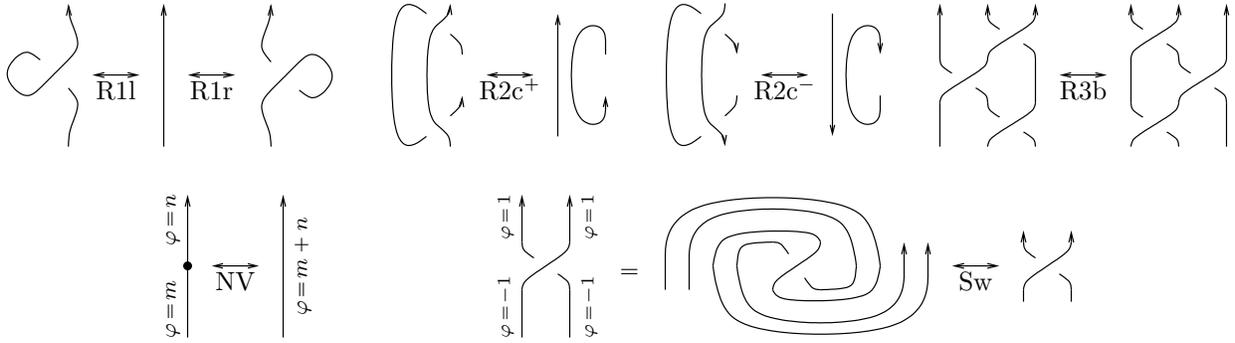
\]
\caption{
  The upright Reidemeister moves: The R1 and R3 moves are already upright
  and remain the same as in Figure~\ref{fig:RMoves}. The crossings in the
  R2 moves of Figure~\ref{fig:RMoves} are rotated to be upright. We also
  need two further moves:
  The null vertex move NV for adding and removing null vertices,
  and the swirl move Sw which then implies that any two ways of turning a crossing upright are the same.
  We sometimes indicate rotation numbers symbolically rather than using complicated spirals.
} \label{fig:UprightRMoves}
\end{figure}

We can now move on to the main part of the proof of our Main Theorem,
Theorem~\ref{thm:Main}. We need to show the invariance of $\theta$ under
the ``upright Reidemeister'' moves of Figure~\ref{fig:UprightRMoves}.

\begin{proposition} \label{prop:UprightRMoves}
The moves in Figure~\ref{fig:UprightRMoves}
are sufficient. If two upright knot diagrams (with null vertices)
represent the same knot, they can be connected by a sequence of moves
as in the figure.
\end{proposition}

\noindent{\em Proof Sketch.} There is an obvious well-defined map
\[
  \frac
    {\text{upright knot diagrams}}
    {\text{relations as in Figure~\ref{fig:UprightRMoves}}}
  \longrightarrow
  \frac
    {\text{oriented knot diagrams}}
    {\text{relations as in Figure~\ref{fig:RMoves}}}
\]
We merely have to construct an inverse to that map. To do that we
have to choose how to turn each crossing in an oriented knot diagram
to be upright. The different ways of doing so differ by instances
of the Sw relation (if deeper spirals need to be swirled away, null
vertices may be inserted using NV and the spirals can be undone
one rotation at a time).  A more detailed version of the proof is
in~\cite{BecerraVanHelden:RotationalReidemeister}. \qed

\begin{proposition} \label{prop:R3} The quantity $\theta_0$ is invariant under R3b.
\end{proposition}

\begin{figure}
\[ {
  \def\i{{$i$}} \def\j{{$j$}} \def\k{{$k$}} \def\m{{$m$}} \def\n{{$n$}} \def\s{{$s$}}
  \def\ip{{$i^+$}} \def\jp{{$j^+$}} \def\kp{{$k^+$}}
  \def\ipp{{$i^{+\!+}$}} \def\jpp{{$j^{+\!+}$}} \def\kpp{{$k^{+\!+}$}}
  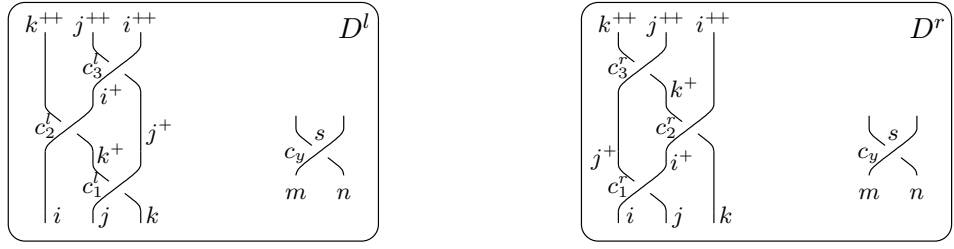
} \]
\caption{
  The two sides $D^l$ and $D^r$ of the R3b move. The left side
  $D^l$ consists of 3 distinguished crossings $c^l_1=(1,j,k)$,
  $c^l_2=(1,i,\kp)$, $c^l_3=(1,\ip,\jp)$ and a collection of further
  crossings $c_y=(s,m,n)\in Y$, where $Y$ is the set of crossings
  not participating in the R3b move. The right side $D^r$ consists of
  $c^r_1=(1,i,j)$, $c^r_2=(1,\ip,k)$, $c^r_3=(1,\jp,\kp)$ and the
  same set $Y$ of further crossings $c_y$.
} \label{fig:R3}
\end{figure}

\noindent{\em Proof.} Let $D_l$ and $D_r$ be two knot diagrams that differ
only by an R3b move, and label their relevant edges and crossings as in
Figure~\ref{fig:R3}. Let $g^l_{\nu\alpha\beta}$ and $g^r_{\nu\alpha\beta}$
be their corresponding Green functions.  Let $F^l_1(c)$, $F^l_2(c_0,c_1)$
and $F^l_3(\varphi,k)$ be defined from $g^l_{\nu\alpha\beta}$ as
in~\eqref{eq:F1}--\eqref{eq:F3}, and similarly make $F^r_1$, $F^r_2$ and $F^r_3$ using
$g^r_{\nu\alpha\beta}$. 

By Theorem~\ref{thm:RelativeInvariant},
$g^l_{\nu\alpha\beta}=g^r_{\nu\alpha\beta}$ so long as
$\alpha,\beta\not\in\{\ip,\jp,\kp\}$. And so the only terms that may
differ in $\theta(D^h)$ between $h=l$ and $h=r$ are the terms
\begin{equation} \label{eq:ABC}
  A^h = \sum_{\makebox[0pt]{$\scriptstyle c\in\{c^h_1,c^h_2,c^h_3\}$}} F^h_1(c)
    \ +\  \sum_{\makebox[0pt]{$\scriptstyle c_0,c_1\in\{c^h_1,c^h_2,c^h_3\}$}} F^h_2(c_0,c_1),
  \quad B^h = 
    \sum_{\makebox[0pt]{$\scriptstyle c_0\in\{c^h_1,c^h_2,c^h_3\},\,c_y\in Y$}} F^h_2(c_0,c_y),
  \quad \text{and} \quad C^h =
    \sum_{\makebox[0pt]{$\scriptstyle c_1\in\{c^h_1,c^h_2,c^h_3\},\,c_y\in Y$}} F^h_2(c_y,c_1).
\end{equation}
We claim that $A^l=A^r$, $B^l=B^r$, and $C^l=C^r$.

To show that $A^l=A^r$, we need to compare polynomials in
$g^l_{\nu\alpha\beta}$ with polynomials in $g^r_{\nu\alpha\beta}$
in which $\alpha$ and $\beta$ may belong to the set $\{\ip,\jp,\kp\}$
on which it may be that $g^l\neq g^r$. Fortunately the $g$-rules of
Equations~\eqref{eq:CarRules} and~\eqref{eq:CounterRules} allow us
to rewrite the offending $g$'s, namely the ones with subscripts in
$\{\ip,\jp,\kp\}$, in terms of other $g$'s whose subscripts are in
$\{i,j,k,\ipp,\jpp,\kpp\}$, where $g^l=g^r$. So it is enough to show that
\begin{equation} \label{eq:R3A}
  \text{under $g^l=g^r$}, \qquad
  A^l\ /.\ \text{(the $g$-rules for $c^l_1$, $c^l_2$, $c^l_3$)}
  = A^r\ /.\ \text{(the $g$-rules for $c^r_1$, $c^r_2$, $c^r_3$)},
\end{equation}
where the symbol $/.$ means ``apply the rules''. This is a finite
computation that can in-principle be carried out by hand. But each $A^h$
is a sum of $3+9=12$ polynomials in the $g^l$'s or the $g^r$'s, these polynomials are
rather unpleasant (see ~\eqref{eq:F1} and~\eqref{eq:F2}), and applying the relevant
$g$-rules adds a bit further to the complexity. Luckily, we can delegate this pages-long
calculation to an entity that works accurately and doesn't complain.

\input Invariance-R3.tex

\begin{proposition} \label{prop:R2c} The quantity $\theta_0$ is invariant under the upright R2c$^+$ and R2c$^-$.
\end{proposition}

\input Invariance-R2c.tex

\parpic[r]{ 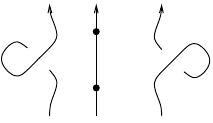 }
\begin{proposition} \label{prop:R1s} The quantity $\theta_0$ is invariant under R1l and R1r.
\end{proposition}

\input Invariance-R1s.tex

\parpic[r]{ \def\r#1{$\varphi\!=\!#1$} 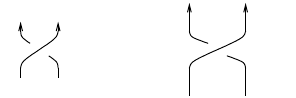 }
\picskip{2}
\begin{proposition} \label{prop:Sw} The quantity $\theta_0$ is invariant under Sw.
\end{proposition}

\input Invariance-Sw.tex

\begin{proposition} \label{prop:NV} The quantity $\theta_0$ is invariant under NV.
\end{proposition}

\noindent{\em Proof.} Indeed, $F_3$ is linear in $\varphi$.
\qed

\vskip 2mm
We are now ready to complete the proof of the first part of the Main Theorem.
\vskip 2mm

\noindent{\em Proof of Invariance.}
The invariance statement in the Main Theorem, Theorem~\ref{thm:Main},
now follows from the invariance of the Alexander polynomial and from
Propositions \ref{prop:UprightRMoves}, \ref{prop:R3}, \ref{prop:R2c},
\ref{prop:R1s}, \ref{prop:Sw}, and~\ref{prop:NV}.  \qed

\draftcut
\subsection{Proof of Polynomiality} \label{ssec:Residue}
We already know (see Comment~\ref{com:normalization}) that the only
obstruction to the polynomiality of $\theta$ comes from the explicit
denominators in Equations~\eqref{eq:F1} and~\eqref{eq:F2}. These
denominators are $(T_2-1)^{-1}$ (if $s,s_1=1$) or $(T_2^{-1}-1)^{-1}
= -T_2(T_2-1)^{-1}$ (if $s,s_1=-1$). So it is enough that we show that
the residue $R$ of $\theta$ at $T_2=1$ vanishes, and this residue comes
solely from the residues of $F_1$ and $F_2$ at $T_2=1$.  Thus $R$
is the knot invariant coming from the same procedure as $\theta$,
only replacing $F_1$, $F_2$, and $F_3$ by their residues $R_1$, $R_2$
and $R_3$ at $T_2=1$. These residues are easily seen to be
\[ R_1(c) = (T^s-1)g_{ji}\left(g_{ii} + 2(T^s-1)g_{ji} - g_{jj}\right), \]
\[ R_2(c_0,c_1) = (T^{s_0}-1)(T^{s_1}-1)g_{j_0i_1}g_{j_1i_0}\left(
    \chi_{i_1\leq i_0} - \chi_{i_1\leq i_0} - \chi_{j_1\leq i_0} + \chi_{j_1\leq j_0}
  \right),
\]
and $R_3=0$, where we have simplified these formulas by making the following observations:
\begin{itemize}
\item $R$ depends only on $T_1$ which we rename to be $T$.
\item At $T_2=1$, $g_{3\alpha\beta}=g_{1\alpha\beta}=g_{\alpha\beta}$.
\item At $T_2=1$, by a simple calculation of the matrices $A$ and $G$
and/or using the traffic interpretation of Comment~\ref{com:traffic},
$g_{2\alpha\beta}$ is the indicator function $\chi_{\alpha\leq\beta}$
of the inequality $\alpha\leq\beta$, which is 1 if the inequality holds
and 0 otherwise.
\end{itemize}

An explicit calculation for some specific knots shows that the sums
corresponding to $R_1$ and to $R_2$ do not vanish individually; instead,
they cancel each other. So we'd better find a technique that relates a
double sum to a single sum. That's the content of the following lemma:

\begin{lemma} \label{lem:del} If there is a function $f(c_0,\gamma)$ that depends
on a crossing $c_0$ and an additional edge label $\gamma$ such that
$(Bf)(c_0)\coloneqq f(c_0,2n+1) - f(c_0,1) = 0$ and such that for any additional
crossing $c_1=(s_1,i_1,j_1)$ we have that
\begin{equation} \label{eq:delf}
  (\partial_{c_1}f)(c_0,c_1) \coloneqq f(c_0,i_1^+) + f(c_0,j_1^+) - f(c_0,i_1) - f(c_0,j_1)
  = R_2(c_0,c_1) + \delta_{c_0,c_1}R_1(c_0),
\end{equation}
then the invariant $R$ vanishes.
\end{lemma}

\begin{proof} Indeed, using the above equation and then telescopic summation over $c_1$ and the vanishing of $Bf$,
\[ R
  = \sum_{c_0,c_1}R_2(c_0,c_1) + \sum_cR_1(c)
  = \sum_{c_0,c_1}(\partial_{c_1}f)(c_0,c_1)
  = \sum_{c_0}(Bf)(c_0) = 0.
\]
\qed
\end{proof}

\vskip 2mm
We can now complete the proof of the second part of the Main Theorem.
\vskip 2mm

\noindent{\em Proof of Polynomiality.}
Take $f(c_0,\gamma)\coloneqq (T^{s_0}-1)g_{\gamma
i_0}g_{j_0^+\gamma}\left(\chi_{\gamma\leq i_0} - \chi_{\gamma\leq
j_0}\right)$.  Use the easily proven facts that $g_{2n+1,i_0} =
0 = g_{j_0^+1}$ to show that $Bf=0$ and then use $g$-rules to verify
Equation~\eqref{eq:delf}. Now using Lemma~\ref{lem:del} we have
that $R=0$ and therefore $\theta$ is a Laurent polynomial. The only
non-integrality for the coefficients of $\theta$ may arise from
the $s/2$ term in Equation~\eqref{eq:F1} and from the $-\varphi_k/2$
terms in Equation~\eqref{eq:F3}. These add up to $(w(D)-\varphi(D)/2$,
using the notation of Equation~\eqref{eq:Delta}. But $w(D)-\varphi(D)$
is always an even number as it is 0 for the long unknot $\uparrow$
and its parity is unchanged by crossing changes and by the moves of
Figure~\ref{fig:UprightRMoves}. \qed

An implementation and a verification of the assertions made in this section is at~\cite[Polynomiality.nb]{Self}.

\draftcut
\section{Strong and Meaningful} \label{sec:SandM}

\subsection{Strong} \label{ssec:Strong}
To illustrate the strength of $\Theta$, Table~\ref{tab:Strong}
summarizes the separation powers of $\Theta$ and of some common knot
invariants and combinations of those knot invariants on prime knots with up
to 15 crossings (up to reflections and reversals).

\begin{table}
\makeatletter
\newcommand{\customlabel}[2]{%
\protected@write \@auxout {}{\string \newlabel {#1}{{#2}{}{}{}{}}}}
\makeatother

\newcounter{tl}
\newcommand{\tl}[1]{%
  \stepcounter{tl}%
  \customlabel{tl:#1}{\thetl}
  \thetl
}
\newcommand\tlref[1]{$^\text{\ref{tl:#1}}$}

\input table.tex

In line~\ref{tl:Ks} of the table we list the total number
of tabulated knots with up to $n$ crossings. For example, there
are 313,230 prime knots up to reflections and reversals with at most 15 crossings. In the
following lines we list the {\em separation deficits} on these knots,
for different invariants or combinations of invariants.  For example,
in line~\ref{tl:D} we can see that on knots with up to 10 crossings, the
Alexander polynomial $\Delta$ has a separation deficit of 38: meaning,
that it attains $249-38=211$ distinct values on the 249 knots with up
to 10 crossings. For deficits, the smaller the better!\footnote{This is
not a political statement.} Thus the deficit of 236,326 for $\Delta$ at
$n\leq 15$ means that the Alexander polynomial is a rather weak invariant,
in as much as separation power is concerned.

In line~\ref{tl:s} we shows the deficits for the
Levine-Tristram signature $\sigma_{LT}$ \cite{Levine:KnotCobordism,
Tristram:CobordismInvariants, Conway:SignaturesSurvey} as computed by the
program in~\cite{Geneva-231201}. We were surprised to find that for knots
with up to 15 crossings these deficits are smaller than those of $\Delta$.

Line~\ref{tl:J} shows the deficits for the Jones polynomial $J$. It is
better than $\Delta$, and better than $\Delta$ and $\sigma_{LT}$ taken together
(deficits not shown) but still rather weak. Line~\ref{tl:Kh} shows
the deficits for Khovanov homology $\Kh$. They are only a bit lower
than those of $J$. On line~\ref{tl:H}, the HOMFLY-PT polynomial $H$ is
noticeably better.

On line~\ref{tl:V} we consider the hyperbolic
volume $\Vol$ of the knot complement, as computed by
SnapPy~\cite{CullerDunfieldGoernerWeeks:SnapPy}. We computed volumes
using SnapPy's \verb$high_precision$ flag, which makes SnapPy compute
to roughly 63 decimal digits, and then truncated the results to 58
decimal digits to account for possible round-off errors within the last
few digits. But then we are unsure if we computed enough\ldots. Hence
the uncertainty symbols ``$\sim$'' on some of the results here and in
the other lines that contain $\Vol$. This said, $\Vol$ seems to be the champion so far.

Line~\ref{tl:KhHV} is ``everything so far, taken together''.  Note that
$\Kh$ dominates $J$ and $H$ dominates both $\Delta$ and $J$, so there's
no point adding $\Delta$ and/or $J$ into the mix. We note that adding
$\sigma_{LT}$ to the triple $(\Kh,H,\Vol)$, or even to the pair $(\Kh,\Vol)$,
does not improve the results; namely, for knots with up to 15 crossings
the pair $(\Kh,\Vol)$ dominates $\sigma_{LT}$, even though each of $\Kh$ and
$\Vol$ does not dominate $\sigma_{LT}$ and the discrepancies start already
at 11 crossings. We don't know if this means anything.

On line~\ref{tl:r1}, the Rozansky-Overbay invariant $\rho_1$
\cite{Rozansky:Contribution, Rozansky:Burau, Rozansky:U1RCC,
Overbay:Thesis}, also discussed by us in~\cite{APAI}, does somewhat
better. Note that the computation of $\Delta$ is a part of the computation
of $\rho_1$, so we always take them together. In line~\ref{tl:r12} we add
$\rho_2$~\cite{Oaxaca-2210} to make the results yet a bit better.

Line~\ref{tl:r12KhHV} is ``everything before $\Theta$''.

Line~\ref{tl:Th} makes our case that $\Theta$ is strong --- the deficit
here, for knots with up to 15 crossings, is about a sixth of the deficit
in line~\ref{tl:r12KhHV}! For the interested, Figure~\ref{fig:DrieParen}
shows the 3 pairs that create the deficit in the column $n\leq 11$
of this line.

\begin{figure}
\[
  \includegraphics[height=2in]{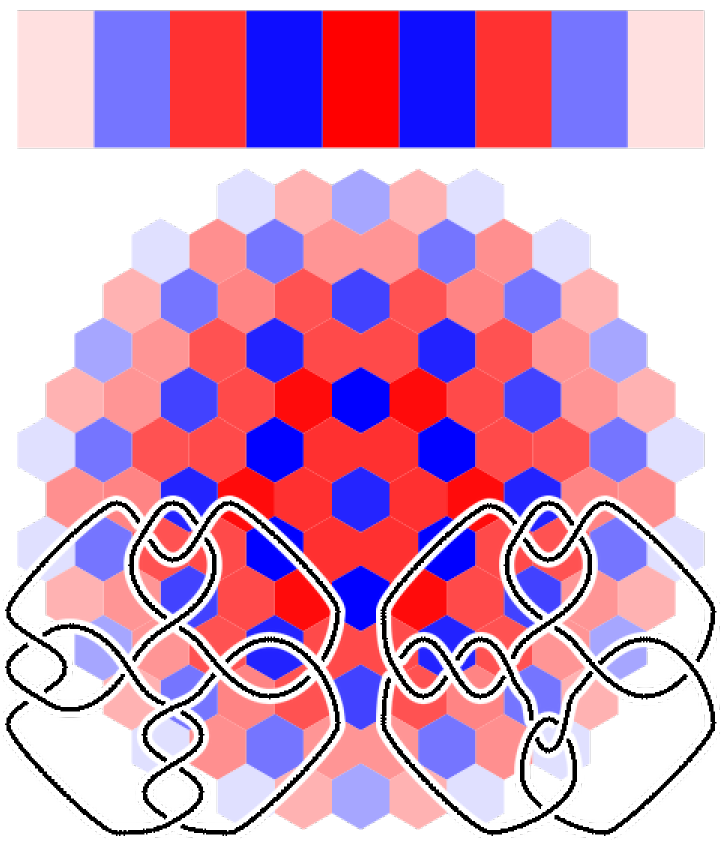}\quad
  \includegraphics[height=2in]{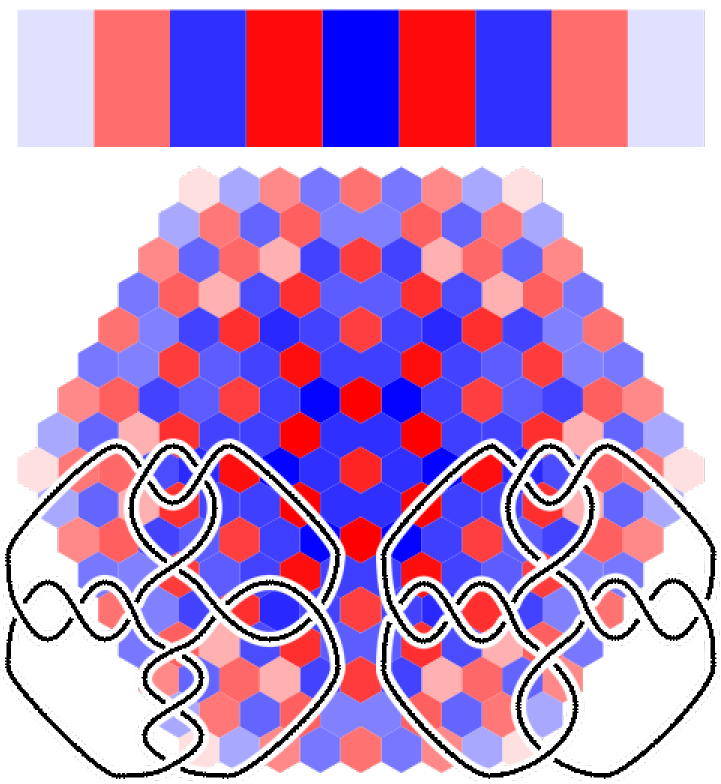}\quad
  \includegraphics[height=2in]{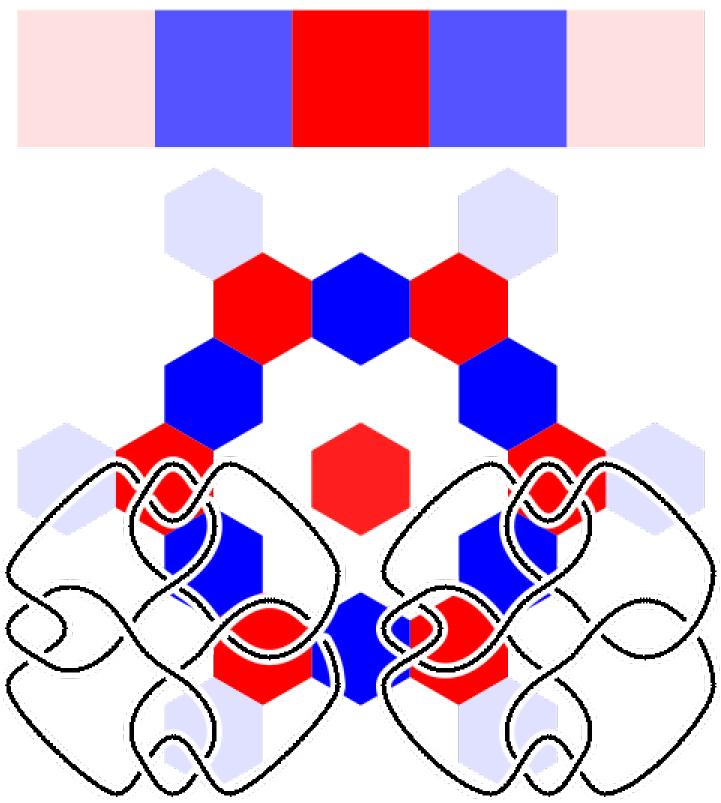}
\]
\caption{
  The three pairs responsible for the deficit of 3 in the column
  $n\leq 11$ of line~\ref{tl:Th} of Table~\ref{tab:Strong}. They
  are $(11_{a44},11_{a47})$, $(11_{a57},11_{a231})$, and
  $(11_{n73},11_{n74})$, and each pair is a pair of mutant Montesinos knots
  (though $\Theta$ sometimes does separate mutant pairs, as was shown
  in Section~\ref{ssec:Examples}).
} \label{fig:DrieParen}
\end{figure}

Line~\ref{tl:Thr2} reinforces our case by just a bit: note
that it makes sense to bundle $\rho_2$ along with $\Theta$, for their
computations are very similar. Note also that Theorem~\ref{thm:rho1}
below means that it is pointless to consider $(\Theta,\rho_1)$.

Line~\ref{tl:Ths} shows that for knots with up to 15 crossings, $\Theta$ dominates
$\sigma_{LT}$. We don't know if this persists.

Lines~\ref{tl:ThKh} through~\ref{tl:ThV} show that at crossing number $\leq
15$ and in the presence of $\Theta$, and especially in the presence of
both $\Theta$ and $\rho_2$, it is pointless to also consider $H$ or $\Kh$,
and only mildly useful to also consider $\Vol$. Line~\ref{tl:Thr2KhHV}
shows that once $\Vol$ has been added to $\Theta$, the other invariants
contribute almost nothing.

We note that of all the invariants considered above, the
only one known to (sometimes) detect knot mutation is $\Theta$ (see
Section~\ref{ssec:Examples}).

We also note that the $V_n$ polynomials of Garoufalidis and
Kashaev~\cite{GaroufalidisKashaev:Multivariable}, and in particular
$V_2$~\cite{GaroufalidisLi:Patterns} share many properties with $\Theta$
and are stronger than $\Theta$ on knots with up to 15 crossings. But
they are not nearly as computable on large knots. It would be very interesting to explore the
relationship between the $V_n$'s and $\Theta$.

\subsection{Meaningful} \label{ssec:Meaningful}
Many knot polynomials have some separation power, some more and some
less, yet they seem to ``see'' almost no other topological properties of
knots. The greatest exception is the Alexander polynomial, which despite
having rather weak separation powers, gives a genus bound, a fiberedness
condition, and a ribbon condition. The definition of $\theta$ is in some
sense ``near'' the definition of $\Delta$, and one may hope that $\theta$
will share some of the good topological properties of $\Delta$.

\subsubsection{The Knot Genus} The following theorem is proven in~\cite{Genus}:

\begin{theorem} \label{thm:Genus}
Let $K$ be a knot and $g(K)$ the genus of $K$. Then $\deg_{T_1}\theta(K)\leq 2g(K)$.
\end{theorem}

The example of the Conway knot
and the Kinoshita-Terasaka knot in Section~\ref{ssec:Examples} shows that
the bound in Theorem~\ref{thm:Genus} can be stronger than the bound
$\deg_T\Delta(K)\leq g(K)$ coming from the Alexander polynomial. Another
such example is the 48-crossing Gompf-Scharlemann-Thompson
$\mathit{GST}_{48}$ knot \cite{GompfScharlemannThompson:Counterexample}
of Figure~\ref{fig:GST48}.  Here's the relevant computation, with $X_{14,1}$
(say) meaning ``the crossing $(1,14,1)$'' and $\bar{X}_{2,29}$ (say)
meaning ``$(-1,2,29)$'':

\begin{figure}
\[ \resizebox{0.75\linewidth}{!}{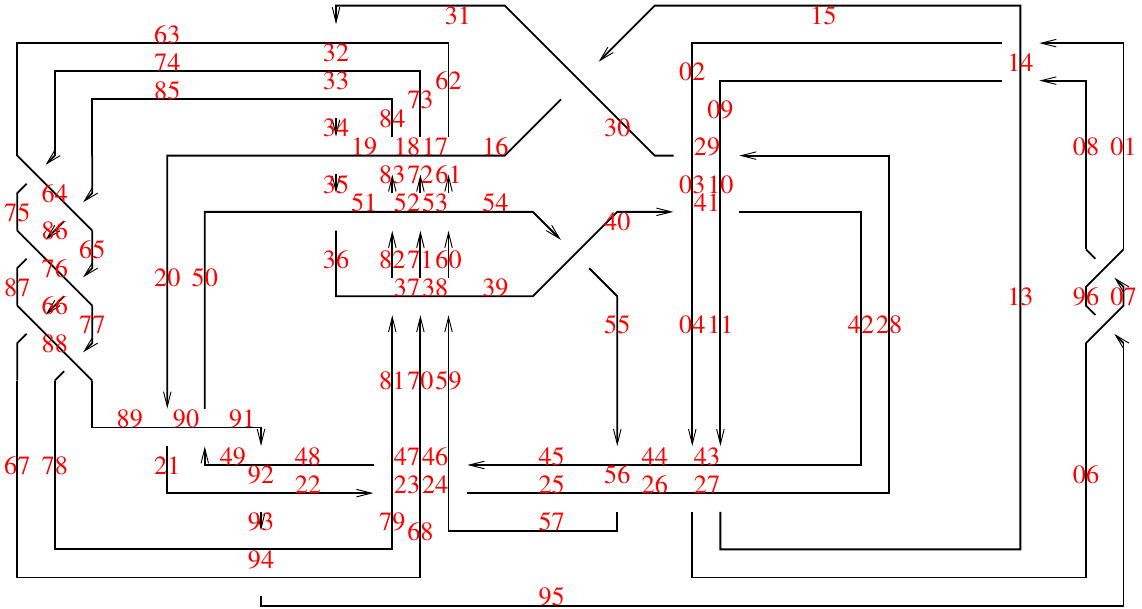} \]
\caption{The 48-crossing Gompf-Scharlemann-Thompson
  $\mathit{GST}_{48}$ knot \cite{GompfScharlemannThompson:Counterexample}.
} \label{fig:GST48}
\end{figure}

\input{GST48.tex}

Thus $\theta$ gives a better lower bound on the genus of
$\mathit{GST}_{48}$, 10, then the lower bound coming from $\Delta$, which
is 8. Seeing that $\mathit{GST}_{48}$ may be a counter-example to the
ribbon-slice conjecture~\cite{GompfScharlemannThompson:Counterexample},
we are happy to have learned more about it. Also see Dream~\ref{dream:Ribbon} below.

The hexagonal QR code of large knots is often a clear hexagon
(e.g. Figure~\ref{fig:300}), but the hexagonal QR code of
$\mathit{GST}_{48}$, displayed above, is rounded at the corners. We
don't know if this is telling us anything about topological properties of
$\mathit{GST}_{48}$.

\subsubsection{Fibered Knots} Upon inspecting the values of $\Theta$
on the Rolfsen table, Figure~\ref{fig:Rolfsen}, we noticed that often
(but not always) the bar code shows the exact same colour sequence
as the top row of the QR code, or exactly its opposite. This and some
experimentation lead us to the following conjecture, for which we do
not have theoretical support. See a similar result on the ADO invariant
at~\cite{LopezNeumannVanDerVeen:FibredLinks}.

\begin{conjecture} \label{conj:Fibered}
If $K$ is a fibered knot and $d$ is the degree of $\Delta(K)$ (the highest power of
$T$), then the coefficient of $T_2^{2d}$ in $\theta(K)$, which is a
polynomial in $T_1$, is an integer multiple of $T_1^d\Delta(K)|_{T\to T_1}$. See examples in
Figure~\ref{fig:FiberedExamples}, where the integer factor is denoted $s(K)$.
\end{conjecture}

Using the available fiberedness data in KnotInfo~\cite{LivingstonMoore:KnotInfo} we found
that the condition in this conjecture holds for all 5,397 fibered knots with
up to 13 crossings, while it fails on all but 48 of the 7,568 non-fibered knots
with up to 13 crossings. See~\cite[FiberedKnots.nb]{Self}.

We note that if $K$ is fibered then degree $d$ of $\Delta(K)$
is the genus of $K$, and $\Delta(K)$ is monic, meaning
that the coefficient of $T^d$ in $\Delta(K)$ is $\pm 1$
(see~\cite[Section~10H]{Rolfsen:KnotsAndLinks}). The latter condition
is an often-used fast-to-compute criterion for a knot to be fibered.

\begin{figure}
  \begin{minipage}[c]{0.49\linewidth}
    \captionsetup{width=0.9\linewidth}
    \caption{
      The invariant $\Theta$ of the fibered knot $12_{n242}$, also known
      as the $(-2,3,7)$ pretzel knot, and of the fibered knot $7_7$. For the first, $s(K)>0$
      and the bar code visibly matches with the top row of the QR code
      (though our screens and printers and eyes may not be good enough
      to detect minor shading differences, so a visual inspection may
      not be enough). For the second, twice the degree of $\Delta$ is visibly greater than the degree of
      $\theta$, so $s(K)=0$.
    } \label{fig:FiberedExamples}
  \end{minipage}
  \hfill
  \begin{minipage}[c]{0.49\linewidth}
    \includegraphics[height=1.7in]{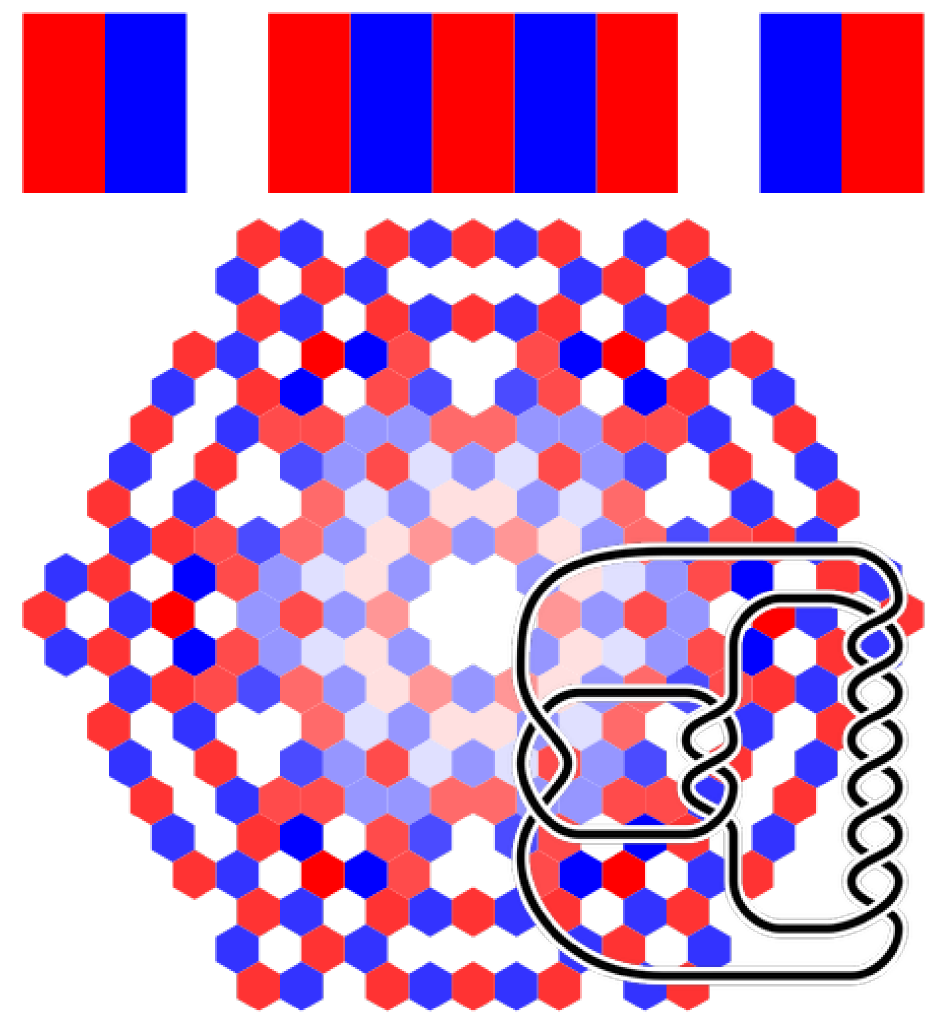}
    \hfill
    \includegraphics[height=1.7in]{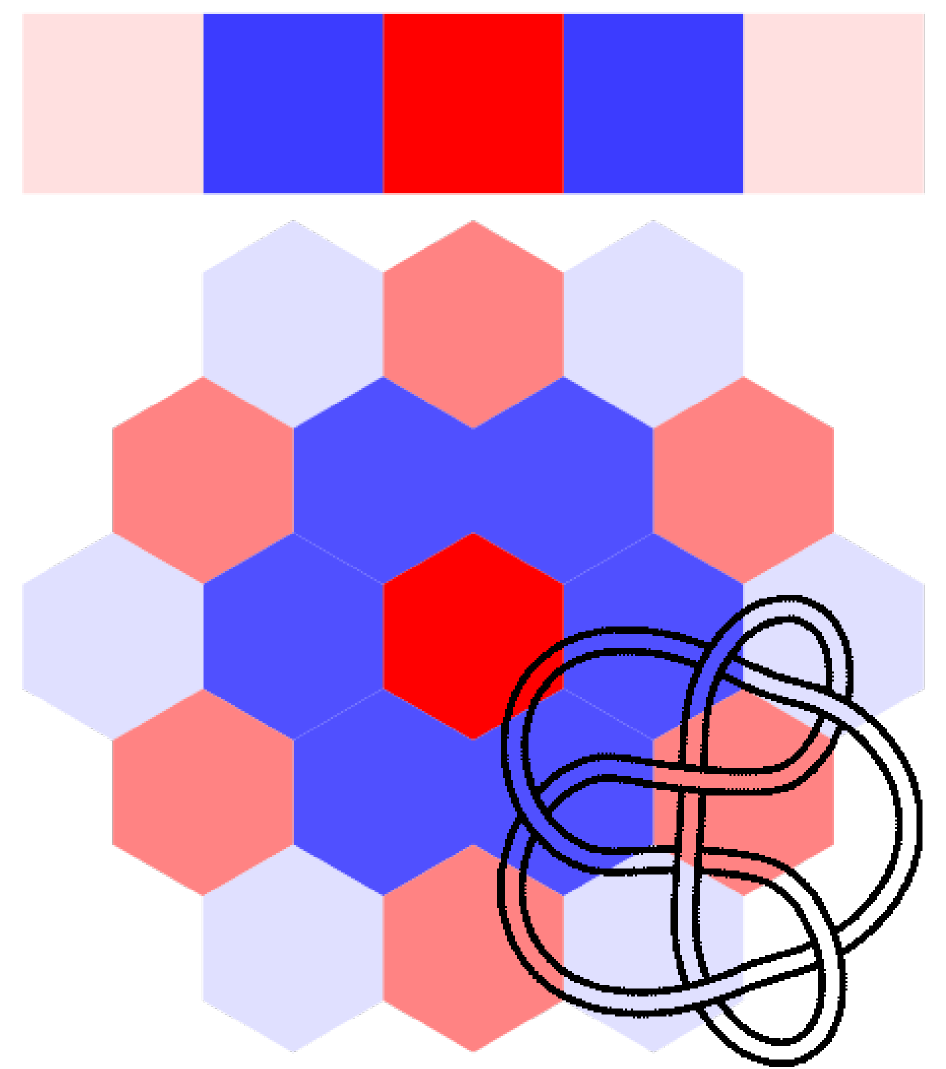}
  \end{minipage}
\end{figure}

If Conjecture~\ref{conj:Fibered} is true then the condition in it is another
fast-to-compute criterion for a knot to be fibered, and this criterion
is sometimes stronger than the Alexander condition. For example, both the Conway and the
Kinoshita-Terasaka knots are not fibered yet their Alexander polynomial is $1$, which is monic. In
both cases the coefficient of $T_2^0$ in $\theta$ is not an integer multiple of $1$ (see
Section~\ref{ssec:Examples}), so the condition
in Conjecture~\ref{conj:Fibered} would detect that these two knots are not fibered.

\draftcut
\section{Stories, Conjectures, and Dreams} \label{sec:CandD}

There is a storyteller in each of us, who wants to tell a coherent story,
with a beginning, a middle, and an end. Unfortunately of us, the $\Theta$
story isn't that neat. Calling the content of the first few sections of
this paper ``the middle'', we are quite unsure about the beginning and
the end. The ``beginning'' can be construed to mean ``the thought process
that lead us here''. But that process was too long and roundabout to
be given in full here (though much of it can be gleaned by reading this
section). What's worse, we believe that ultimately, our peculiar thought
process will be replaced by much more solid foundations and motivations,
perhaps along the lines of Dreams~\ref{dream:Seifert}
and~\ref{dream:UniversalSeifert}. But this solid foundation is not
available yet, even if we are working hard to expose it. As for the end
of the story, it is clearly in the future.

Hence this section is a bit sketchy and disorganized. Those facts that we
already know, those conjectures we believe in, and the dreams we dream,
are here in some random order. But the narrative is lacking.

Many of the statements below continue a theme from
Section~\ref{ssec:Meaningful}, that $\theta$ shares many of the properties
of $\Delta$, and sometimes sharpens them.

\begin{conjecture} \label{conj:D6}
$\theta$ has hexagonal symmetry. That is, for any knot $K$, $\theta(K)$ is
invariant under the substitutions $(T_1\to T_1,T_2\to T_1^{-1}T_2^{-1})$
(``the QR code is invariant under reflection about a horizontal line''),
and $(T_1\to T_1T_2, T_2\to T_2^{-1})$ (``the QR code is invariant under
reflection about the line of slope $30^\circ$'').
\end{conjecture}

The Alexander polynomial $\Delta$ is invariant under a simpler symmetry,
$T\to T^{-1}$. It is rather difficult to deduce the symmetry
of $\Delta$ from the formula in this paper, Equation~\eqref{eq:Delta}
(though it is possible; once notational differences are overcome, the
proof is e.g.\ in~\cite[Chapter~IX]{CrowellFox:KnotTheory}). Instead,
the standard proof of the symmetry of $\Delta$ uses the Seifert surface
formula for $\Delta$ (e.g.~\cite[Chapter~6]{Lickorish:KnotTheory}). We
expect that Conjecture~\ref{conj:D6} will be proven as soon as a Seifert
formula is found for $\theta$. See Dream~\ref{dream:Seifert} below.

\begin{figure}
  \begin{minipage}[c]{0.6\linewidth}
    \captionsetup{width=0.9\linewidth}
    \caption{
      A long version of the rotational virtual knot $KS$ from~\cite{Kauffman:RotationalVirtualKnots}.
      It has $X=\{(-1,1,6),(-1,2,4),(1,9,3),(-1,7,5),(1,10,8)\}$ and
      $\varphi=(-1,0,0,1,0,-1,0,0,1,0,0)$.
    } \label{fig:KS}
  \end{minipage}
  \begin{minipage}[c]{0.35\linewidth}\[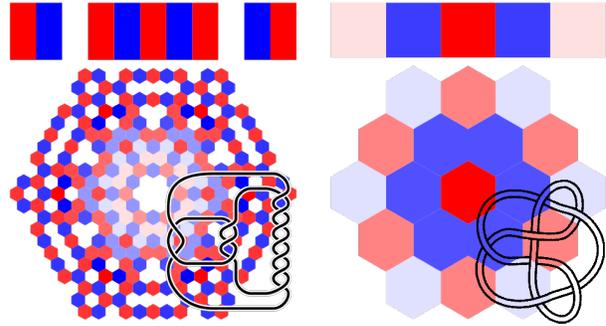\]\end{minipage}
\end{figure}

A {\em rotational virtual knot} is a virtual knot
diagram~\cite{Kauffman:VirtualKnotTheory} whose edges\footnote{Ignoring
``virtual crossings''. See~\cite[Section~4]{BDV:OU}.} are marked
with ``rotation numbers'' $\varphi_k$, modulo the same moves as in
Figure~\ref{fig:UprightRMoves}.\footnote{This definition is slightly
different than the original in~\cite{Kauffman:RotationalVirtualKnots}
but the equivalence is easy to show.} Clearly, $\Theta$ extends to
long rotational virtual knots, and the proof of the Main Theorem,
Theorem~\ref{thm:Main}, extends nearly verbatim\footnote{The only
exception is that some of the coefficients of $\theta$ may be half
integers, as $w(D)-\varphi(D)$ may be odd for a rotational virtual knot
diagram.}. Yet as shown below, on the long rotational virtual knot $KS$
of Figure~\ref{fig:KS} (and indeed, on almost any other long rotational
virtual knot which is not a classical knot), the hexagonal symmetry of
$\theta$ fails. So something non-local must happen within any proof of
Conjecture~\ref{conj:D6}.

\input KS.tex

\begin{conjecture} \label{conj:Mirror} If $\bar{K}$ denotes the mirror image of a knot $K$, then
$\theta(\bar{K}) = -\theta(K)$.
\end{conjecture}

\begin{conjecture} \label{conj:Reverse} If $-K$ denotes the reverse
of a knot $K$ (namely, $K$ taken with the opposite orientation), then
$\theta(-K)=\theta(K)$.
\end{conjecture}

\begin{fact} \label{fact:ConnectedSum}
$\theta_0(K)$ is additive under the connected sum operation of knots:
$\theta_0(K_l\#K_r)=\theta_0(K_l)+\theta_0(K_r)$. Equivalently, using the known multiplicativity of $\Delta$,
\[ \theta(K_l\#K_r)
    = \theta(K_l)\Delta_1(K_r)\Delta_2(K_r)\Delta_3(K_r)
    + \theta(K_r)\Delta_1(K_l)\Delta_2(K_l)\Delta_3(K_l).
\]
\end{fact}

Oddly, Fact~\ref{fact:ConnectedSum} is easier to prove than Conjectures~\ref{conj:Mirror}
and~\ref{conj:Reverse}:

\noindent{\em Proof Sketch.} The $F_1$ and $F_3$ summations
in~Equation~\eqref{eq:Main} are clearly additive, and so is the part
of the $F_2$ summation in which $c_0$ and $c_1$ fall within the same
component. It remains to consider the case where $c_0$ and $c_1$ fall
within different components. But in that case, the factor $g_{1j_1i_0}
g_{3j_0i_1}$ within the definition of $F_2$ in~\eqref{eq:F2} vanishes
because cars only drive forward, and either $g_{1j_1i_0}$ or $g_{3j_0i_1}$
measures traffic going backwards. \qed

\begin{theorem}[Murugesan,~\cite{Murugesan:ThetaRecovers}]
\label{thm:rho1} $\theta$ dominates the Rozansky-Overbay invariant
$\rho_1$ \cite{Rozansky:Contribution, Rozansky:Burau, Rozansky:U1RCC,
Overbay:Thesis}, also discussed by us in~\cite{APAI}. In fact, $\rho_1 =
-\theta|_{T_1\to 1, T_2\to T}$.
\end{theorem}

\begin{conjecture} \label{conj:TwoLoop} $\theta$ is equal
to the ``two-loop polynomial'' studied extensively by Ohtsuki
\cite{Ohtsuki:TwoLoop}, continuing Rozansky, Garoufalidis, and Kricker
\cite{GaroufalidisRozansky:LoopExpansion, Rozansky:Contribution,
Rozansky:Burau, Rozansky:U1RCC, Kricker:Lines}.
\end{conjecture}

\begin{discussion} \label{disc:TwoLoop}
People who are already familiar with ``the loop expansion'' may consider
the above conjecture an ``explanation'' of $\theta$. We differ. An
elementary construction ought to have a simple explanation, and the loop
expansion is too complicated to be that.

Be it as it may, Ohtsuki~\cite{Ohtsuki:TwoLoop}
shows that Conjecture~\ref{conj:TwoLoop} implies
Conjectures \ref{conj:D6}, \ref{conj:Mirror},
and~\ref{conj:Reverse} as well as Theorem~\ref{thm:Genus} and
Fact~\ref{fact:ConnectedSum}. Conjecture~\ref{conj:TwoLoop} would
also predict the behaviour of $\theta$ under Whitehead doubles as
in~\cite{Garoufalidis:Whitehead} and under cabling operations as
in~\cite{Ohtsuki:Cabling}.
\endpar{\ref{disc:TwoLoop}}
\end{discussion}

Next, let us briefly sketch some key points from~\cite{DPG, PG}, where we
explain how to obtain poly-time computable knot invariants from certain
Lie algebraic constructions.

\begin{discussion} \label{disc:SolvApp}
Let $\frakg$ be a semi-simple Lie algebra, let $\frakb$ be its upper Borel
subalgebra, and let $\frakh$ be its Cartan subalgebra. Then $\frakb$ has
a Lie bracket $\beta$ and, as the dual of the lower Borel subalgebra,
it also has a cobracket $\delta$. It turns out that $\frakg$ can be
recovered from the triple $(\frakb,\beta,\delta)$; in fact, $\frakg^+
\coloneqq \frakg\oplus\frakh \simeq \calD(\frakb,\beta,\delta)$, where
$\calD$ denotes the Manin double construction\footnotemark. We now
set $\frakg^+_\eps \coloneqq \calD(\frakb,\beta,\eps\delta)$, where
$\eps$ is a formal ``small'' parameter. The family $\frakg^+_\eps$
is a 1-parameter family of Lie algebras all defined on the same
underlying vector space $\frakb\oplus\frakb^*$. If $\eps$ is invertible
then $\frakg^+_\eps$ is independent of $\eps$ and is always isomorphic to
$\frakg^+=\frakg^+_1$. Yet at $\eps=0$, $\frakg^+_0$ is solvable, and as
the name ``solvable'' suggests, computations in $\frakg^+_0$ can be ``solved'',
meaning, can be carried out efficiently in closed form.

\footnotetext{We are unsure about naming. $\calD$ is also known as
``the Drinfeld double'' construction for Lie bialgebras (as opposed to
Hopf algebras). Yet when Drinfeld first refers to this construction
in~\cite{Drinfeld:QuantumGroups}, in reference to Lie bialgebras,
he repeatedly names it after Manin (under the less clear name ``Manin
triples''), yet without providing a reference. Our choice is to use
``Manin double'' when doubling Lie bialgebras and ``Drinfeld double''
when doubling a Hopf algebra, as we found no indication that Manin knew
about the latter process.}

Hence in~\cite{DPG, PG}, mostly in the case where $\frakg=sl_2$,
we use standard techniques to quantize the universal
enveloping algebra $\calU(\frakg^+_\eps)$ and use it to define
a ``universal quantum invariant'' $Z^\frakg_\eps$ (in the sense
of~\cite{Lawrence:UniversalUsingQG, Ohtsuki:QuantumInvariants}). We then
expand $Z^\frakg_\eps$ near where it's easy; namely, as a power series
around $\eps=0$. In the case of $\frakg=sl_2$, and almost certainly in
general, we write $Z^\frakg_\eps = \rho_0^\frakg\exp\left(\sum_{d\geq
1}\rho^\frakg_d\eps^d\right)$ and find that we can interpret the
$\rho^\frakg_d$ as polynomials in as many variables as the rank of
$\frakg$. It turns out that $\rho^\frakg_0$ is always determined by the
Alexander polynomial and the $\rho^\frakg_d$ are always computable in
polynomial time (with polynomials whose exponents and coefficients get
worse as $d$ grows bigger and $\frakg$ gets more complicated).

Our papers and talks \cite{APAI, PG, Oaxaca-2210} carry out the above
procedure in the case where $\frakg=sl_2$, calling the resulting
invariants $\rho_d$, for $d\geq 1$. They are the same as $\rho_1$ and
$\rho_2$ of Section~\ref{ssec:Strong}. \endpar{\ref{disc:SolvApp}}

\end{discussion}

Following some preliminary work by Schaveling~\cite{Schaveling:Thesis},
in the summer of 2024 we've set out to find good formulas for
$\rho_1^{sl_3}$. Tracing Discussion~\ref{disc:SolvApp} seemed technically
hard, so instead, we extracted from the procedure the ``shape'' of the
formulas we could expect to get and, and then we found the invariant
$\theta$ by the method of undetermined coefficients assisted by some
difficult-to-formulate intuition (more in Comment~\ref{com:UC}
below). Thus our formulas for $\theta$ arose from our expectations for
$\rho_1^{sl_3}$, and yet we have not proved that they are equal!

\begin{conjecture} \label{conj:sl3} Up to conventions and normalizations, $\theta=\rho_1^{sl_3}$.
\end{conjecture}

\begin{comment} $\theta$ and the invariants cosidered by Harper~\cite{Harper:Non-Abelian} share
their relationship with $sl_3$ and with the Alexander polynomial, but there does not seem to be an
immediate relationship between them.
\end{comment}

\begin{discussion} \label{disc:LieVersed}
People who are versed with Lie algebras and their quantizations may
consider the above an ``explanation'' of $\theta$, and may be looking
forward to a more detailed exposition of $\rho^\frakg_d$. We differ,
for the same reasons as in Discussion~\ref{disc:TwoLoop}. We expect the
eventual ``origin story'' of $\theta$ to be simpler and more
natural.\endpar{\ref{disc:LieVersed}}
\end{discussion}

\begin{discussion} \label{disc:Satellites} Seeing that the coproduct of
the quantized algebras of Discussion~\ref{disc:SolvApp} correspond to strand
doubling, and also noting Ohtsuki's~\cite{Ohtsuki:Cabling}, we expect that
there should be cabling and satellite formulas for all the invariants of
the type $\rho^\frakg_d$, and in particular for $\Theta$. In particular, it should not be 
possible to increase the separation power of $\Theta$ by 
pre-composing it with cabling or satellite operations. \endpar{\ref{disc:Satellites}}
\end{discussion}

\begin{discussion} \label{disc:FPGI}
It is the basis of the theory of ``Feynman diagrams'',
and hence it is extremely well known in the physics community, that
perturbed Gaussian integrals, when convergent,  can be computed
(as asymptotic series) efficiently using ``Feynman diagrams'' (see
e.g.~\cite{Polyak:FeynmanDiagrams}). Physicists use this routinely
in infinite dimensions; yet the finite dimensional formulation can be
sketched as follows:
\begin{equation} \label{eq:PGI}
  \int_{\bbR^d}e^{Q+\eps P} \sim C\sum_{n\geq 0}\eps^n\sum_F\calE(F),
\end{equation}
where $Q$ is a non-degenerate quadratic on $\bbR^d$, $P$ is a ``smaller''
perturbation, $C$ is some constant involving $\pi$'s and the determinant
of $Q$, the summation $\sum_F$ is over ``Feynman diagrams'' of complexity
$n$, and $F\mapsto\calE(F)$ is some procedure, which can be specified in
full but we will not do it here, which assigns to every Feynman diagram
$F$ an algebraic sum which in itself depends only on the coefficients
of $P$ and the entries of the inverse of $Q$.

In fact, one may take the right-hand-side of Equation~\eqref{eq:PGI}
to be the definition of the left-hand-side, especially if
the left-hand-side is not convergent, or does not make sense for some other reason. Namely, one may
set
\begin{equation} \label{eq:FPGI}
  \Gint_{\bbR^d}e^{Q+\eps P} \coloneqq C\sum_{n\geq 0}\eps^n\sum_F\calE(F).
\end{equation}
The result is an integration theory defined on perturbed Gaussians in
fully algebraic terms, and which shares some of the properties of ``ordinary'' integration, such as
having a version of Fubini's theorem. In a sense, that's what physicists do: path integrals don't
quite make sense, so instead they are defined using Feynman diagrams and the right-hand-side of
Equation~\eqref{eq:FPGI}. Another example is the ``\AA{}rhus integral'' of
\cite{Bar-NatanGaroufalidisRozanskyThurston:Aarhus}, where the integral in itself is diagrammatic, as
is the output of the integration procedure.\endpar{\ref{disc:FPGI}}
\end{discussion}

\begin{fact} \label{fact:PerturbedGaussian}
There is a perturbed Gaussian formula for $\Theta$. More
precisely, one can assign a 6-dimensional Euclidean space $\bbR^6_e$ with coordinates
$p_{1e},p_{2e},p_{3e},x_{1e},x_{2e},x_{3e}$ to each edge $e$ of a knot diagram $D$ and then form
$\bbR_{6E}\coloneqq\prod_e\bbR^6_e$, a space whose dimension is 6 times the number of edges in $E$. 
One can then form a ``Lagrangian'' $L_D=Q_D+\eps P_D$ by summing over all the crossings of $D$ local
contributions that involve only the variables associated with the four edges around each crossings,
and adding a ``correction'' which is a sum over the edges $e$ of $D$ of terms that depend only on the
rotation number of $e$ and on the variables in $\bbR^6_e$, such that
\[
  \Gint_{\bbR_{6E}} e^{L_D}
  = \Gint_{\bbR_{6E}} e^{Q_D+\eps P_D}
  = \frac{(2\pi)^{3|E|}}{\Delta_1\Delta_2\Delta_3}\exp(\eps\theta_0) + O(\eps^2),
\]
and such that the Feynman diagram expansion of the left-hand-side of the above equation becomes
precisely formula~\eqref{eq:Main} for $\theta$. See more about all this in~\cite{IType}.
\end{fact}

\begin{comment} \label{com:UC}
In fact, Fact~\ref{fact:PerturbedGaussian} is what
we initially predicted based on Discussion~\ref{disc:SolvApp}, along
with some further information about the ``shape'' of $P_D$. We used the
method of undetermined coefficients to find precise formulas for $P_D$,
and then the technique of Feynman diagrams to derive our main formula,
Equation~\ref{eq:Main}.
\end{comment}

\begin{dream} \label{dream:Seifert}
There is a ``Seifert formula'' for $\Theta$. More precisely,
let $K$ be a knot, let $\Sigma$ be a Seifert surface for
$K$, let $H\coloneq H_1(\Sigma;\bbR)$, and let $6H$ denote $H\oplus
H\oplus H\oplus H\oplus H\oplus H$. Let $Q_\Sigma$ denote 3 copies of the
standard Seifert form on $H\oplus H$, taken with parameters $T_1$, $T_2$,
and $T_3$; so $Q_\Sigma$ is a quadratic on $6H$. We dream that there a
``perturbation term'' $P_\Sigma$, a polynomial function on $6H$ defined
in terms of some low degree finite type invariants of various knotted
graphs formed by representatives of classes in $H$ (also taking account
of their intersections), such that
\[
  \Gint_{6H} e^{L_\Sigma}
  = \Gint_{6H} e^{Q_\Sigma+\eps P_\Sigma}
  = \frac{(2\pi)^{3\dim(H)}}{\Delta_1\Delta_2\Delta_3}\exp(\eps\theta_0) + O(\eps^2).
\]
\end{dream}

If this dream is true, it will probably re-prove Theorem~\ref{thm:Genus}
and prove Conjectures \ref{conj:D6}, \ref{conj:Mirror},
and~\ref{conj:Reverse} much as the Seifert formula for $\Delta$ can
be used to prove the genus bound provided by $\Delta$ and its basic
symmetry properties.

We note the relationship between this dream and~\cite[Theorem~4.4]{Ohtsuki:TwoLoop}.

\begin{dream} \label{dream:UniversalSeifert}
All the invariants from Discussion~\ref{disc:SolvApp} have Seifert formulas in the style of
Dream~\ref{dream:Seifert}. In fact, there ought to be a characterization of those Lagrangians
$L_\Sigma$ for which $\Gint e^{L_\Sigma}$ is a knot invariant, and there may be a construction of all
those Lagrangians which is intrinsic to topology and does not rely on the theory of Lie algebras.
\end{dream}

\parpic[r]{
  \adjustbox{valign=m}{\includegraphics[height=20mm]{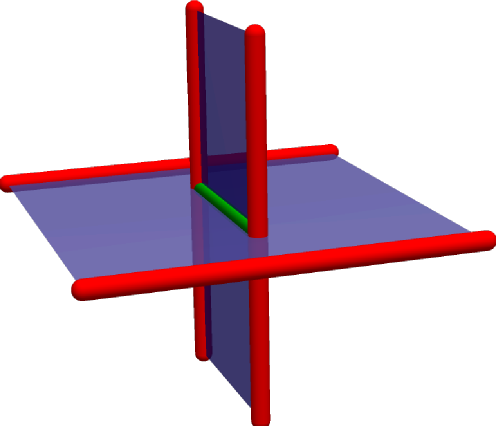}}
  $\to$
  \adjustbox{valign=m}{\includegraphics[height=20mm]{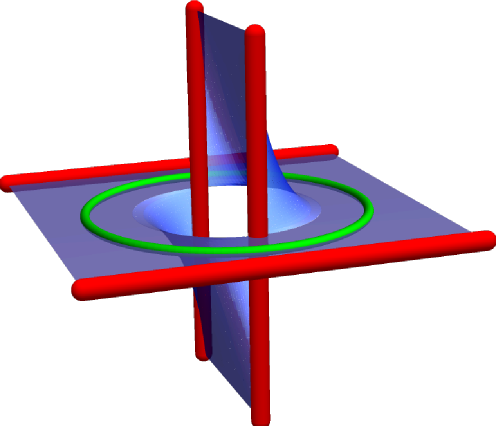}}
}
If a knot $K$ is ribbon then for some $g$ it has a Seifert surface $\Sigma$
of genus $g$ such that $g$ of the generators of $H_1(\Sigma)$ can be represented by a $g$-component
unlink (see the hint on the right, and see further details in \cite[Chapter~VIII]{Kauffman:OnKnots}
or in \cite[Section~3.4]{Bai:Alexander}). This implies that the Seifert matrix $M$ of $\Sigma$ has the
form $\begin{pmatrix}0&A\\A^*&B\end{pmatrix}$, which implies that the determinant of $M$, the
Alexander polynomial $\Delta$, satisfies the Fox-Milnor condition:

\begin{theorem}[Fox and Milnor,~\cite{FoxMilnor:CobordismOfKnots}]
If $K$ is a ribbon knot, then there exists some polynomial $f(T)$ such
that $\Delta=f(T)f(T^{-1})$.
\end{theorem}

\begin{dream} \label{dream:Ribbon}
Dream~\ref{dream:Seifert}, along with the fact that half the homology of a Seifert
surface of a ribbon knot can be represented by an unlink, will imply that $\theta$ takes a special
form on ribbon knots, giving us stronger powers to detect knots that are not ribbon.
\end{dream}

\begin{discussion} \label{disc:Links}
In this paper we concentrated on knots, yet at least partially,
$\Theta$ can be generalized also to links. Indeed, the definitions in
Section~\ref{sec:MainTheorem} and the proof in Section~\ref{sec:Proof}
go through provided the matrix $A$ is invertible; namely, provided the
Alexander polynomial $\Delta$ is non-zero (for knots, this is always
the case), and provided we choose one component of the link to cut open.

The programs of Section~\ref{sec:Implementation} fail for minor reasons,
and a fix is in~\cite[Theta4Links.nb]{Self}. Some results are in
Figure~\ref{fig:Theta4Links}. Preliminary testing using these programs
suggests that the resulting invariant is independent of the choice of
the cut component, but we did not prove that.

\begin{figure}
\adjustbox{valign=m}{\resizebox*{!}{4.67in}{
\def\l#1#2#3#4{{\parbox{0.5in}{\centering
  \vskip 2pt
  \includegraphics[width=\linewidth,height=\linewidth]{LinkFigs/L_#1#2_#4.pdf}
  \newline \href{https://katlas.org/wiki/L#1#2#3}{L#1#2\_#3}
  \vskip 2pt
}}}
\begin{tabular}{|c|c|c|c|c|c|c|c|c|c|}
  \hline \l{2}{a}{1}{1} & \l{4}{a}{1}{1} & \l{5}{a}{1}{1} & \l{6}{a}{1}{1} & \l{6}{a}{2}{2} & \l{6}{a}{3}{3} & \l{6}{a}{4}{4} & \l{6}{a}{5}{5} & \l{6}{n}{1}{1} & \l{7}{a}{1}{1} \\
  \hline \l{7}{a}{2}{2} & \l{7}{a}{3}{3} & \l{7}{a}{4}{4} & \l{7}{a}{5}{5} & \l{7}{a}{6}{6} & \l{7}{a}{7}{7} & \l{7}{n}{1}{1} & \l{7}{n}{2}{2} & \l{8}{a}{1}{01} & \l{8}{a}{2}{02} \\ 
  \hline \l{8}{a}{3}{03} & \l{8}{a}{4}{04} & \l{8}{a}{5}{05} & \l{8}{a}{6}{06} & \l{8}{a}{7}{07} & \l{8}{a}{8}{08} & \l{8}{a}{9}{09} & \l{8}{a}{10}{10} & \l{8}{a}{11}{11} & \l{8}{a}{12}{12} \\ 
  \hline \l{8}{a}{13}{13} & \l{8}{a}{14}{14} & \l{8}{a}{15}{15} & \l{8}{a}{16}{16} & \l{8}{a}{17}{17} & \l{8}{a}{18}{18} & \l{8}{a}{19}{19} & \l{8}{a}{20}{20} & \l{8}{a}{21}{21} & \l{8}{n}{1}{01} \\ 
  \hline \l{8}{n}{2}{02} & \l{8}{n}{3}{03} & \l{8}{n}{4}{04} & \l{8}{n}{5}{05} & \l{8}{n}{6}{06} & \l{8}{n}{7}{07} & \l{8}{n}{8}{08} & \l{9}{a}{1}{01} & \l{9}{a}{2}{02} & \l{9}{a}{3}{03} \\ 
  \hline \l{9}{a}{4}{04} & \l{9}{a}{5}{05} & \l{9}{a}{6}{06} & \l{9}{a}{7}{07} & \l{9}{a}{8}{08} & \l{9}{a}{9}{09} & \l{9}{a}{10}{10} & \l{9}{a}{11}{11} & \l{9}{a}{12}{12} & \l{9}{a}{13}{13} \\ 
  \hline \l{9}{a}{14}{14} & \l{9}{a}{15}{15} & \l{9}{a}{16}{16} & \l{9}{a}{17}{17} & \l{9}{a}{18}{18} & \l{9}{a}{19}{19} & \l{9}{a}{20}{20} & \l{9}{a}{21}{21} & \l{9}{a}{22}{22} & \l{9}{a}{23}{23} \\ 
  \hline \l{9}{a}{24}{24} & \l{9}{a}{25}{25} & \l{9}{a}{26}{26} & \l{9}{a}{27}{27} & \l{9}{a}{28}{28} & \l{9}{a}{29}{29} & \l{9}{a}{30}{30} & \l{9}{a}{31}{31} & \l{9}{a}{32}{32} & \l{9}{a}{33}{33} \\ 
  \hline \l{9}{a}{34}{34} & \l{9}{a}{35}{35} & \l{9}{a}{36}{36} & \l{9}{a}{37}{37} & \l{9}{a}{38}{38} & \l{9}{a}{39}{39} & \l{9}{a}{40}{40} & \l{9}{a}{41}{41} & \l{9}{a}{42}{42} & \l{9}{a}{43}{43} \\ 
  \hline \l{9}{a}{44}{44} & \l{9}{a}{45}{45} & \l{9}{a}{46}{46} & \l{9}{a}{47}{47} & \l{9}{a}{48}{48} & \l{9}{a}{49}{49} & \l{9}{a}{50}{50} & \l{9}{a}{51}{51} & \l{9}{a}{52}{52} & \l{9}{a}{53}{53} \\ 
  \hline \l{9}{a}{54}{54} & \l{9}{a}{55}{55} & \l{9}{n}{1}{01} & \l{9}{n}{2}{02} & \l{9}{n}{3}{03} & \l{9}{n}{4}{04} & \l{9}{n}{5}{05} & \l{9}{n}{6}{06} & \l{9}{n}{7}{07} & \l{9}{n}{8}{08} \\ 
  \hline \l{9}{n}{9}{09} & \l{9}{n}{10}{10} & \l{9}{n}{11}{11} & \l{9}{n}{12}{12} & \l{9}{n}{13}{13} & \l{9}{n}{14}{14} & \l{9}{n}{15}{15} & \l{9}{n}{16}{16} & \l{9}{n}{17}{17} & \l{9}{n}{18}{18} \\ 
  \hline \l{9}{n}{19}{19} & \l{9}{n}{20}{20} & \l{9}{n}{21}{21} & \l{9}{n}{22}{22} & \l{9}{n}{23}{23} & \l{9}{n}{24}{24} & \l{9}{n}{25}{25} & \l{9}{n}{26}{26} & \l{9}{n}{27}{27} & \l{9}{n}{28}{28} \\ 
  \hline
\end{tabular}
}}
$\underset{\Theta}{\rightarrow}$
\includegraphics[height=4.67in,valign=m]{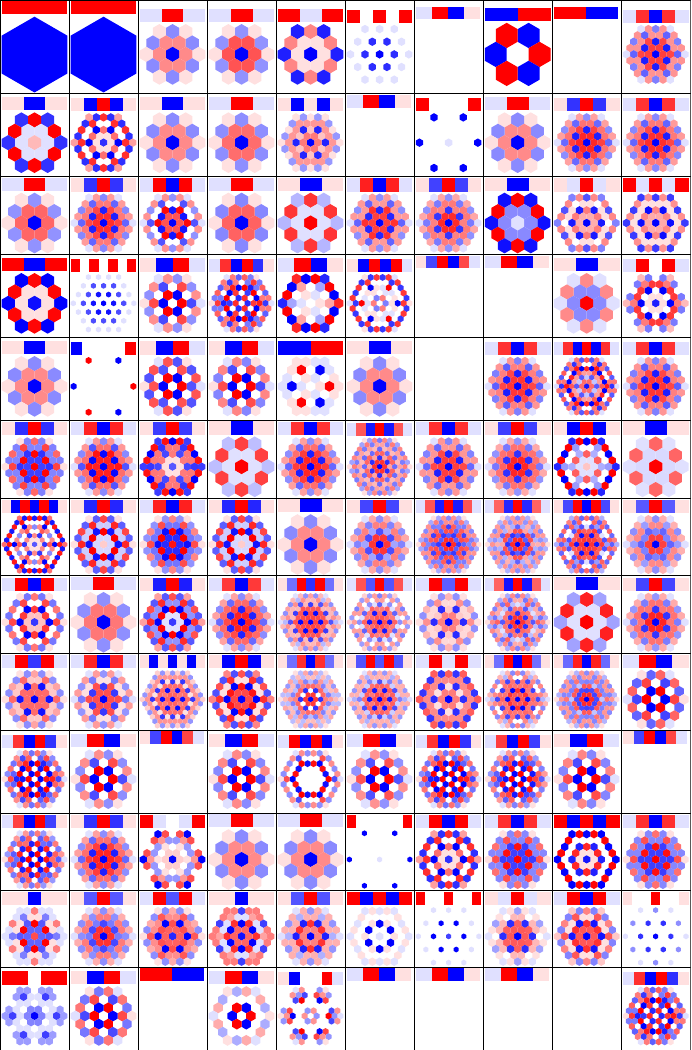}
\caption{
  $\Theta$ for all the prime links with up to 9 crossings, up to reflections and with arbitrary
  choices of strand orientations. Empty boxes correspond to links for which $\Delta=0$.
} \label{fig:Theta4Links}
\end{figure}

If $\Delta=0$, one may contemplate replacing $G=A^{-1}$ by the
adjugate matrix $\adj(A)$ of $A$ (the matrix of codimension 1
minors, which satisfies $A\cdot\adj(A)=\det(A)I$).\footnote{Similar
``adjugate'' reasoning shows that $\theta$ is always divisible
by $\Delta^{(2)}(T_1)\Delta^{(2)}(T_2)\Delta^{(2)}(T_3)$, where
$\Delta^{(2)}(T)$ is the second Alexander polynomial (e.g.~
\cite[Definition~8.10]{BurdeZieschang:Knots}).} Some preliminary testing is also
in~\cite[Theta4Links.nb]{Self}. Yet if $G$ is replaced with $\adj(A)$,
its equivalence with the $g$-rules (Equations~\eqref{eq:CarRules}
and~\eqref{eq:CounterRules}) breaks, and so we have no proof of
invariance. We may attempt to fix that in a future work, but it is not
done yet.

We note that the loop expansion of Conjecture~\ref{conj:TwoLoop}
does not predict that $\Theta$ should extend to links. We
also note that the solvable approximation technique of
Discussion~\ref{disc:SolvApp} does predict such an extension,
and in fact, it predicts more: that much like the Gassner
representation~\cite{Gassner:OnBraidGroups} and the multi-variable
Alexander polynomial (e.g.~\cite[Chapter~7]{Kawauchi:Srvery}),
there should be a multi-variable version of $\Theta$ which would
be a polynomial in $2m$ variables when evaluated on an $m$-component
link. We did not attempt to find explicit formulas for the multi-variable
$\Theta$. \endpar{\ref{disc:Links}}

\end{discussion}

Ever since Khovanov homology~\cite{Khovanov:Categorification, Bar-Natan:Categorification} it is almost mandatory to ask about
anything, ``does it categorify?''. $\Theta$ is not exempt:

\begin{question} Is there a categorification of $\theta$? Is there a finite triply-graded chain complex whose Euler
characteristic is $\theta$ and whose homology is invariant?
\end{question}

We note that $\theta$ is a neighbor of $\Delta$ (indeed they live
together within $\Theta$), and that $\Delta$ is categorified by knot Floer
homology~\cite{OzsvathSzabo:IntroductionToHF, Manolescu:KnotFloerHomology,
Juhasz:Survey}. Thus one may wonder if a categorification of $\theta$
will end up a neighbor of Floer knot homology. This applies even more
to a possible categorication of $g_{\alpha\beta}$:

\begin{question} Is there a categorification of
$\Delta\cdot\tilg_{ab}$? Is there a finite doubly-graded chain
complex whose Euler characteristic is $\Delta\cdot\tilg_{ab}$
and whose homology is a relative invariant in the sense of
Theorem~\ref{thm:RelativeInvariant}?
\end{question}

The latter seems likely: $\Delta\cdot\tilg_{ab}$ is, after all, a minor
of a matrix whose determinant is $\Delta$.

\section{Acknowledgement} We wish to thank Jorge Becerra and Stavros
Garoufalidis for comments and suggestions. This work was partially
supported by NSERC grants RGPIN-2018-04350 and RGPIN-2025-06718 and by
the Chu Family Foundation (NYC).

%% file: figs/SampleKnot.pdf_t
\begin{picture}(0,0)%
\includegraphics{figs/SampleKnot.pdf}%
\end{picture}%
%
%
\setlength{\unitlength}{3947sp}%
\begingroup\makeatletter\ifx\SetFigFont\undefined%
\gdef\SetFigFont#1#2#3#4#5{%
  \reset@font\fontsize{#1}{#2pt}%
  \fontfamily{#3}\fontseries{#4}\fontshape{#5}%
  \selectfont}%
\fi\endgroup%
\begin{picture}(1202,1749)(14,-898)
\put( 76,-798){\makebox(0,0)[b]{\smash{{\SetFigFont{10}{12.0}{\rmdefault}{\mddefault}{\itdefault}{\color[rgb]{0,0,0}1}%
}}}}
\put(526,-312){\makebox(0,0)[b]{\smash{{\SetFigFont{10}{12.0}{\rmdefault}{\mddefault}{\itdefault}{\color[rgb]{0,0,0}2}%
}}}}
\put( 76,130){\makebox(0,0)[b]{\smash{{\SetFigFont{10}{12.0}{\rmdefault}{\mddefault}{\itdefault}{\color[rgb]{0,0,0}3}%
}}}}
\put(976,-135){\makebox(0,0)[b]{\smash{{\SetFigFont{10}{12.0}{\rmdefault}{\mddefault}{\itdefault}{\color[rgb]{0,0,0}4}%
}}}}
\put(526,130){\makebox(0,0)[b]{\smash{{\SetFigFont{10}{12.0}{\rmdefault}{\mddefault}{\itdefault}{\color[rgb]{0,0,0}6}%
}}}}
\put( 76,617){\makebox(0,0)[b]{\smash{{\SetFigFont{10}{12.0}{\rmdefault}{\mddefault}{\itdefault}{\color[rgb]{0,0,0}7}%
}}}}
\put( 76,-312){\makebox(0,0)[b]{\smash{{\SetFigFont{10}{12.0}{\rmdefault}{\mddefault}{\itdefault}{\color[rgb]{0,0,0}5}%
}}}}
\put(1201,706){\makebox(0,0)[b]{\smash{{\SetFigFont{11}{13.2}{\familydefault}{\mddefault}{\updefault}{\color[rgb]{0,0,0}$D$}%
}}}}
\put(1126,-85){\makebox(0,0)[lb]{\smash{{\SetFigFont{10}{12.0}{\rmdefault}{\mddefault}{\itdefault}{\color[rgb]{0,0,0}$\varphi_4=-1$}%
}}}}
\end{picture}%

%% file: figs/Xings.pdf_t
\begin{picture}(0,0)%
\includegraphics{figs/Xings.pdf}%
\end{picture}%
%
%
\setlength{\unitlength}{3947sp}%
\begingroup\makeatletter\ifx\SetFigFont\undefined%
\gdef\SetFigFont#1#2#3#4#5{%
  \reset@font\fontsize{#1}{#2pt}%
  \fontfamily{#3}\fontseries{#4}\fontshape{#5}%
  \selectfont}%
\fi\endgroup%
\begin{picture}(2355,1042)(436,-191)
\put(2476,-136){\makebox(0,0)[b]{\smash{{\SetFigFont{10}{12.0}{\familydefault}{\mddefault}{\updefault}{\color[rgb]{0,0,0}$s=-1$}%
}}}}
\put(2176, 89){\makebox(0,0)[rb]{\smash{{\SetFigFont{10}{12.0}{\familydefault}{\mddefault}{\updefault}{\color[rgb]{0,0,0}$j$}%
}}}}
\put(2176,645){\makebox(0,0)[rb]{\smash{{\SetFigFont{10}{12.0}{\familydefault}{\mddefault}{\updefault}{\color[rgb]{0,0,0}$i+1$}%
}}}}
\put(2776,645){\makebox(0,0)[lb]{\smash{{\SetFigFont{10}{12.0}{\familydefault}{\mddefault}{\updefault}{\color[rgb]{0,0,0}$j+1$}%
}}}}
\put(2776, 89){\makebox(0,0)[lb]{\smash{{\SetFigFont{10}{12.0}{\familydefault}{\mddefault}{\updefault}{\color[rgb]{0,0,0}$i$}%
}}}}
\put(451, 89){\makebox(0,0)[rb]{\smash{{\SetFigFont{10}{12.0}{\familydefault}{\mddefault}{\updefault}{\color[rgb]{0,0,0}$i$}%
}}}}
\put(751,-136){\makebox(0,0)[b]{\smash{{\SetFigFont{10}{12.0}{\familydefault}{\mddefault}{\updefault}{\color[rgb]{0,0,0}$s=+1$}%
}}}}
\put(451,645){\makebox(0,0)[rb]{\smash{{\SetFigFont{10}{12.0}{\familydefault}{\mddefault}{\updefault}{\color[rgb]{0,0,0}$j+1$}%
}}}}
\put(1051,645){\makebox(0,0)[lb]{\smash{{\SetFigFont{10}{12.0}{\familydefault}{\mddefault}{\updefault}{\color[rgb]{0,0,0}$i+1$}%
}}}}
\put(1051, 89){\makebox(0,0)[lb]{\smash{{\SetFigFont{10}{12.0}{\familydefault}{\mddefault}{\updefault}{\color[rgb]{0,0,0}$j$}%
}}}}
\end{picture}%

%% file: figs/CarRules2.pdf_t
\begin{picture}(0,0)%
\includegraphics{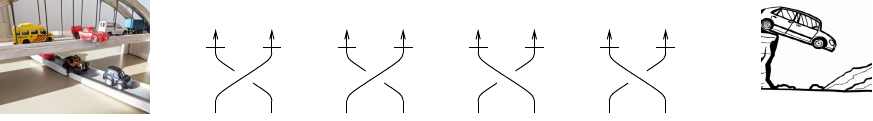}%
\end{picture}%
%
%
\setlength{\unitlength}{3947sp}%
\begingroup\makeatletter\ifx\SetFigFont\undefined%
\gdef\SetFigFont#1#2#3#4#5{%
  \reset@font\fontsize{#1}{#2pt}%
  \fontfamily{#3}\fontseries{#4}\fontshape{#5}%
  \selectfont}%
\fi\endgroup%
\begin{picture}(6975,962)(-299,-123)
\put(258,-19){\makebox(0,0)[b]{\smash{{\SetFigFont{9}{10.8}{\rmdefault}{\mddefault}{\updefault}{\color[rgb]{0,0,0}$p=1-T^s$}%
}}}}
\put(6185, 30){\makebox(0,0)[b]{\smash{{\SetFigFont{11}{13.2}{\rmdefault}{\mddefault}{\updefault}{\color[rgb]{0,0,0}\DallEA}%
}}}}
\put(6185,-108){\makebox(0,0)[b]{\smash{{\SetFigFont{11}{13.2}{\rmdefault}{\mddefault}{\updefault}{\color[rgb]{0,0,0}\DallEB}%
}}}}
\put(1426,-61){\makebox(0,0)[lb]{\smash{{\SetFigFont{10}{12.0}{\rmdefault}{\mddefault}{\updefault}{\color[rgb]{0,0,0}\rac}%
}}}}
\put(2926,-61){\makebox(0,0)[rb]{\smash{{\SetFigFont{10}{12.0}{\rmdefault}{\mddefault}{\updefault}{\color[rgb]{0,0,0}\car}%
}}}}
\put(3526,-61){\makebox(0,0)[lb]{\smash{{\SetFigFont{10}{12.0}{\rmdefault}{\mddefault}{\updefault}{\color[rgb]{0,0,0}\rac}%
}}}}
\put(5026,-61){\makebox(0,0)[rb]{\smash{{\SetFigFont{10}{12.0}{\rmdefault}{\mddefault}{\updefault}{\color[rgb]{0,0,0}\car}%
}}}}
\put(1426,689){\makebox(0,0)[b]{\smash{{\SetFigFont{9}{10.8}{\rmdefault}{\mddefault}{\updefault}{\color[rgb]{0,0,0}$1\!-\!T$}%
}}}}
\put(1876,689){\makebox(0,0)[b]{\smash{{\SetFigFont{9}{10.8}{\rmdefault}{\mddefault}{\updefault}{\color[rgb]{0,0,0}$T$}%
}}}}
\put(2476,689){\makebox(0,0)[b]{\smash{{\SetFigFont{9}{10.8}{\rmdefault}{\mddefault}{\updefault}{\color[rgb]{0,0,0}$1$}%
}}}}
\put(2926,689){\makebox(0,0)[b]{\smash{{\SetFigFont{9}{10.8}{\rmdefault}{\mddefault}{\updefault}{\color[rgb]{0,0,0}$0$}%
}}}}
\put(3526,689){\makebox(0,0)[b]{\smash{{\SetFigFont{9}{10.8}{\rmdefault}{\mddefault}{\updefault}{\color[rgb]{0,0,0}$0$}%
}}}}
\put(3976,689){\makebox(0,0)[b]{\smash{{\SetFigFont{9}{10.8}{\rmdefault}{\mddefault}{\updefault}{\color[rgb]{0,0,0}$1$}%
}}}}
\put(4576,689){\makebox(0,0)[b]{\smash{{\SetFigFont{9}{10.8}{\rmdefault}{\mddefault}{\updefault}{\color[rgb]{0,0,0}$T^{-1}$}%
}}}}
\put(5101,689){\makebox(0,0)[b]{\smash{{\SetFigFont{9}{10.8}{\rmdefault}{\mddefault}{\updefault}{\color[rgb]{0,0,0}$1\!-\!T^{-1}$}%
}}}}
\put(5878,239){\makebox(0,0)[lb]{\smash{{\SetFigFont{8}{9.6}{\rmdefault}{\mddefault}{\updefault}{\color[rgb]{0,0,0}$2n+1$}%
}}}}
\end{picture}%

%% file: figs/FeynmanDiagrams.pdf_t
\begin{picture}(0,0)%
\includegraphics{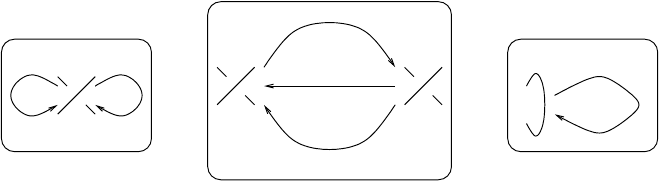}%
\end{picture}%
%
%
\setlength{\unitlength}{3947sp}%
\begingroup\makeatletter\ifx\SetFigFont\undefined%
\gdef\SetFigFont#1#2#3#4#5{%
  \reset@font\fontsize{#1}{#2pt}%
  \fontfamily{#3}\fontseries{#4}\fontshape{#5}%
  \selectfont}%
\fi\endgroup%
\begin{picture}(5274,1449)(-311,-373)
\put(-149,614){\makebox(0,0)[b]{\smash{{\SetFigFont{10}{12.0}{\rmdefault}{\mddefault}{\updefault}{\color[rgb]{0,0,0}$D$}%
}}}}
\put(151, 14){\makebox(0,0)[b]{\smash{{\SetFigFont{10}{12.0}{\rmdefault}{\mddefault}{\updefault}{\color[rgb]{0,0,0}$i$}%
}}}}
\put(751, 14){\makebox(0,0)[b]{\smash{{\SetFigFont{10}{12.0}{\rmdefault}{\mddefault}{\updefault}{\color[rgb]{0,0,0}$\car_\mu$}%
}}}}
\put(-149, 14){\makebox(0,0)[b]{\smash{{\SetFigFont{10}{12.0}{\rmdefault}{\mddefault}{\updefault}{\color[rgb]{0,0,0}$\car_\nu$}%
}}}}
\put(451, 14){\makebox(0,0)[b]{\smash{{\SetFigFont{10}{12.0}{\rmdefault}{\mddefault}{\updefault}{\color[rgb]{0,0,0}$j$}%
}}}}
\put(3901,614){\makebox(0,0)[b]{\smash{{\SetFigFont{10}{12.0}{\rmdefault}{\mddefault}{\updefault}{\color[rgb]{0,0,0}$D$}%
}}}}
\put(4051,-61){\makebox(0,0)[b]{\smash{{\SetFigFont{10}{12.0}{\rmdefault}{\mddefault}{\updefault}{\color[rgb]{0,0,0}$k$}%
}}}}
\put(4726,-61){\makebox(0,0)[b]{\smash{{\SetFigFont{10}{12.0}{\rmdefault}{\mddefault}{\updefault}{\color[rgb]{0,0,0}$\car_3$}%
}}}}
\put(3976,239){\makebox(0,0)[b]{\smash{{\SetFigFont{10}{12.0}{\rmdefault}{\mddefault}{\updefault}{\color[rgb]{0,0,0}$\varphi$}%
}}}}
\put(1726, 89){\makebox(0,0)[b]{\smash{{\SetFigFont{10}{12.0}{\rmdefault}{\mddefault}{\updefault}{\color[rgb]{0,0,0}$j_0$}%
}}}}
\put(2926, 89){\makebox(0,0)[b]{\smash{{\SetFigFont{10}{12.0}{\rmdefault}{\mddefault}{\updefault}{\color[rgb]{0,0,0}$i_1$}%
}}}}
\put(1501,614){\makebox(0,0)[b]{\smash{{\SetFigFont{10}{12.0}{\rmdefault}{\mddefault}{\updefault}{\color[rgb]{0,0,0}$D$}%
}}}}
\put(1426, 89){\makebox(0,0)[b]{\smash{{\SetFigFont{10}{12.0}{\rmdefault}{\mddefault}{\updefault}{\color[rgb]{0,0,0}$i_0$}%
}}}}
\put(3226, 89){\makebox(0,0)[b]{\smash{{\SetFigFont{10}{12.0}{\rmdefault}{\mddefault}{\updefault}{\color[rgb]{0,0,0}$j_1$}%
}}}}
\put(2401,-286){\makebox(0,0)[b]{\smash{{\SetFigFont{10}{12.0}{\rmdefault}{\mddefault}{\updefault}{\color[rgb]{0,0,0}$\rac_2$}%
}}}}
\put(2326,239){\makebox(0,0)[b]{\smash{{\SetFigFont{10}{12.0}{\rmdefault}{\mddefault}{\updefault}{\color[rgb]{0,0,0}$\rac_1$}%
}}}}
\put(2026,914){\makebox(0,0)[b]{\smash{{\SetFigFont{10}{12.0}{\rmdefault}{\mddefault}{\updefault}{\color[rgb]{0,0,0}$\car_3$}%
}}}}
\end{picture}%

%% file: Implementation.tex
We start by loading the package \verb$KnotTheory`$ --- it is only needed because it has many specific knots pre-defined. In this Section and in the next, $\input{figs/face.pdf_t}$ and \scalebox{0.66}{$\input{figs/human.pdf_t}$} mean ``human input'' while $\input{figs/machine.pdf_t}$ means ``computer output'':
{\def\nbpdfPostInput{\hfill}

\nbpdfInput{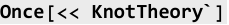}
\nbpdfPrint{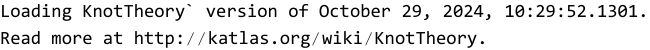}
}
Next we quietly define the modules \verb$Rot$, used to compute rotation numbers, and \verb$PolyPlot$, used to plot polynomials as bar codes and as hexagonal QR codes. Neither is a part of the core of the computation of $\Theta$, so neither is shown; yet we do show one usage example for each.

\nbpdfInput{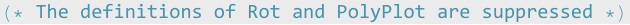}
{\def\nbpdfPostInput{\hfill}

\nbpdfInput{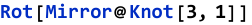}
\nbpdfMessage{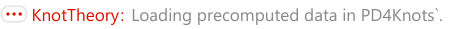}
\nbpdfOutput{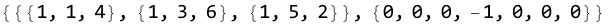}
}
We urge the reader to compare the above output with the knot diagram in Figure~\ref{fig:SampleDiagram}.

{\def\nbpdfPostInput{\hfill}

\nbpdfInput{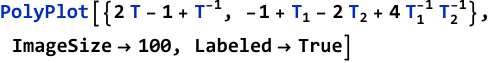}
\nbpdfOutput{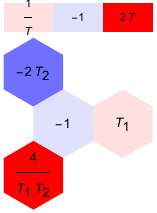}
}
The definition of \verb$CF$ below is a technicality telling the computer how to best store polynomials in the $g_{\nu\alpha\beta}$'s such as $F_1$ and $F_2$. The programs would run just the same without it, albeit a bit more slowly:

\nbpdfInput{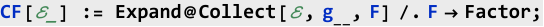}
Next, we decree that $T_3=T_1T_2$ and define the three ``Feynman Diagram'' polynomials $F_1$, $F_2$, and $F_3$:

\nbpdfInput{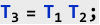}
\nbpdfInput{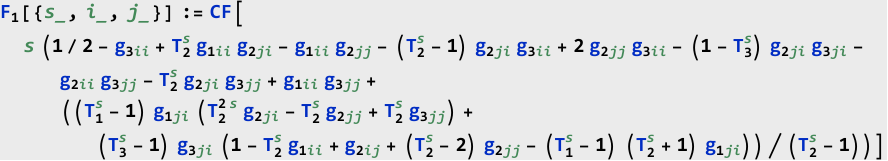}
\nbpdfInput{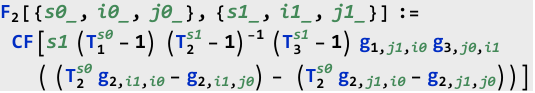}
\nbpdfInput{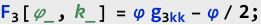}
Next comes the main program computing $\Theta(K)$. Fortunately, it matches perfectly with the mathematical description in Section~\ref{sec:MainTheorem}. In line 1 below we use \verb$Rot$ to let $X$ and $\varphi$ be the crossings and rotation numbers of $K$. In addition we let $n$ be the length of $X$, namely, the number of crossings in $K$, and we let the starting value of $A$ be the $(2n+1)\times(2n+1)$ identity matrix. Then in line 2, for each crossing in $X$ we add to $A$ a $2\times 2$ block, in rows $i$ and $j$ and columns $i+1$ and $j+1$, as explain in Equation~\eqref{eq:A}. In line 3 we compute the normalized Alexander polynomial $\Delta$ as in~\eqref{eq:Delta}. In line 4 we let $G$ be the inverse of $A$. In line 5 we declare what it means to evaluate, \verb$ev$, a formula $\calE$ that may contain symbols of the form $g_{\nu\alpha\beta}$: each such symbol is to be replaced by the entry in position $\alpha,\beta$ of $G$, but with $T$ replaced with $T_\nu$. In line 6 we start computing $\theta$ by computing the first summand in~\eqref{eq:Main}, which in itself, is a sum over the crossings of the knot. In line 7 we add to $\theta$ the double sum corresponding to the second term in~\eqref{eq:Main}, and in line 8, we add the third summand of~\eqref{eq:Main}. Finally, line 9 outputs a pair: $\Delta$, and the re-normalized version of $\theta$.

\nbpdfInput{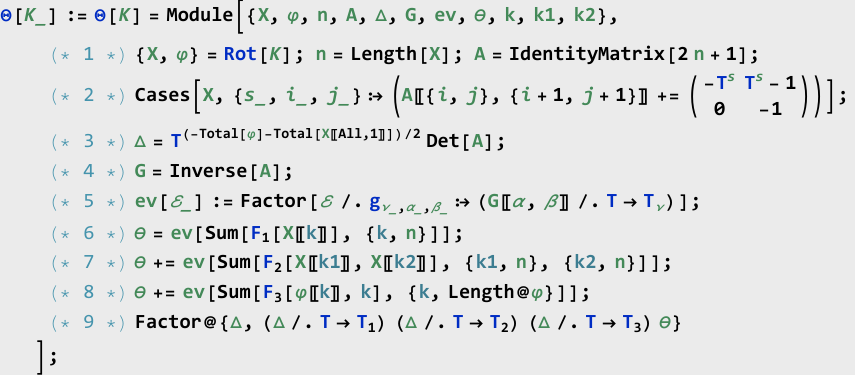}
\subsection{Examples} \label{ssec:Examples}
On to examples! Starting with the trefoil knot.

\nbpdfInput{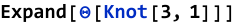}
\nbpdfMessage{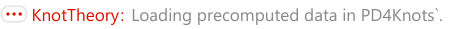}
\nbpdfOutput{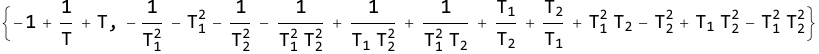}
{\def\nbpdfPostInput{\hfill}

\nbpdfInput{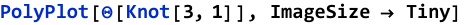}
\nbpdfMessage{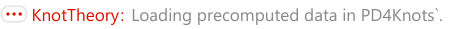}
\nbpdfOutput{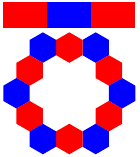}
}

\parpic[r]{\parbox{42mm}{
  \includegraphics[width=20mm]{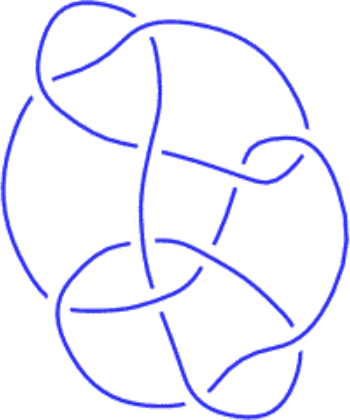}\hfill\includegraphics[width=20mm]{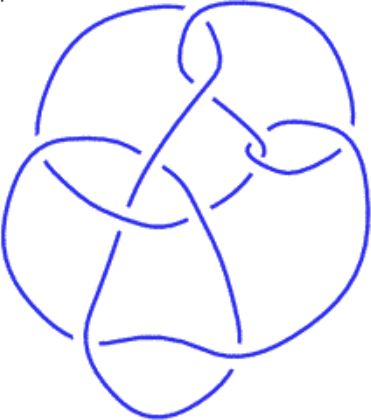}
}}
Next are the Conway knot $11_{n34}$ and the Kinoshita-Terasaka knot $11_{n42}$. The two are mutants and famously hard to separate: they both have $\Delta=1$ (as evidenced by their one-bar Alexander bar codes below), and they have the same hyperbolic volume, HOMFLY-PT polynomial, and Khovanov homology. Yet their $\theta$ invariants are different. Note that the genus of the Conway knot is 3, while the genus of the Kinoshita-Terasaka knot is 2. This agrees with the higher complexity of the QR code of the Conway polynomial and with Theorem~\ref{thm:Genus} below.
{\def\nbpdfPostInput{\hfill}

\nbpdfInput{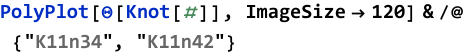}
\nbpdfMessage{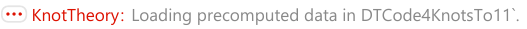}
\nbpdfMessage{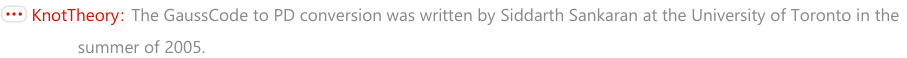}
\nbpdfOutput{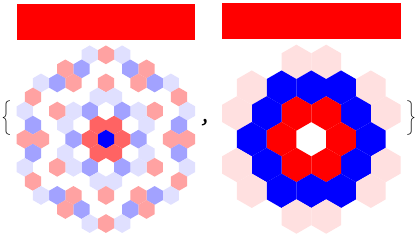}
}

\begin{figure}

\nbpdfInput{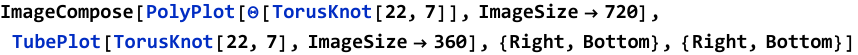}
\nbpdfOutput{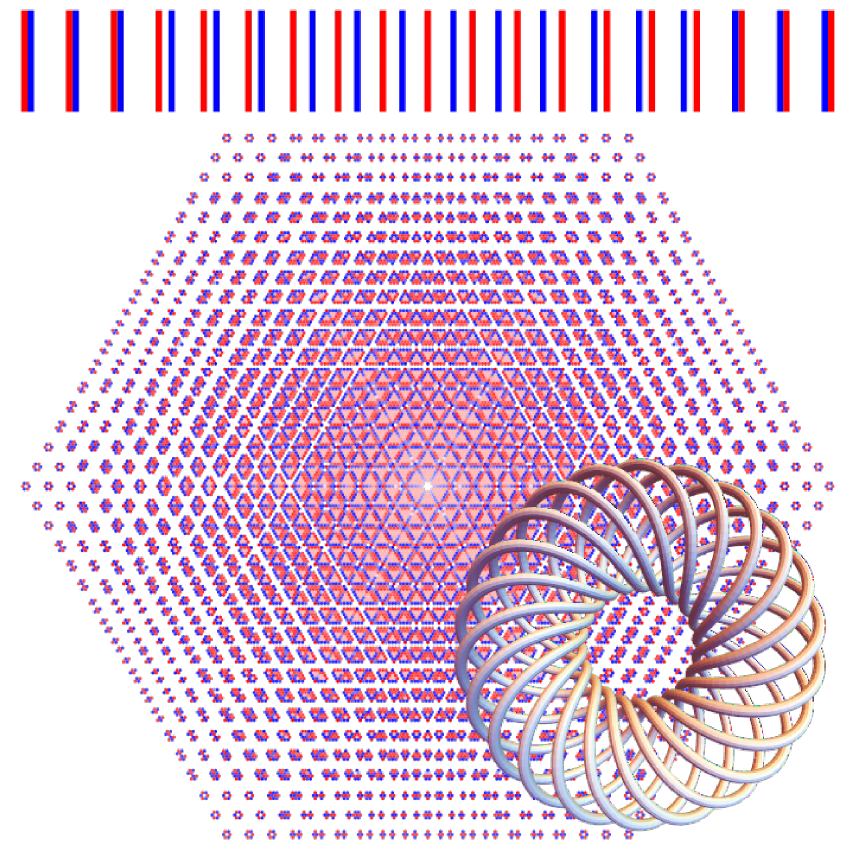}
\caption{The 132-crossing torus knot $T_{22/7}$ and a plot of its $\Theta$ invariant} \label{fig:T227}
\end{figure}

Torus knots have particularly nice-looking $\Theta$ invariants. Here are the torus knots $T_{13/2}$, $T_{17/3}$, $T_{13/5}$, and $T_{7/6}$:

\nbpdfInput{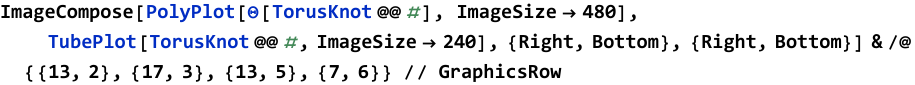}
\nbpdfOutput{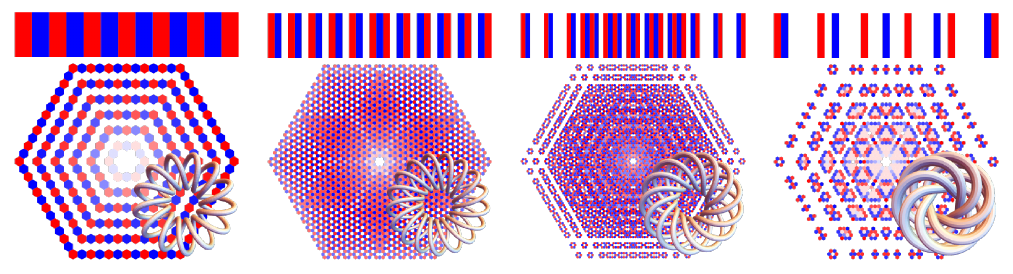}
The next line shows the computation time in seconds for the 132-crossing torus knot $T_{22/7}$ on a 2024 laptop, without actually showing the output. The output plot is in Figure~\ref{fig:T227}.
{\def\nbpdfPostInput{\hfill}

\nbpdfInput{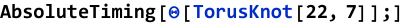}
\nbpdfOutput{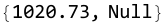}
}

We note that if $T_1$ and $T_2$ are assigned specific rational numbers and if the program for $\Theta$ is slightly modified so as to compute each $G_\nu$ separately (rather than computing $G$ symbolically and then substituting $T\to T_\nu$), then the program becomes significantly more efficient, for inverting a numerical matrix is cheaper than inverting a symbolic matrix (but then one obtains numerical answers and the beauty and the topological significance (Section~\ref{sec:SandM}) are lost). The Mathematica notebook that accompanies this paper, \cite[Theta.nb]{Self}, contains the required modified program as well as a few computational examples. One finds that with $T_1=22/7$ and $T_2=21/13$, the invariant $\Theta$ can be computed for knots with 600 crossings, and that for knots with up to 15 crossings, its separation power remains the same.

If $T_1$ and $T_2$ are assigned approximate real values, say $\pi$ and $e$ computed to 100 decimal digits, then $\Theta$ can be computed on knots with 1,000 crossings and, for knots with up to 15 crossings it remains very strong. But approximate real numbers are a bit thorny. It is hard to know how far one needs to compute before deciding that two such numbers are equal, and when two such numbers appear unequal, it is hard to tell if that is merely because they were computed differently and different roundings were applied. Thorns and snares are in the way of the perverse; He who guards his soul will be far from them (Proverbs 22:5)\footnote{\quad.\cjRL{/swmr np/sw yr.hq}}.

%% file: figs/Xingsp.pdf_t
\begin{picture}(0,0)%
\includegraphics{figs/Xingsp.pdf}%
\end{picture}%
%
%
\setlength{\unitlength}{3947sp}%
\begingroup\makeatletter\ifx\SetFigFont\undefined%
\gdef\SetFigFont#1#2#3#4#5{%
  \reset@font\fontsize{#1}{#2pt}%
  \fontfamily{#3}\fontseries{#4}\fontshape{#5}%
  \selectfont}%
\fi\endgroup%
\begin{picture}(2205,1042)(436,-191)
\put(2326,-136){\makebox(0,0)[b]{\smash{{\SetFigFont{10}{12.0}{\familydefault}{\mddefault}{\updefault}{\color[rgb]{0,0,0}$s=-1$}%
}}}}
\put(2026, 89){\makebox(0,0)[rb]{\smash{{\SetFigFont{10}{12.0}{\familydefault}{\mddefault}{\updefault}{\color[rgb]{0,0,0}$j$}%
}}}}
\put(2026,645){\makebox(0,0)[rb]{\smash{{\SetFigFont{10}{12.0}{\familydefault}{\mddefault}{\updefault}{\color[rgb]{0,0,0}$\ip$}%
}}}}
\put(2626,645){\makebox(0,0)[lb]{\smash{{\SetFigFont{10}{12.0}{\familydefault}{\mddefault}{\updefault}{\color[rgb]{0,0,0}$\jp$}%
}}}}
\put(2626, 89){\makebox(0,0)[lb]{\smash{{\SetFigFont{10}{12.0}{\familydefault}{\mddefault}{\updefault}{\color[rgb]{0,0,0}$i$}%
}}}}
\put(451, 89){\makebox(0,0)[rb]{\smash{{\SetFigFont{10}{12.0}{\familydefault}{\mddefault}{\updefault}{\color[rgb]{0,0,0}$i$}%
}}}}
\put(751,-136){\makebox(0,0)[b]{\smash{{\SetFigFont{10}{12.0}{\familydefault}{\mddefault}{\updefault}{\color[rgb]{0,0,0}$s=+1$}%
}}}}
\put(451,645){\makebox(0,0)[rb]{\smash{{\SetFigFont{10}{12.0}{\familydefault}{\mddefault}{\updefault}{\color[rgb]{0,0,0}$\jp$}%
}}}}
\put(1051,645){\makebox(0,0)[lb]{\smash{{\SetFigFont{10}{12.0}{\familydefault}{\mddefault}{\updefault}{\color[rgb]{0,0,0}$\ip$}%
}}}}
\put(1051, 89){\makebox(0,0)[lb]{\smash{{\SetFigFont{10}{12.0}{\familydefault}{\mddefault}{\updefault}{\color[rgb]{0,0,0}$j$}%
}}}}
\end{picture}%

%% file: figs/NullVertex.pdf_t
\begin{picture}(0,0)%
\includegraphics{figs/NullVertex.pdf}%
\end{picture}%
%
%
\setlength{\unitlength}{3947sp}%
\begingroup\makeatletter\ifx\SetFigFont\undefined%
\gdef\SetFigFont#1#2#3#4#5{%
  \reset@font\fontsize{#1}{#2pt}%
  \fontfamily{#3}\fontseries{#4}\fontshape{#5}%
  \selectfont}%
\fi\endgroup%
\begin{picture}(1149,343)(-11,508)
\put(226,614){\makebox(0,0)[b]{\smash{{\SetFigFont{10}{12.0}{\rmdefault}{\mddefault}{\updefault}{\color[rgb]{0,0,0}$j$}%
}}}}
\put(826,614){\makebox(0,0)[b]{\smash{{\SetFigFont{10}{12.0}{\rmdefault}{\mddefault}{\updefault}{\color[rgb]{0,0,0}$k$}%
}}}}
\end{picture}%

%% file: figs/RelativeInvariance.pdf_t
\begin{picture}(0,0)%
\includegraphics{figs/RelativeInvariance.pdf}%
\end{picture}%
%
%
\setlength{\unitlength}{4342sp}%
\begingroup\makeatletter\ifx\SetFigFont\undefined%
\gdef\SetFigFont#1#2#3#4#5{%
  \reset@font\fontsize{#1}{#2pt}%
  \fontfamily{#3}\fontseries{#4}\fontshape{#5}%
  \selectfont}%
\fi\endgroup%
\begin{picture}(4483,924)(30,-73)
\put(1756,695){\makebox(0,0)[b]{\smash{{\SetFigFont{9}{10.8}{\rmdefault}{\mddefault}{\updefault}{\color[rgb]{0,0,0}R1}%
}}}}
\put(1744,199){\makebox(0,0)[b]{\smash{{\SetFigFont{9}{10.8}{\rmdefault}{\mddefault}{\updefault}{\color[rgb]{0,0,0}R2}%
}}}}
\put(637,482){\makebox(0,0)[b]{\smash{{\SetFigFont{9}{10.8}{\rmdefault}{\mddefault}{\updefault}{\color[rgb]{0,0,0}R3}%
}}}}
\put(3226, 89){\makebox(0,0)[b]{\smash{{\SetFigFont{11}{13.2}{\rmdefault}{\mddefault}{\updefault}{\color[rgb]{0,0,0}$a$}%
}}}}
\put(4276,314){\makebox(0,0)[b]{\smash{{\SetFigFont{11}{13.2}{\rmdefault}{\mddefault}{\updefault}{\color[rgb]{0,0,0}$b$}%
}}}}
\put(2701,689){\makebox(0,0)[b]{\smash{{\SetFigFont{11}{13.2}{\rmdefault}{\mddefault}{\updefault}{\color[rgb]{0,0,0}$D$}%
}}}}
\put(3601,552){\makebox(0,0)[b]{\smash{{\SetFigFont{11}{13.2}{\rmdefault}{\mddefault}{\updefault}{\color[rgb]{0,0,0}$\tilg_{ab}$}%
}}}}
\end{picture}%

%% file: figs/RMoves.pdf_t
\begin{picture}(0,0)%
\includegraphics{figs/RMoves.pdf}%
\end{picture}%
%
%
\setlength{\unitlength}{3947sp}%
\begingroup\makeatletter\ifx\SetFigFont\undefined%
\gdef\SetFigFont#1#2#3#4#5{%
  \reset@font\fontsize{#1}{#2pt}%
  \fontfamily{#3}\fontseries{#4}\fontshape{#5}%
  \selectfont}%
\fi\endgroup%
\begin{picture}(6790,2706)(56,-1855)
\put(1351,239){\makebox(0,0)[b]{\smash{{\SetFigFont{10}{12.0}{\rmdefault}{\mddefault}{\updefault}{\color[rgb]{0,0,0}R1r}%
}}}}
\put(751,239){\makebox(0,0)[b]{\smash{{\SetFigFont{10}{12.0}{\rmdefault}{\mddefault}{\updefault}{\color[rgb]{0,0,0}R1l}%
}}}}
\put(1051,-361){\makebox(0,0)[b]{\smash{{\SetFigFont{11}{13.2}{\rmdefault}{\mddefault}{\updefault}{\color[rgb]{0,0,0}R1-left and R1-right}%
}}}}
\put(3601,164){\makebox(0,0)[b]{\smash{{\SetFigFont{10}{12.0}{\rmdefault}{\mddefault}{\updefault}{\color[rgb]{0,0,0}R3b}%
}}}}
\put(3601,-361){\makebox(0,0)[b]{\smash{{\SetFigFont{11}{13.2}{\rmdefault}{\mddefault}{\updefault}{\color[rgb]{0,0,0}braid-like R3}%
}}}}
\put(1351,-1261){\makebox(0,0)[b]{\smash{{\SetFigFont{10}{12.0}{\rmdefault}{\mddefault}{\updefault}{\color[rgb]{0,0,0}R2c$^+$}%
}}}}
\put(3301,-1261){\makebox(0,0)[b]{\smash{{\SetFigFont{10}{12.0}{\rmdefault}{\mddefault}{\updefault}{\color[rgb]{0,0,0}R2c$^-$}%
}}}}
\put(2251,-1786){\makebox(0,0)[b]{\smash{{\SetFigFont{11}{13.2}{\rmdefault}{\mddefault}{\updefault}{\color[rgb]{0,0,0}counterclockwise and clockwise cyclic R2}%
}}}}
\put(5176,689){\makebox(0,0)[lb]{\smash{{\SetFigFont{11}{13.2}{\rmdefault}{\mddefault}{\updefault}{\color[rgb]{0,0,0}Aside 1:}%
}}}}
\put(5176,-736){\makebox(0,0)[lb]{\smash{{\SetFigFont{11}{13.2}{\rmdefault}{\mddefault}{\updefault}{\color[rgb]{0,0,0}Aside 2:}%
}}}}
\end{picture}%

%% file: figs/RIR3.pdf_t
\begin{picture}(0,0)%
\includegraphics{figs/RIR3.pdf}%
\end{picture}%
%
%
\setlength{\unitlength}{3947sp}%
\begingroup\makeatletter\ifx\SetFigFont\undefined%
\gdef\SetFigFont#1#2#3#4#5{%
  \reset@font\fontsize{#1}{#2pt}%
  \fontfamily{#3}\fontseries{#4}\fontshape{#5}%
  \selectfont}%
\fi\endgroup%
\begin{picture}(7689,2401)(-26,-1550)
\put(451,-1111){\makebox(0,0)[b]{\smash{{\SetFigFont{10}{12.0}{\familydefault}{\mddefault}{\updefault}{\color[rgb]{0,0,0}$\cdots$}%
}}}}
\put(2326,-1111){\makebox(0,0)[b]{\smash{{\SetFigFont{10}{12.0}{\familydefault}{\mddefault}{\updefault}{\color[rgb]{0,0,0}$\cdots$}%
}}}}
\put(4351,-1111){\makebox(0,0)[b]{\smash{{\SetFigFont{10}{12.0}{\familydefault}{\mddefault}{\updefault}{\color[rgb]{0,0,0}$\cdots$}%
}}}}
\put(6226,-1111){\makebox(0,0)[b]{\smash{{\SetFigFont{10}{12.0}{\familydefault}{\mddefault}{\updefault}{\color[rgb]{0,0,0}$\cdots$}%
}}}}
\put(2326,-1336){\makebox(0,0)[b]{\smash{{\SetFigFont{10}{12.0}{\familydefault}{\mddefault}{\updefault}{\color[rgb]{0,0,0}further}%
}}}}
\put(2326,-1486){\makebox(0,0)[b]{\smash{{\SetFigFont{10}{12.0}{\familydefault}{\mddefault}{\updefault}{\color[rgb]{0,0,0}$g$-rules}%
}}}}
\put(4351,-1336){\makebox(0,0)[b]{\smash{{\SetFigFont{10}{12.0}{\familydefault}{\mddefault}{\updefault}{\color[rgb]{0,0,0}further}%
}}}}
\put(4351,-1486){\makebox(0,0)[b]{\smash{{\SetFigFont{10}{12.0}{\familydefault}{\mddefault}{\updefault}{\color[rgb]{0,0,0}crossings}%
}}}}
\put(6226,-1336){\makebox(0,0)[b]{\smash{{\SetFigFont{10}{12.0}{\familydefault}{\mddefault}{\updefault}{\color[rgb]{0,0,0}further}%
}}}}
\put(6226,-1486){\makebox(0,0)[b]{\smash{{\SetFigFont{10}{12.0}{\familydefault}{\mddefault}{\updefault}{\color[rgb]{0,0,0}$g'$-rules}%
}}}}
\put(451,-1336){\makebox(0,0)[b]{\smash{{\SetFigFont{10}{12.0}{\familydefault}{\mddefault}{\updefault}{\color[rgb]{0,0,0}further}%
}}}}
\put(451,-1486){\makebox(0,0)[b]{\smash{{\SetFigFont{10}{12.0}{\familydefault}{\mddefault}{\updefault}{\color[rgb]{0,0,0}crossings}%
}}}}
\put(1111,-106){\makebox(0,0)[lb]{\smash{{\SetFigFont{11}{13.2}{\rmdefault}{\mddefault}{\updefault}{\color[rgb]{0,0,0}\grulesB}%
}}}}
\put(1111,494){\makebox(0,0)[lb]{\smash{{\SetFigFont{11}{13.2}{\rmdefault}{\mddefault}{\updefault}{\color[rgb]{0,0,0}\grulesC}%
}}}}
\put(1111,-706){\makebox(0,0)[lb]{\smash{{\SetFigFont{11}{13.2}{\rmdefault}{\mddefault}{\updefault}{\color[rgb]{0,0,0}\grulesA}%
}}}}
\put(5011,-106){\makebox(0,0)[lb]{\smash{{\SetFigFont{11}{13.2}{\rmdefault}{\mddefault}{\updefault}{\color[rgb]{0,0,0}\gprulesB}%
}}}}
\put(5011,494){\makebox(0,0)[lb]{\smash{{\SetFigFont{11}{13.2}{\rmdefault}{\mddefault}{\updefault}{\color[rgb]{0,0,0}\gprulesC}%
}}}}
\put(5011,-706){\makebox(0,0)[lb]{\smash{{\SetFigFont{11}{13.2}{\rmdefault}{\mddefault}{\updefault}{\color[rgb]{0,0,0}\gprulesA}%
}}}}
\put(  1,689){\makebox(0,0)[lb]{\smash{{\SetFigFont{10}{12.0}{\rmdefault}{\mddefault}{\updefault}{\color[rgb]{0,0,0}$\kpp$}%
}}}}
\put(451,689){\makebox(0,0)[lb]{\smash{{\SetFigFont{10}{12.0}{\rmdefault}{\mddefault}{\updefault}{\color[rgb]{0,0,0}$\jpp$}%
}}}}
\put(901,689){\makebox(0,0)[lb]{\smash{{\SetFigFont{10}{12.0}{\rmdefault}{\mddefault}{\updefault}{\color[rgb]{0,0,0}$\ipp$}%
}}}}
\put(  1,-811){\makebox(0,0)[lb]{\smash{{\SetFigFont{10}{12.0}{\rmdefault}{\mddefault}{\updefault}{\color[rgb]{0,0,0}$i$}%
}}}}
\put(451,-811){\makebox(0,0)[lb]{\smash{{\SetFigFont{10}{12.0}{\rmdefault}{\mddefault}{\updefault}{\color[rgb]{0,0,0}$j$}%
}}}}
\put(901,-811){\makebox(0,0)[lb]{\smash{{\SetFigFont{10}{12.0}{\rmdefault}{\mddefault}{\updefault}{\color[rgb]{0,0,0}$k$}%
}}}}
\put(451, 89){\makebox(0,0)[lb]{\smash{{\SetFigFont{10}{12.0}{\rmdefault}{\mddefault}{\updefault}{\color[rgb]{0,0,0}$\ip$}%
}}}}
\put(901,-61){\makebox(0,0)[lb]{\smash{{\SetFigFont{10}{12.0}{\rmdefault}{\mddefault}{\updefault}{\color[rgb]{0,0,0}$\jp$}%
}}}}
\put(451,-286){\makebox(0,0)[lb]{\smash{{\SetFigFont{10}{12.0}{\rmdefault}{\mddefault}{\updefault}{\color[rgb]{0,0,0}$\kp$}%
}}}}
\put(3901,689){\makebox(0,0)[lb]{\smash{{\SetFigFont{10}{12.0}{\rmdefault}{\mddefault}{\updefault}{\color[rgb]{0,0,0}$\kpp$}%
}}}}
\put(4351,689){\makebox(0,0)[lb]{\smash{{\SetFigFont{10}{12.0}{\rmdefault}{\mddefault}{\updefault}{\color[rgb]{0,0,0}$\jpp$}%
}}}}
\put(4801,689){\makebox(0,0)[lb]{\smash{{\SetFigFont{10}{12.0}{\rmdefault}{\mddefault}{\updefault}{\color[rgb]{0,0,0}$\ipp$}%
}}}}
\put(3901,-811){\makebox(0,0)[lb]{\smash{{\SetFigFont{10}{12.0}{\rmdefault}{\mddefault}{\updefault}{\color[rgb]{0,0,0}$i$}%
}}}}
\put(4351,-811){\makebox(0,0)[lb]{\smash{{\SetFigFont{10}{12.0}{\rmdefault}{\mddefault}{\updefault}{\color[rgb]{0,0,0}$j$}%
}}}}
\put(4801,-811){\makebox(0,0)[lb]{\smash{{\SetFigFont{10}{12.0}{\rmdefault}{\mddefault}{\updefault}{\color[rgb]{0,0,0}$k$}%
}}}}
\put(4351,-361){\makebox(0,0)[lb]{\smash{{\SetFigFont{10}{12.0}{\rmdefault}{\mddefault}{\updefault}{\color[rgb]{0,0,0}$\ip$}%
}}}}
\put(3901,-136){\makebox(0,0)[lb]{\smash{{\SetFigFont{10}{12.0}{\rmdefault}{\mddefault}{\updefault}{\color[rgb]{0,0,0}$\jp$}%
}}}}
\put(4351, 89){\makebox(0,0)[lb]{\smash{{\SetFigFont{10}{12.0}{\rmdefault}{\mddefault}{\updefault}{\color[rgb]{0,0,0}$\kp$}%
}}}}
\end{picture}%

%% file: figs/R2cp.pdf_t
\begin{picture}(0,0)%
\includegraphics{figs/R2cp.pdf}%
\end{picture}%
%
%
\setlength{\unitlength}{3947sp}%
\begingroup\makeatletter\ifx\SetFigFont\undefined%
\gdef\SetFigFont#1#2#3#4#5{%
  \reset@font\fontsize{#1}{#2pt}%
  \fontfamily{#3}\fontseries{#4}\fontshape{#5}%
  \selectfont}%
\fi\endgroup%
\begin{picture}(7377,1849)(-14,-950)
\put(6676,-511){\makebox(0,0)[b]{\smash{{\SetFigFont{10}{12.0}{\familydefault}{\mddefault}{\updefault}{\color[rgb]{0,0,0}$\cdots$}%
}}}}
\put(6676,-736){\makebox(0,0)[b]{\smash{{\SetFigFont{10}{12.0}{\familydefault}{\mddefault}{\updefault}{\color[rgb]{0,0,0}further}%
}}}}
\put(6676,-886){\makebox(0,0)[b]{\smash{{\SetFigFont{10}{12.0}{\familydefault}{\mddefault}{\updefault}{\color[rgb]{0,0,0}$g'$-rules}%
}}}}
\put(5251,-511){\makebox(0,0)[b]{\smash{{\SetFigFont{10}{12.0}{\familydefault}{\mddefault}{\updefault}{\color[rgb]{0,0,0}$\cdots$}%
}}}}
\put(5251,-736){\makebox(0,0)[b]{\smash{{\SetFigFont{10}{12.0}{\familydefault}{\mddefault}{\updefault}{\color[rgb]{0,0,0}further}%
}}}}
\put(5251,-886){\makebox(0,0)[b]{\smash{{\SetFigFont{10}{12.0}{\familydefault}{\mddefault}{\updefault}{\color[rgb]{0,0,0}crossings}%
}}}}
\put(5140,239){\makebox(0,0)[b]{\smash{{\SetFigFont{10}{12.0}{\rmdefault}{\mddefault}{\updefault}{\color[rgb]{0,0,0}$\ip$}%
}}}}
\put(5604,239){\makebox(0,0)[b]{\smash{{\SetFigFont{10}{12.0}{\rmdefault}{\mddefault}{\updefault}{\color[rgb]{0,0,0}$\jp$}%
}}}}
\put(1051,539){\makebox(0,0)[lb]{\smash{{\SetFigFont{11}{13.2}{\rmdefault}{\mddefault}{\updefault}{\color[rgb]{0,0,0}\grulesB}%
}}}}
\put(1051,-61){\makebox(0,0)[lb]{\smash{{\SetFigFont{11}{13.2}{\rmdefault}{\mddefault}{\updefault}{\color[rgb]{0,0,0}\grulesA}%
}}}}
\put(2401,-511){\makebox(0,0)[b]{\smash{{\SetFigFont{10}{12.0}{\familydefault}{\mddefault}{\updefault}{\color[rgb]{0,0,0}$\cdots$}%
}}}}
\put(2401,-736){\makebox(0,0)[b]{\smash{{\SetFigFont{10}{12.0}{\familydefault}{\mddefault}{\updefault}{\color[rgb]{0,0,0}further}%
}}}}
\put(2401,-886){\makebox(0,0)[b]{\smash{{\SetFigFont{10}{12.0}{\familydefault}{\mddefault}{\updefault}{\color[rgb]{0,0,0}$g$-rules}%
}}}}
\put(5101,-211){\makebox(0,0)[b]{\smash{{\SetFigFont{10}{12.0}{\rmdefault}{\mddefault}{\updefault}{\color[rgb]{0,0,0}$i$}%
}}}}
\put(5851,-61){\makebox(0,0)[lb]{\smash{{\SetFigFont{11}{13.2}{\rmdefault}{\mddefault}{\updefault}{\color[rgb]{0,0,0}\grulesD}%
}}}}
\put(5851,539){\makebox(0,0)[lb]{\smash{{\SetFigFont{11}{13.2}{\rmdefault}{\mddefault}{\updefault}{\color[rgb]{0,0,0}\grulesC}%
}}}}
\put(5551,764){\makebox(0,0)[b]{\smash{{\SetFigFont{10}{12.0}{\rmdefault}{\mddefault}{\updefault}{\color[rgb]{0,0,0}$j$}%
}}}}
\put(5626,-211){\makebox(0,0)[b]{\smash{{\SetFigFont{10}{12.0}{\rmdefault}{\mddefault}{\updefault}{\color[rgb]{0,0,0}$\jpp$}%
}}}}
\put(5176,764){\makebox(0,0)[b]{\smash{{\SetFigFont{10}{12.0}{\rmdefault}{\mddefault}{\updefault}{\color[rgb]{0,0,0}$\ipp$}%
}}}}
\put(376,-511){\rotatebox{360.0}{\makebox(0,0)[b]{\smash{{\SetFigFont{10}{12.0}{\familydefault}{\mddefault}{\updefault}{\color[rgb]{0,0,0}$\cdots$}%
}}}}}
\put(376,-736){\rotatebox{360.0}{\makebox(0,0)[b]{\smash{{\SetFigFont{10}{12.0}{\familydefault}{\mddefault}{\updefault}{\color[rgb]{0,0,0}further}%
}}}}}
\put(376,-886){\rotatebox{360.0}{\makebox(0,0)[b]{\smash{{\SetFigFont{10}{12.0}{\familydefault}{\mddefault}{\updefault}{\color[rgb]{0,0,0}crossings}%
}}}}}
\put(  1,764){\makebox(0,0)[b]{\smash{{\SetFigFont{10}{12.0}{\rmdefault}{\mddefault}{\updefault}{\color[rgb]{0,0,0}$\ipp$}%
}}}}
\put(676,764){\makebox(0,0)[b]{\smash{{\SetFigFont{10}{12.0}{\rmdefault}{\mddefault}{\updefault}{\color[rgb]{0,0,0}$j$}%
}}}}
\put( 76,-211){\makebox(0,0)[b]{\smash{{\SetFigFont{10}{12.0}{\rmdefault}{\mddefault}{\updefault}{\color[rgb]{0,0,0}$i$}%
}}}}
\put(751,-211){\makebox(0,0)[b]{\smash{{\SetFigFont{10}{12.0}{\rmdefault}{\mddefault}{\updefault}{\color[rgb]{0,0,0}$\jpp$}%
}}}}
\put(729,239){\makebox(0,0)[b]{\smash{{\SetFigFont{10}{12.0}{\rmdefault}{\mddefault}{\updefault}{\color[rgb]{0,0,0}$\ip$}%
}}}}
\put( 21,239){\makebox(0,0)[b]{\smash{{\SetFigFont{10}{12.0}{\rmdefault}{\mddefault}{\updefault}{\color[rgb]{0,0,0}$\jp$}%
}}}}
\end{picture}%

%% file: figs/R2cm.pdf_t
\begin{picture}(0,0)%
\includegraphics{figs/R2cm.pdf}%
\end{picture}%
%
%
\setlength{\unitlength}{3947sp}%
\begingroup\makeatletter\ifx\SetFigFont\undefined%
\gdef\SetFigFont#1#2#3#4#5{%
  \reset@font\fontsize{#1}{#2pt}%
  \fontfamily{#3}\fontseries{#4}\fontshape{#5}%
  \selectfont}%
\fi\endgroup%
\begin{picture}(7077,1170)(-14,-271)
\put(4840,239){\makebox(0,0)[b]{\smash{{\SetFigFont{10}{12.0}{\rmdefault}{\mddefault}{\updefault}{\color[rgb]{0,0,0}$\ip$}%
}}}}
\put(5304,239){\makebox(0,0)[b]{\smash{{\SetFigFont{10}{12.0}{\rmdefault}{\mddefault}{\updefault}{\color[rgb]{0,0,0}$\jp$}%
}}}}
\put(5551,-61){\makebox(0,0)[lb]{\smash{{\SetFigFont{11}{13.2}{\rmdefault}{\mddefault}{\updefault}{\color[rgb]{0,0,0}\grulesD}%
}}}}
\put(5551,539){\makebox(0,0)[lb]{\smash{{\SetFigFont{11}{13.2}{\rmdefault}{\mddefault}{\updefault}{\color[rgb]{0,0,0}\grulesC}%
}}}}
\put(4801,764){\makebox(0,0)[b]{\smash{{\SetFigFont{10}{12.0}{\rmdefault}{\mddefault}{\updefault}{\color[rgb]{0,0,0}$i$}%
}}}}
\put(4876,-211){\makebox(0,0)[b]{\smash{{\SetFigFont{10}{12.0}{\rmdefault}{\mddefault}{\updefault}{\color[rgb]{0,0,0}$\ipp$}%
}}}}
\put(5251,-211){\makebox(0,0)[b]{\smash{{\SetFigFont{10}{12.0}{\rmdefault}{\mddefault}{\updefault}{\color[rgb]{0,0,0}$j$}%
}}}}
\put(5326,764){\makebox(0,0)[b]{\smash{{\SetFigFont{10}{12.0}{\rmdefault}{\mddefault}{\updefault}{\color[rgb]{0,0,0}$\jpp$}%
}}}}
\put(729,239){\makebox(0,0)[b]{\smash{{\SetFigFont{10}{12.0}{\rmdefault}{\mddefault}{\updefault}{\color[rgb]{0,0,0}$\ip$}%
}}}}
\put( 21,239){\makebox(0,0)[b]{\smash{{\SetFigFont{10}{12.0}{\rmdefault}{\mddefault}{\updefault}{\color[rgb]{0,0,0}$\jp$}%
}}}}
\put(1051,539){\makebox(0,0)[lb]{\smash{{\SetFigFont{11}{13.2}{\rmdefault}{\mddefault}{\updefault}{\color[rgb]{0,0,0}\grulesB}%
}}}}
\put(1051,-61){\makebox(0,0)[lb]{\smash{{\SetFigFont{11}{13.2}{\rmdefault}{\mddefault}{\updefault}{\color[rgb]{0,0,0}\grulesA}%
}}}}
\put( 76,764){\makebox(0,0)[b]{\smash{{\SetFigFont{10}{12.0}{\rmdefault}{\mddefault}{\updefault}{\color[rgb]{0,0,0}$i$}%
}}}}
\put(  1,-211){\makebox(0,0)[b]{\smash{{\SetFigFont{10}{12.0}{\rmdefault}{\mddefault}{\updefault}{\color[rgb]{0,0,0}$\ipp$}%
}}}}
\put(676,-211){\makebox(0,0)[b]{\smash{{\SetFigFont{10}{12.0}{\rmdefault}{\mddefault}{\updefault}{\color[rgb]{0,0,0}$j$}%
}}}}
\put(751,764){\makebox(0,0)[b]{\smash{{\SetFigFont{10}{12.0}{\rmdefault}{\mddefault}{\updefault}{\color[rgb]{0,0,0}$\jpp$}%
}}}}
\end{picture}%

%% file: figs/R1s.pdf_t
\begin{picture}(0,0)%
\includegraphics{figs/R1s.pdf}%
\end{picture}%
%
%
\setlength{\unitlength}{3947sp}%
\begingroup\makeatletter\ifx\SetFigFont\undefined%
\gdef\SetFigFont#1#2#3#4#5{%
  \reset@font\fontsize{#1}{#2pt}%
  \fontfamily{#3}\fontseries{#4}\fontshape{#5}%
  \selectfont}%
\fi\endgroup%
\begin{picture}(7307,967)(56,-116)
\put(301,689){\makebox(0,0)[b]{\smash{{\SetFigFont{10}{12.0}{\rmdefault}{\mddefault}{\updefault}{\color[rgb]{0,0,0}$\ipp$}%
}}}}
\put(376,-61){\makebox(0,0)[b]{\smash{{\SetFigFont{10}{12.0}{\rmdefault}{\mddefault}{\updefault}{\color[rgb]{0,0,0}$i$}%
}}}}
\put(226,314){\makebox(0,0)[b]{\smash{{\SetFigFont{10}{12.0}{\rmdefault}{\mddefault}{\updefault}{\color[rgb]{0,0,0}$\ip$}%
}}}}
\put(601,314){\makebox(0,0)[lb]{\smash{{\SetFigFont{11}{13.2}{\rmdefault}{\mddefault}{\updefault}{\color[rgb]{0,0,0}\grulesA}%
}}}}
\put(2926,689){\makebox(0,0)[b]{\smash{{\SetFigFont{10}{12.0}{\rmdefault}{\mddefault}{\updefault}{\color[rgb]{0,0,0}$\ipp$}%
}}}}
\put(2851,-61){\makebox(0,0)[b]{\smash{{\SetFigFont{10}{12.0}{\rmdefault}{\mddefault}{\updefault}{\color[rgb]{0,0,0}$i$}%
}}}}
\put(2899,314){\makebox(0,0)[b]{\smash{{\SetFigFont{10}{12.0}{\rmdefault}{\mddefault}{\updefault}{\color[rgb]{0,0,0}$\ip$}%
}}}}
\put(3001,314){\makebox(0,0)[lb]{\smash{{\SetFigFont{11}{13.2}{\rmdefault}{\mddefault}{\updefault}{\color[rgb]{0,0,0}\grulesB}%
}}}}
\put(5251,314){\makebox(0,0)[b]{\smash{{\SetFigFont{10}{12.0}{\rmdefault}{\mddefault}{\updefault}{\color[rgb]{0,0,0}$\ip$}%
}}}}
\put(5476,314){\makebox(0,0)[lb]{\smash{{\SetFigFont{11}{13.2}{\rmdefault}{\mddefault}{\updefault}{\color[rgb]{0,0,0}\grulesC}%
}}}}
\put(5101,689){\makebox(0,0)[b]{\smash{{\SetFigFont{10}{12.0}{\rmdefault}{\mddefault}{\updefault}{\color[rgb]{0,0,0}$\ipp$}%
}}}}
\put(5026,-61){\makebox(0,0)[b]{\smash{{\SetFigFont{10}{12.0}{\rmdefault}{\mddefault}{\updefault}{\color[rgb]{0,0,0}$i$}%
}}}}
\end{picture}%

%% file: figs/UprightRMoves.pdf_t
\begin{picture}(0,0)%
\includegraphics{figs/UprightRMoves.pdf}%
\end{picture}%
%
%
\setlength{\unitlength}{3947sp}%
\begingroup\makeatletter\ifx\SetFigFont\undefined%
\gdef\SetFigFont#1#2#3#4#5{%
  \reset@font\fontsize{#1}{#2pt}%
  \fontfamily{#3}\fontseries{#4}\fontshape{#5}%
  \selectfont}%
\fi\endgroup%
\begin{picture}(7696,2127)(356,-1276)
\put(1651,239){\makebox(0,0)[b]{\smash{{\SetFigFont{10}{12.0}{\rmdefault}{\mddefault}{\updefault}{\color[rgb]{0,0,0}R1r}%
}}}}
\put(1051,239){\makebox(0,0)[b]{\smash{{\SetFigFont{10}{12.0}{\rmdefault}{\mddefault}{\updefault}{\color[rgb]{0,0,0}R1l}%
}}}}
\put(1426,-1261){\rotatebox{90.0}{\makebox(0,0)[lb]{\smash{{\SetFigFont{8}{9.6}{\rmdefault}{\mddefault}{\updefault}{\color[rgb]{0,0,0}\r{m}}%
}}}}}
\put(1426,-361){\rotatebox{90.0}{\makebox(0,0)[rb]{\smash{{\SetFigFont{8}{9.6}{\rmdefault}{\mddefault}{\updefault}{\color[rgb]{0,0,0}\r{n}}%
}}}}}
\put(2251,-811){\rotatebox{90.0}{\makebox(0,0)[b]{\smash{{\SetFigFont{8}{9.6}{\rmdefault}{\mddefault}{\updefault}{\color[rgb]{0,0,0}\r{m+n}}%
}}}}}
\put(1801,-961){\makebox(0,0)[b]{\smash{{\SetFigFont{10}{12.0}{\rmdefault}{\mddefault}{\updefault}{\color[rgb]{0,0,0}NV}%
}}}}
\put(6451,-961){\makebox(0,0)[b]{\smash{{\SetFigFont{10}{12.0}{\rmdefault}{\mddefault}{\updefault}{\color[rgb]{0,0,0}Sw}%
}}}}
\put(4276,-886){\makebox(0,0)[b]{\smash{{\SetFigFont{10}{12.0}{\rmdefault}{\mddefault}{\updefault}{\color[rgb]{0,0,0}$=$}%
}}}}
\put(3526,-1261){\rotatebox{90.0}{\makebox(0,0)[lb]{\smash{{\SetFigFont{7}{8.4}{\rmdefault}{\mddefault}{\updefault}{\color[rgb]{0,0,0}\r{-1}}%
}}}}}
\put(4051,-1261){\rotatebox{90.0}{\makebox(0,0)[lb]{\smash{{\SetFigFont{7}{8.4}{\rmdefault}{\mddefault}{\updefault}{\color[rgb]{0,0,0}\r{-1}}%
}}}}}
\put(3526,-361){\rotatebox{90.0}{\makebox(0,0)[rb]{\smash{{\SetFigFont{7}{8.4}{\rmdefault}{\mddefault}{\updefault}{\color[rgb]{0,0,0}\r{1}}%
}}}}}
\put(4051,-361){\rotatebox{90.0}{\makebox(0,0)[rb]{\smash{{\SetFigFont{7}{8.4}{\rmdefault}{\mddefault}{\updefault}{\color[rgb]{0,0,0}\r{1}}%
}}}}}
\put(3526,239){\makebox(0,0)[b]{\smash{{\SetFigFont{10}{12.0}{\rmdefault}{\mddefault}{\updefault}{\color[rgb]{0,0,0}R2c$^+$}%
}}}}
\put(5251,239){\makebox(0,0)[b]{\smash{{\SetFigFont{10}{12.0}{\rmdefault}{\mddefault}{\updefault}{\color[rgb]{0,0,0}R2c$^-$}%
}}}}
\put(7126,239){\makebox(0,0)[b]{\smash{{\SetFigFont{10}{12.0}{\rmdefault}{\mddefault}{\updefault}{\color[rgb]{0,0,0}R3b}%
}}}}
\end{picture}%

%% file: figs/R3.pdf_t
\begin{picture}(0,0)%
\includegraphics{figs/R3.pdf}%
\end{picture}%
\setlength{\unitlength}{3947sp}%
\begingroup\makeatletter\ifx\SetFigFont\undefined%
\gdef\SetFigFont#1#2#3#4#5{%
  \reset@font\fontsize{#1}{#2pt}%
  \fontfamily{#3}\fontseries{#4}\fontshape{#5}%
  \selectfont}%
\fi\endgroup%
\begin{picture}(5949,1524)(-86,-673)
\put(2101,614){\makebox(0,0)[b]{\smash{{\SetFigFont{11}{13.2}{\rmdefault}{\mddefault}{\updefault}{\color[rgb]{0,0,0}$D^l$}%
}}}}
\put(5701,614){\makebox(0,0)[b]{\smash{{\SetFigFont{11}{13.2}{\rmdefault}{\mddefault}{\updefault}{\color[rgb]{0,0,0}$D^r$}%
}}}}
\put(5333,-398){\makebox(0,0)[b]{\smash{{\SetFigFont{9}{10.8}{\rmdefault}{\mddefault}{\updefault}{\color[rgb]{0,0,0}\m}%
}}}}
\put(5633,-398){\makebox(0,0)[b]{\smash{{\SetFigFont{9}{10.8}{\rmdefault}{\mddefault}{\updefault}{\color[rgb]{0,0,0}\n}%
}}}}
\put(5483,-23){\makebox(0,0)[b]{\smash{{\SetFigFont{9}{10.8}{\rmdefault}{\mddefault}{\updefault}{\color[rgb]{0,0,0}\s}%
}}}}
\put(1733,-398){\makebox(0,0)[b]{\smash{{\SetFigFont{9}{10.8}{\rmdefault}{\mddefault}{\updefault}{\color[rgb]{0,0,0}\m}%
}}}}
\put(2033,-398){\makebox(0,0)[b]{\smash{{\SetFigFont{9}{10.8}{\rmdefault}{\mddefault}{\updefault}{\color[rgb]{0,0,0}\n}%
}}}}
\put(1883,-23){\makebox(0,0)[b]{\smash{{\SetFigFont{9}{10.8}{\rmdefault}{\mddefault}{\updefault}{\color[rgb]{0,0,0}\s}%
}}}}
\put(526,-361){\makebox(0,0)[rb]{\smash{{\SetFigFont{9}{10.8}{\rmdefault}{\mddefault}{\updefault}{\color[rgb]{0,0,0}$c^l_1$}%
}}}}
\put(226, 14){\makebox(0,0)[rb]{\smash{{\SetFigFont{9}{10.8}{\rmdefault}{\mddefault}{\updefault}{\color[rgb]{0,0,0}$c^l_2$}%
}}}}
\put(526,389){\makebox(0,0)[rb]{\smash{{\SetFigFont{9}{10.8}{\rmdefault}{\mddefault}{\updefault}{\color[rgb]{0,0,0}$c^l_3$}%
}}}}
\put(1801,-136){\makebox(0,0)[rb]{\smash{{\SetFigFont{9}{10.8}{\rmdefault}{\mddefault}{\updefault}{\color[rgb]{0,0,0}$c_y$}%
}}}}
\put(233,-548){\makebox(0,0)[b]{\smash{{\SetFigFont{9}{10.8}{\rmdefault}{\mddefault}{\updefault}{\color[rgb]{0,0,0}\i}%
}}}}
\put(533,-548){\makebox(0,0)[b]{\smash{{\SetFigFont{9}{10.8}{\rmdefault}{\mddefault}{\updefault}{\color[rgb]{0,0,0}\j}%
}}}}
\put(833,-548){\makebox(0,0)[b]{\smash{{\SetFigFont{9}{10.8}{\rmdefault}{\mddefault}{\updefault}{\color[rgb]{0,0,0}\k}%
}}}}
\put(158,652){\makebox(0,0)[b]{\smash{{\SetFigFont{9}{10.8}{\rmdefault}{\mddefault}{\updefault}{\color[rgb]{0,0,0}\kpp}%
}}}}
\put(458,652){\makebox(0,0)[b]{\smash{{\SetFigFont{9}{10.8}{\rmdefault}{\mddefault}{\updefault}{\color[rgb]{0,0,0}\jpp}%
}}}}
\put(758,652){\makebox(0,0)[b]{\smash{{\SetFigFont{9}{10.8}{\rmdefault}{\mddefault}{\updefault}{\color[rgb]{0,0,0}\ipp}%
}}}}
\put(3833,-548){\makebox(0,0)[b]{\smash{{\SetFigFont{9}{10.8}{\rmdefault}{\mddefault}{\updefault}{\color[rgb]{0,0,0}\i}%
}}}}
\put(4133,-548){\makebox(0,0)[b]{\smash{{\SetFigFont{9}{10.8}{\rmdefault}{\mddefault}{\updefault}{\color[rgb]{0,0,0}\j}%
}}}}
\put(4433,-548){\makebox(0,0)[b]{\smash{{\SetFigFont{9}{10.8}{\rmdefault}{\mddefault}{\updefault}{\color[rgb]{0,0,0}\k}%
}}}}
\put(3758,652){\makebox(0,0)[b]{\smash{{\SetFigFont{9}{10.8}{\rmdefault}{\mddefault}{\updefault}{\color[rgb]{0,0,0}\kpp}%
}}}}
\put(4058,652){\makebox(0,0)[b]{\smash{{\SetFigFont{9}{10.8}{\rmdefault}{\mddefault}{\updefault}{\color[rgb]{0,0,0}\jpp}%
}}}}
\put(4358,652){\makebox(0,0)[b]{\smash{{\SetFigFont{9}{10.8}{\rmdefault}{\mddefault}{\updefault}{\color[rgb]{0,0,0}\ipp}%
}}}}
\put(5401,-136){\makebox(0,0)[rb]{\smash{{\SetFigFont{9}{10.8}{\rmdefault}{\mddefault}{\updefault}{\color[rgb]{0,0,0}$c_y$}%
}}}}
\put(3826,389){\makebox(0,0)[rb]{\smash{{\SetFigFont{9}{10.8}{\rmdefault}{\mddefault}{\updefault}{\color[rgb]{0,0,0}$c^r_3$}%
}}}}
\put(3826,-361){\makebox(0,0)[rb]{\smash{{\SetFigFont{9}{10.8}{\rmdefault}{\mddefault}{\updefault}{\color[rgb]{0,0,0}$c^r_1$}%
}}}}
\put(793,-23){\makebox(0,0)[lb]{\smash{{\SetFigFont{9}{10.8}{\rmdefault}{\mddefault}{\updefault}{\color[rgb]{0,0,0}\jp}%
}}}}
\put(500,202){\makebox(0,0)[lb]{\smash{{\SetFigFont{9}{10.8}{\rmdefault}{\mddefault}{\updefault}{\color[rgb]{0,0,0}\ip}%
}}}}
\put(480,-173){\makebox(0,0)[lb]{\smash{{\SetFigFont{9}{10.8}{\rmdefault}{\mddefault}{\updefault}{\color[rgb]{0,0,0}\kp}%
}}}}
\put(4126, 14){\makebox(0,0)[rb]{\smash{{\SetFigFont{9}{10.8}{\rmdefault}{\mddefault}{\updefault}{\color[rgb]{0,0,0}$c^r_2$}%
}}}}
\put(3673,-183){\makebox(0,0)[b]{\smash{{\SetFigFont{9}{10.8}{\rmdefault}{\mddefault}{\updefault}{\color[rgb]{0,0,0}\jp}%
}}}}
\put(4160,-213){\makebox(0,0)[b]{\smash{{\SetFigFont{9}{10.8}{\rmdefault}{\mddefault}{\updefault}{\color[rgb]{0,0,0}\ip}%
}}}}
\put(4173,242){\makebox(0,0)[b]{\smash{{\SetFigFont{9}{10.8}{\rmdefault}{\mddefault}{\updefault}{\color[rgb]{0,0,0}\kp}%
}}}}
\end{picture}%

%% file: Invariance-R3.tex
First, we implement the Kronecker $\delta$-function, the $g$-rules for a crossing $(s,i,j)$, and the $g$-rules for a list of crossings $X$:

\nbpdfInput{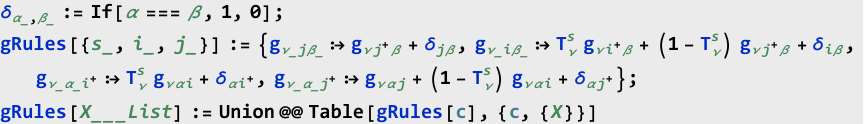}
We then let \verb$Xl$ be the three crossings in the left-hand-side of the R3b move, as in Figure~\ref{fig:R3}, we let \verb$Al$ be the $A^l$ term of~\eqref{eq:ABC}, and we let \verb$lhs$ be the result of applying the $g$-rules for the crossings in \verb$Xl$ to \verb$Al$. We print only a ``\verb$Short$'' version of \verb$lhs$ because the full thing would cover about 2.5 pages: 

\nbpdfInput{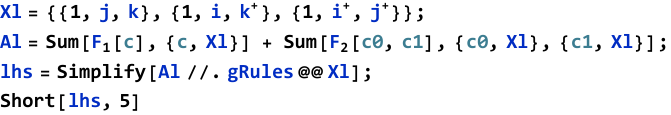}
\nbpdfOutput{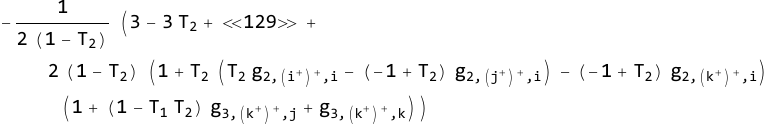}
We do the same for $A^r$, except this time, without printing at all:

\nbpdfInput{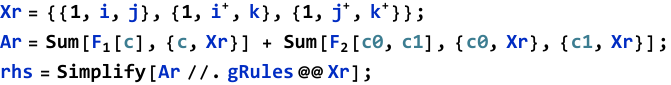}
We then compare \verb$lhs$ with \verb$rhs$. The output, \verb$True$, tells us that we have proven~\eqref{eq:R3A}:
{\def\nbpdfPostInput{\hfill}

\nbpdfInput{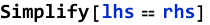}
\nbpdfOutput{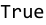}
We show that $B^l=B^r$ by following exactly the same procedure. Note that we ignore the summation over $c_y$ and instead treat $c_y$ as a fixed crossing $(s,m,n)$. If an equality is proven for every fixed $c_y$, it is of course also proven for the sum over $c_y\in Y$.

\nbpdfInput{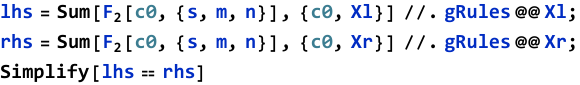}
\nbpdfOutput{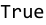}
Similarly we prove that $C^l=C^r$, and this concludes the proof of Proposition~\ref{prop:R3}.

\nbpdfInput{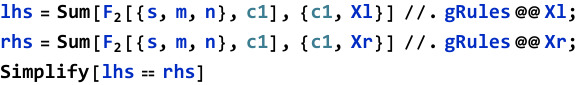}
\nbpdfOutput{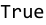}
}
\vskip -7mm
\ \endpar{\ref{prop:R3}}

\begin{remark} \label{rem:E} The computations above were carried out for generic $g_{\nu\alpha\beta}$ and for a generic $c_y=(s,m,n)$; namely, without specifying the knot diagrams in full, and hence without assigning specific values to $g_{\nu\alpha\beta}$, and without specifying $m$ and $n$. Under these conditions the three parts of~\eqref{eq:ABC} cannot mix (namely, terms from, say, $A^h$ cannot cancel terms in $B^h$ or $C^h$), and so it would have been enough to show that $E^l=E^r$, where $E^h$ combines $A^h$ and $B^h$ and $C^h$ (and a few harmless further terms) by adding $c_y$ to the summation corresponding to $A^h$:
\[
  E^h = \sum_{c\in\{c^h_1,c^h_2,c^h_3,c_y\}} F^h_1(c)
    + \sum_{c_0,c_1\in\{c^h_1,c^h_2,c^h_3,c_y\}} F^h_2(c_0,c_1).
\]
But that's a simpler computation:

\nbpdfInput{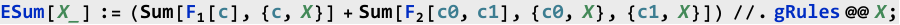}
{\def\nbpdfPostInput{\hfill}

\nbpdfInput{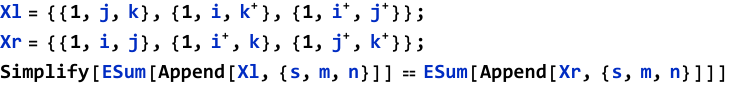}
\nbpdfOutput{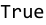}
}\end{remark}\endpar{\ref{rem:E}}

%% file: Invariance-R2c.tex
\noindent{\em Proof.} For R2c$^+$ we follow the same logic as in the proof of Proposition~\ref{prop:R3}, as simplified by Remark~\ref{rem:E}. We start with the figure that replaces Figure~\ref{fig:R3} (note the null vertices in $D^r$ and their minimal effect as in Lemma~\ref{lem:NullVertices} and Remark~\ref{rem:F4Null}):
\[ {
  \def\i{{$i$}} \def\j{{$j$}} \def\k{{$k$}} \def\m{{$m$}} \def\n{{$n$}} \def\s{{$s$}}
  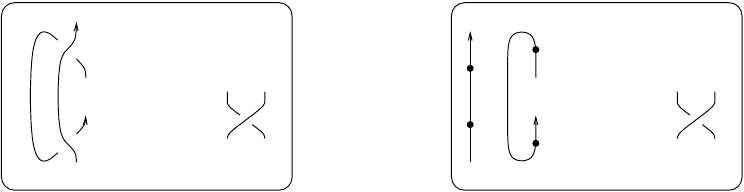 
} \]

As in Remark~\ref{rem:E}, we let $E^l$ and $E^r$ be the sums corresponding to the diagrams $D^l$ and $D^r$ above:
\[
  E^l = \sum_{\makebox[0pt]{$\scriptstyle c\in\{c^l_1,c^l_2,c_y\}$}} F^l_1(c)
    \ +\ \sum_{\makebox[0pt]{$\scriptstyle c_0,c_1\in\{c^l_1,c^l_2,c_y\}$}} F^l_2(c_0,c_1)
    \ +\ F^l_3(\jp)|_{\varphi_{\jp}=1},
  \quad
    E^r = F^r_1(c_y) \ +\ F^r_2(c_y,c_y) \ +\ F^r_3(\jp)|_{\varphi_{\jp}=1}.
\]
We need to show that $E^l=E^r$ after all relevant $g$-rules are applied to both sides.

To compute these $E$ sums we first have to extend the \verb$ESum$ routine to accept also a list $R$ of pairs $(\varphi,k)$ of the form (rotation number, edge label):

\nbpdfInput{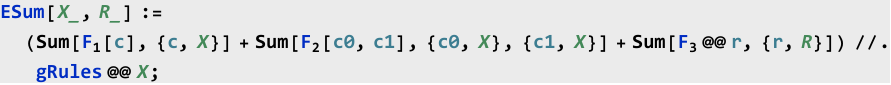}
We then compute $E^l$ (and apply the relevant $g$-rules) by calling \verb$ESum$ with crossings $(-1,i,\jp)$, $(1,\ip,j)$, and $(s,m,n)$, and a rotation number of $1$ on edge $\jp$:

\nbpdfInput{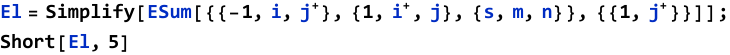}
\nbpdfOutput{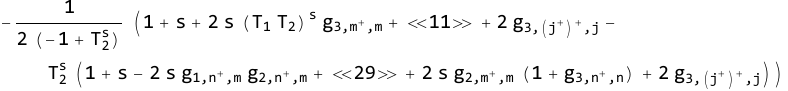}
The computation of $E^r$ is simpler, as it only involves the generic $(s,m,n)$ and the rotation $(1,\jp)$. We implement the $g$-rules for null vertices as in Equations~\eqref{eq:NullCarRules} and~\eqref{eq:NullCounterRules}, compute $E^r$, and then compare $E^l$ with $E^r$ to conclude the invariance under R2c$^+$:

\nbpdfInput{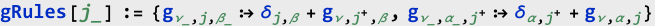}
\nbpdfInput{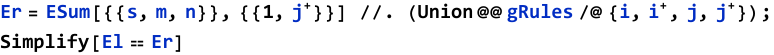}
\nbpdfOutput{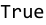}
For R2c$^-$ we allow ourselves to be even more concise:

\parpic[r]{ \def\i{{$i$}} \def\j{{$j$}}  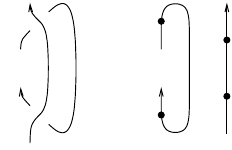 }
\picskip{1}

\nbpdfInput{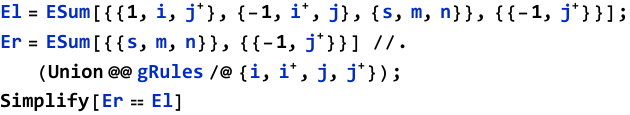}
\nbpdfOutput{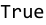}
\endpar{\ref{prop:R2c}}

%% file: figs/R2cE.pdf_t
\begin{picture}(0,0)%
\includegraphics{figs/R2cE.pdf}%
\end{picture}%
\setlength{\unitlength}{3947sp}%
\begingroup\makeatletter\ifx\SetFigFont\undefined%
\gdef\SetFigFont#1#2#3#4#5{%
  \reset@font\fontsize{#1}{#2pt}%
  \fontfamily{#3}\fontseries{#4}\fontshape{#5}%
  \selectfont}%
\fi\endgroup%
\begin{picture}(5949,1524)(-11,-673)
\put(751,239){\makebox(0,0)[b]{\smash{{\SetFigFont{10}{12.0}{\rmdefault}{\mddefault}{\updefault}{\color[rgb]{0,0,0}$j$}%
}}}}
\put(751,539){\makebox(0,0)[b]{\smash{{\SetFigFont{10}{12.0}{\rmdefault}{\mddefault}{\updefault}{\color[rgb]{0,0,0}$\ipp$}%
}}}}
\put(676,-436){\makebox(0,0)[b]{\smash{{\SetFigFont{10}{12.0}{\rmdefault}{\mddefault}{\updefault}{\color[rgb]{0,0,0}$i$}%
}}}}
\put(544, 14){\makebox(0,0)[b]{\smash{{\SetFigFont{10}{12.0}{\rmdefault}{\mddefault}{\updefault}{\color[rgb]{0,0,0}$\ip$}%
}}}}
\put(151, 14){\makebox(0,0)[b]{\smash{{\SetFigFont{10}{12.0}{\rmdefault}{\mddefault}{\updefault}{\color[rgb]{0,0,0}$\jp$}%
}}}}
\put(826,-136){\makebox(0,0)[b]{\smash{{\SetFigFont{10}{12.0}{\rmdefault}{\mddefault}{\updefault}{\color[rgb]{0,0,0}$\jpp$}%
}}}}
\put(2176,614){\makebox(0,0)[b]{\smash{{\SetFigFont{11}{13.2}{\rmdefault}{\mddefault}{\updefault}{\color[rgb]{0,0,0}$D^l$}%
}}}}
\put(5776,614){\makebox(0,0)[b]{\smash{{\SetFigFont{11}{13.2}{\rmdefault}{\mddefault}{\updefault}{\color[rgb]{0,0,0}$D^r$}%
}}}}
\put(3826,-436){\makebox(0,0)[b]{\smash{{\SetFigFont{10}{12.0}{\rmdefault}{\mddefault}{\updefault}{\color[rgb]{0,0,0}$i$}%
}}}}
\put(4351,239){\makebox(0,0)[b]{\smash{{\SetFigFont{10}{12.0}{\rmdefault}{\mddefault}{\updefault}{\color[rgb]{0,0,0}$j$}%
}}}}
\put(4426,-136){\makebox(0,0)[b]{\smash{{\SetFigFont{10}{12.0}{\rmdefault}{\mddefault}{\updefault}{\color[rgb]{0,0,0}$\jpp$}%
}}}}
\put(3901,539){\makebox(0,0)[b]{\smash{{\SetFigFont{10}{12.0}{\rmdefault}{\mddefault}{\updefault}{\color[rgb]{0,0,0}$\ipp$}%
}}}}
\put(4175, 14){\makebox(0,0)[b]{\smash{{\SetFigFont{10}{12.0}{\rmdefault}{\mddefault}{\updefault}{\color[rgb]{0,0,0}$\jp$}%
}}}}
\put(3865, 14){\makebox(0,0)[b]{\smash{{\SetFigFont{10}{12.0}{\rmdefault}{\mddefault}{\updefault}{\color[rgb]{0,0,0}$\ip$}%
}}}}
\put(5408,-398){\makebox(0,0)[b]{\smash{{\SetFigFont{9}{10.8}{\rmdefault}{\mddefault}{\updefault}{\color[rgb]{0,0,0}\m}%
}}}}
\put(5708,-398){\makebox(0,0)[b]{\smash{{\SetFigFont{9}{10.8}{\rmdefault}{\mddefault}{\updefault}{\color[rgb]{0,0,0}\n}%
}}}}
\put(5558,-23){\makebox(0,0)[b]{\smash{{\SetFigFont{9}{10.8}{\rmdefault}{\mddefault}{\updefault}{\color[rgb]{0,0,0}\s}%
}}}}
\put(1808,-398){\makebox(0,0)[b]{\smash{{\SetFigFont{9}{10.8}{\rmdefault}{\mddefault}{\updefault}{\color[rgb]{0,0,0}\m}%
}}}}
\put(2108,-398){\makebox(0,0)[b]{\smash{{\SetFigFont{9}{10.8}{\rmdefault}{\mddefault}{\updefault}{\color[rgb]{0,0,0}\n}%
}}}}
\put(1958,-23){\makebox(0,0)[b]{\smash{{\SetFigFont{9}{10.8}{\rmdefault}{\mddefault}{\updefault}{\color[rgb]{0,0,0}\s}%
}}}}
\put(1876,-136){\makebox(0,0)[rb]{\smash{{\SetFigFont{9}{10.8}{\rmdefault}{\mddefault}{\updefault}{\color[rgb]{0,0,0}$c_y$}%
}}}}
\put(5476,-136){\makebox(0,0)[rb]{\smash{{\SetFigFont{9}{10.8}{\rmdefault}{\mddefault}{\updefault}{\color[rgb]{0,0,0}$c_y$}%
}}}}
\put(451,-286){\makebox(0,0)[rb]{\smash{{\SetFigFont{9}{10.8}{\rmdefault}{\mddefault}{\updefault}{\color[rgb]{0,0,0}$c^l_1$}%
}}}}
\put(451,389){\makebox(0,0)[rb]{\smash{{\SetFigFont{9}{10.8}{\rmdefault}{\mddefault}{\updefault}{\color[rgb]{0,0,0}$c^l_2$}%
}}}}
\end{picture}%

%% file: figs/uR2cm.pdf_t
\begin{picture}(0,0)%
\includegraphics{figs/uR2cm.pdf}%
\end{picture}%
\setlength{\unitlength}{3947sp}%
\begingroup\makeatletter\ifx\SetFigFont\undefined%
\gdef\SetFigFont#1#2#3#4#5{%
  \reset@font\fontsize{#1}{#2pt}%
  \fontfamily{#3}\fontseries{#4}\fontshape{#5}%
  \selectfont}%
\fi\endgroup%
\begin{picture}(1847,1192)(211,-416)
\put(697,164){\makebox(0,0)[b]{\smash{{\SetFigFont{10}{12.0}{\rmdefault}{\mddefault}{\updefault}{\color[rgb]{0,0,0}$\ip$}%
}}}}
\put(922,164){\makebox(0,0)[b]{\smash{{\SetFigFont{10}{12.0}{\rmdefault}{\mddefault}{\updefault}{\color[rgb]{0,0,0}$\jp$}%
}}}}
\put(376,-361){\makebox(0,0)[b]{\smash{{\SetFigFont{10}{12.0}{\rmdefault}{\mddefault}{\updefault}{\color[rgb]{0,0,0}$i$}%
}}}}
\put(301,389){\makebox(0,0)[b]{\smash{{\SetFigFont{10}{12.0}{\rmdefault}{\mddefault}{\updefault}{\color[rgb]{0,0,0}$j$}%
}}}}
\put(301,614){\makebox(0,0)[b]{\smash{{\SetFigFont{10}{12.0}{\rmdefault}{\mddefault}{\updefault}{\color[rgb]{0,0,0}$\ipp$}%
}}}}
\put(226,-61){\makebox(0,0)[b]{\smash{{\SetFigFont{10}{12.0}{\rmdefault}{\mddefault}{\updefault}{\color[rgb]{0,0,0}$\jpp$}%
}}}}
\put(1426,389){\makebox(0,0)[b]{\smash{{\SetFigFont{10}{12.0}{\rmdefault}{\mddefault}{\updefault}{\color[rgb]{0,0,0}$j$}%
}}}}
\put(1651,164){\makebox(0,0)[b]{\smash{{\SetFigFont{10}{12.0}{\rmdefault}{\mddefault}{\updefault}{\color[rgb]{0,0,0}$\jp$}%
}}}}
\put(1951,164){\makebox(0,0)[b]{\smash{{\SetFigFont{10}{12.0}{\rmdefault}{\mddefault}{\updefault}{\color[rgb]{0,0,0}$\ip$}%
}}}}
\put(1951,614){\makebox(0,0)[b]{\smash{{\SetFigFont{10}{12.0}{\rmdefault}{\mddefault}{\updefault}{\color[rgb]{0,0,0}$\ipp$}%
}}}}
\put(1351,-61){\makebox(0,0)[b]{\smash{{\SetFigFont{10}{12.0}{\rmdefault}{\mddefault}{\updefault}{\color[rgb]{0,0,0}$\jpp$}%
}}}}
\put(1951,-286){\makebox(0,0)[b]{\smash{{\SetFigFont{10}{12.0}{\rmdefault}{\mddefault}{\updefault}{\color[rgb]{0,0,0}$i$}%
}}}}
\end{picture}%

%% file: figs/R1sE.pdf_t
\begin{picture}(0,0)%
\includegraphics{figs/R1sE.pdf}%
\end{picture}%
%
%
\setlength{\unitlength}{3947sp}%
\begingroup\makeatletter\ifx\SetFigFont\undefined%
\gdef\SetFigFont#1#2#3#4#5{%
  \reset@font\fontsize{#1}{#2pt}%
  \fontfamily{#3}\fontseries{#4}\fontshape{#5}%
  \selectfont}%
\fi\endgroup%
\begin{picture}(1690,967)(56,-116)
\put(301,689){\makebox(0,0)[b]{\smash{{\SetFigFont{10}{12.0}{\rmdefault}{\mddefault}{\updefault}{\color[rgb]{0,0,0}$\ipp$}%
}}}}
\put(376,-61){\makebox(0,0)[b]{\smash{{\SetFigFont{10}{12.0}{\rmdefault}{\mddefault}{\updefault}{\color[rgb]{0,0,0}$i$}%
}}}}
\put(226,314){\makebox(0,0)[b]{\smash{{\SetFigFont{10}{12.0}{\rmdefault}{\mddefault}{\updefault}{\color[rgb]{0,0,0}$\ip$}%
}}}}
\put(976,689){\makebox(0,0)[b]{\smash{{\SetFigFont{10}{12.0}{\rmdefault}{\mddefault}{\updefault}{\color[rgb]{0,0,0}$\ipp$}%
}}}}
\put(901,-61){\makebox(0,0)[b]{\smash{{\SetFigFont{10}{12.0}{\rmdefault}{\mddefault}{\updefault}{\color[rgb]{0,0,0}$i$}%
}}}}
\put(949,314){\makebox(0,0)[b]{\smash{{\SetFigFont{10}{12.0}{\rmdefault}{\mddefault}{\updefault}{\color[rgb]{0,0,0}$\ip$}%
}}}}
\put(1651,314){\makebox(0,0)[b]{\smash{{\SetFigFont{10}{12.0}{\rmdefault}{\mddefault}{\updefault}{\color[rgb]{0,0,0}$\ip$}%
}}}}
\put(1501,689){\makebox(0,0)[b]{\smash{{\SetFigFont{10}{12.0}{\rmdefault}{\mddefault}{\updefault}{\color[rgb]{0,0,0}$\ipp$}%
}}}}
\put(1426,-61){\makebox(0,0)[b]{\smash{{\SetFigFont{10}{12.0}{\rmdefault}{\mddefault}{\updefault}{\color[rgb]{0,0,0}$i$}%
}}}}
\end{picture}%

%% file: Invariance-R1s.tex
\noindent{\em Proof.} We aim to use the same approach and conventions as in the previous two proofs but hit a minor snag. The $g$-rules for R1l include
\[
  g_{\ip\beta}  = \delta_{\ip\beta} + Tg_{\ipp,\beta}+(1-T)g_{\ip,\beta}
  \qquad\text{and}\qquad
  g_{\alpha,\ip} = g_{\alpha i} + (1-T)g_{\alpha\ip} + \delta_{\alpha,\ip},
\]
and if these are implemented as simple left to right replacement rules, they lead to infinite recursion. Fortunately, these rules can be rewritten in the form
\[
  g_{\ip\beta}  = T^{-1}\delta_{\ip\beta} +g_{\ipp,\beta}
  \qquad\text{and}\qquad
  g_{\alpha,\ip} = T^{-1}g_{\alpha i} + T^{-1}\delta_{\alpha,\ip},
\]
which makes perfectly valid replacement rules. We thus redefine:

\nbpdfInput{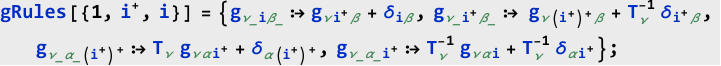}
The same issue does not arise for R1r (!), and thus the following lines conclude the proof:
{\def\nbpdfPostInput{\hfill}

\nbpdfInput{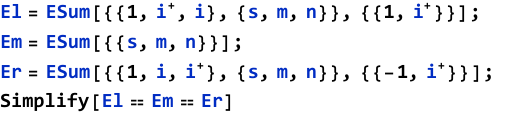}
\nbpdfOutput{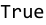}
}\endpar{\ref{prop:R1s}}

%% file: figs/Sw.pdf_t
\begin{picture}(0,0)%
\includegraphics{figs/Sw.pdf}%
\end{picture}%
%
%
\setlength{\unitlength}{3947sp}%
\begingroup\makeatletter\ifx\SetFigFont\undefined%
\gdef\SetFigFont#1#2#3#4#5{%
  \reset@font\fontsize{#1}{#2pt}%
  \fontfamily{#3}\fontseries{#4}\fontshape{#5}%
  \selectfont}%
\fi\endgroup%
\begin{picture}(2352,774)(211,2)
\put(751, 89){\makebox(0,0)[b]{\smash{{\SetFigFont{10}{12.0}{\rmdefault}{\mddefault}{\updefault}{\color[rgb]{0,0,0}$j$}%
}}}}
\put(226,539){\makebox(0,0)[b]{\smash{{\SetFigFont{10}{12.0}{\rmdefault}{\mddefault}{\updefault}{\color[rgb]{0,0,0}$\jp$}%
}}}}
\put(826,539){\makebox(0,0)[b]{\smash{{\SetFigFont{10}{12.0}{\rmdefault}{\mddefault}{\updefault}{\color[rgb]{0,0,0}$\ip$}%
}}}}
\put(301, 89){\makebox(0,0)[b]{\smash{{\SetFigFont{10}{12.0}{\rmdefault}{\mddefault}{\updefault}{\color[rgb]{0,0,0}$i$}%
}}}}
\put(2251,539){\makebox(0,0)[lb]{\smash{{\SetFigFont{7}{8.4}{\rmdefault}{\mddefault}{\updefault}{\color[rgb]{0,0,0}\r{1}}%
}}}}
\put(1651,539){\makebox(0,0)[rb]{\smash{{\SetFigFont{7}{8.4}{\rmdefault}{\mddefault}{\updefault}{\color[rgb]{0,0,0}\r{1}}%
}}}}
\put(1801, 89){\makebox(0,0)[b]{\smash{{\SetFigFont{10}{12.0}{\rmdefault}{\mddefault}{\updefault}{\color[rgb]{0,0,0}$i$}%
}}}}
\put(2101, 89){\makebox(0,0)[b]{\smash{{\SetFigFont{10}{12.0}{\rmdefault}{\mddefault}{\updefault}{\color[rgb]{0,0,0}$j$}%
}}}}
\put(1852,614){\makebox(0,0)[b]{\smash{{\SetFigFont{10}{12.0}{\rmdefault}{\mddefault}{\updefault}{\color[rgb]{0,0,0}$\jp$}%
}}}}
\put(2064,614){\makebox(0,0)[b]{\smash{{\SetFigFont{10}{12.0}{\rmdefault}{\mddefault}{\updefault}{\color[rgb]{0,0,0}$\ip$}%
}}}}
\put(1651,164){\makebox(0,0)[rb]{\smash{{\SetFigFont{7}{8.4}{\rmdefault}{\mddefault}{\updefault}{\color[rgb]{0,0,0}\r{-1}}%
}}}}
\put(2251,164){\makebox(0,0)[lb]{\smash{{\SetFigFont{7}{8.4}{\rmdefault}{\mddefault}{\updefault}{\color[rgb]{0,0,0}\r{-1}}%
}}}}
\end{picture}%

%% file: Invariance-Sw.tex
\noindent{\em Proof.} This one is routine:

\nbpdfInput{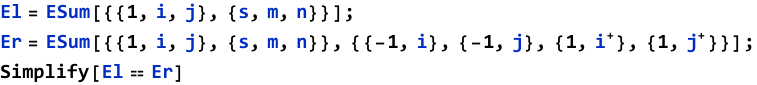}
\nbpdfOutput{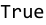}
\endpar{\ref{prop:Sw}}

%% file: table.tex
{\def\s{$\sim$} \begin{center}\begin{tabular}{|r|c|c|c|c|c|c|c|}
\hline
\tl{n}&		$n$&				$\leq 10$&	$\leq 11$&	$\leq 12$&	$\leq 13$&	$\leq 14$&	$\leq 15$ \\ \hline
\tl{Ks}&	knots&				249&		801&		2,977&		12,965&		59,937&		313,230 \\ \hline
\tl{D}&		$\Delta$&			(38)&		(250)&		(1,204)&	(7,326)&	(39,741)&	(236,326) \\ \hline
\tl{s}&		$\sigma_{LT}$&			(108)&		(356)&		(1,525)&	(7,736)&	(40,101)&	(230,592) \\ \hline
\tl{J}&		$J$&				(7)&		(70)&		(482)&		(3,434)&	(21,250)&	(138,591) \\ \hline
\tl{Kh}&	$\Kh$&				(6)&		(65)&		(452)&		(3,226)&	(19,754)&	(127,261) \\ \hline
\tl{H}&		$H$&				(2)&		(31)&		(222)&		(1,839)&	(11,251)&	(73,892) \\ \hline
\tl{V}&		$\Vol$&				(\s6)&		(\s25)&		(\s113)&	(\s1,012)&	(\s6,353)&	(\s43,607) \\ \hline
\tl{KhHV}&	$(\Kh,H,\Vol)$&			(\s0)&		(\s14)&		(\s84)&		(\s911)&	(\s5,917)&	(\s41,434) \\ \hline
\tl{r1}&	$(\Delta,\rho_1)$&		(0)&		(14)&		(95)&		(959)&		(6,253)&	(42,914) \\ \hline
\tl{r12}&	$(\Delta,\rho_1,\rho_2)$&	(0)&		(14)&		(84)&		(911)&		(5,926)&	(41,469) \\ \hline
\tl{r12KhHV}&	$(\rho_1,\rho_2,\Kh,H,\Vol)$&	(0)&		(\s14)&		(\s84)&		(\s911)&	(\s5,916)&	(\s41,432) \\ \hline
\rowcolor{yellow}
\tl{Th}&	$\Theta$&			(0)&		(3)&		(19)&		(194)&		(1,118)&	(6,758) \\ \hline
\tl{Thr2}&	$(\Theta,\rho_2)$&		(0)&		(3)&		(10)&		(169)&		(982)&		(6,341) \\ \hline
\tl{Ths}&	$(\Theta,\sigma_{LT})$&		(0)&		(3)&		(19)&		(194)&		(1,118)&	(6,758) \\ \hline
\tl{ThKh}&	$(\Theta,\Kh)$&			(0)&		(3)&		(18)&		(185)&		(1,062)&	(6,555) \\ \hline
\tl{ThH}&	$(\Theta,H)$&			(0)&		(3)&		(18)&		(185)&		(1,064)&	(6,563) \\ \hline
\tl{ThV}&	$(\Theta,\Vol)$&		(0)&		(\s3)&		(\s10)&		(\s169)&	(\s973)&	(\s6,308) \\ \hline
\tl{Thr2KhHV}&	$(\Theta,\rho_2,\Kh,H,\Vol)$&	(0)&		(\s3)& 		(\s10)&		(\s169)&	(\s972)&	(\s6,304) \\ \hline
\end{tabular}\end{center}}
\caption{
  The separation powers of some knot invariants and combinations of knot
  invariants (in lines \ref{tl:D}--\ref{tl:Thr2KhHV}, smaller numbers are
  better). The data in this table was assembled by \cite[Stats.nb]{Self}.
} \label{tab:Strong}
\end{table}

%% file: figs/GST48-Marked.pdf_t
\begin{picture}(0,0)%
\includegraphics{figs/GST48-Marked.pdf}%
\end{picture}%
%
%
\setlength{\unitlength}{3947sp}%
\begingroup\makeatletter\ifx\SetFigFont\undefined%
\gdef\SetFigFont#1#2#3#4#5{%
  \reset@font\fontsize{#1}{#2pt}%
  \fontfamily{#3}\fontseries{#4}\fontshape{#5}%
  \selectfont}%
\fi\endgroup%
\begin{picture}(9127,4852)(87,-3983)
\end{picture}%

%% file: GST48.tex
{\def\nbpdfPostInput{\hfill}
  {\def\nbpdfPreOutput{\noindent\hspace{-7mm}}

\nbpdfInput{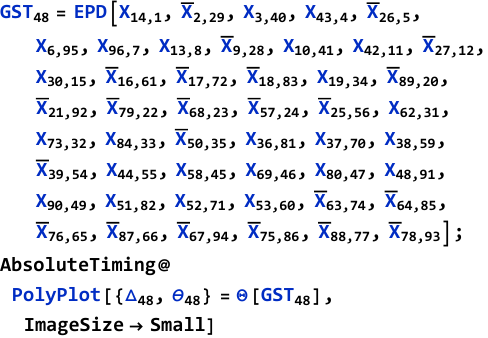}
\nbpdfOutput{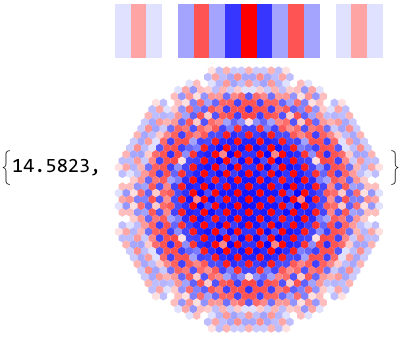}
  }

\nbpdfInput{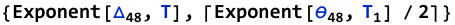}
\nbpdfOutput{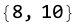}
}

%% file: figs/KS.pdf_t
\begin{picture}(0,0)%
\includegraphics{figs/KS.pdf}%
\end{picture}%
%
%
\setlength{\unitlength}{2368sp}%
\begingroup\makeatletter\ifx\SetFigFont\undefined%
\gdef\SetFigFont#1#2#3#4#5{%
  \reset@font\fontsize{#1}{#2pt}%
  \fontfamily{#3}\fontseries{#4}\fontshape{#5}%
  \selectfont}%
\fi\endgroup%
\begin{picture}(3984,3197)(15,-2028)
\put(2999,-1053){\makebox(0,0)[b]{\smash{{\SetFigFont{6}{7.2}{\rmdefault}{\mddefault}{\updefault}{\color[rgb]{0,0,0}1}%
}}}}
\put(2179,-903){\makebox(0,0)[b]{\smash{{\SetFigFont{6}{7.2}{\rmdefault}{\mddefault}{\updefault}{\color[rgb]{0,0,0}2}%
}}}}
\put(1899,-387){\makebox(0,0)[b]{\smash{{\SetFigFont{6}{7.2}{\rmdefault}{\mddefault}{\updefault}{\color[rgb]{0,0,0}3}%
}}}}
\put(1666,-969){\makebox(0,0)[b]{\smash{{\SetFigFont{6}{7.2}{\rmdefault}{\mddefault}{\updefault}{\color[rgb]{0,0,0}4}%
}}}}
\put(2186,-546){\makebox(0,0)[b]{\smash{{\SetFigFont{6}{7.2}{\rmdefault}{\mddefault}{\updefault}{\color[rgb]{0,0,0}5}%
}}}}
\put(2565,-1409){\makebox(0,0)[b]{\smash{{\SetFigFont{6}{7.2}{\rmdefault}{\mddefault}{\updefault}{\color[rgb]{0,0,0}6}%
}}}}
\put(2506,-426){\makebox(0,0)[b]{\smash{{\SetFigFont{6}{7.2}{\rmdefault}{\mddefault}{\updefault}{\color[rgb]{0,0,0}7}%
}}}}
\put(2186,-43){\makebox(0,0)[b]{\smash{{\SetFigFont{6}{7.2}{\rmdefault}{\mddefault}{\updefault}{\color[rgb]{0,0,0}8}%
}}}}
\put(1586,-467){\makebox(0,0)[b]{\smash{{\SetFigFont{6}{7.2}{\rmdefault}{\mddefault}{\updefault}{\color[rgb]{0,0,0}9}%
}}}}
\put(1739,-37){\makebox(0,0)[b]{\smash{{\SetFigFont{6}{7.2}{\rmdefault}{\mddefault}{\updefault}{\color[rgb]{0,0,0}10}%
}}}}
\put(1946,200){\makebox(0,0)[b]{\smash{{\SetFigFont{6}{7.2}{\rmdefault}{\mddefault}{\updefault}{\color[rgb]{0,0,0}$+$}%
}}}}
\put(1493,-260){\makebox(0,0)[b]{\smash{{\SetFigFont{6}{7.2}{\rmdefault}{\mddefault}{\updefault}{\color[rgb]{0,0,0}$+$}%
}}}}
\put(2180,-313){\makebox(0,0)[b]{\smash{{\SetFigFont{6}{7.2}{\rmdefault}{\mddefault}{\updefault}{\color[rgb]{0,0,0}$-$}%
}}}}
\put(1840,-613){\makebox(0,0)[b]{\smash{{\SetFigFont{6}{7.2}{\rmdefault}{\mddefault}{\updefault}{\color[rgb]{0,0,0}$-$}%
}}}}
\put(2472,-1046){\makebox(0,0)[b]{\smash{{\SetFigFont{6}{7.2}{\rmdefault}{\mddefault}{\updefault}{\color[rgb]{0,0,0}$-$}%
}}}}
\put(352,440){\makebox(0,0)[rb]{\smash{{\SetFigFont{6}{7.2}{\rmdefault}{\mddefault}{\updefault}{\color[rgb]{0,0,0}$\varphi_9=1$}%
}}}}
\put(432,-706){\makebox(0,0)[lb]{\smash{{\SetFigFont{6}{7.2}{\rmdefault}{\mddefault}{\updefault}{\color[rgb]{0,0,0}$\varphi_4=1$}%
}}}}
\put(3207,-43){\makebox(0,0)[lb]{\smash{{\SetFigFont{6}{7.2}{\rmdefault}{\mddefault}{\updefault}{\color[rgb]{0,0,0}$\varphi_6=-1$}%
}}}}
\put(2572,721){\makebox(0,0)[b]{\smash{{\SetFigFont{6}{7.2}{\rmdefault}{\mddefault}{\updefault}{\color[rgb]{0,0,0}11}%
}}}}
\end{picture}%

%% file: KS.tex
{\def\nbpdfPostInput{\hfill}
  {\def\nbpdfPreOutput{\noindent\hspace{-7mm}}

\nbpdfInput{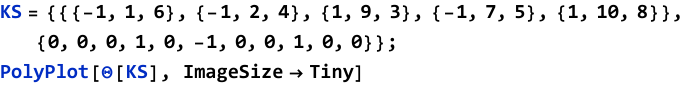}
\nbpdfOutput{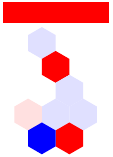}
  }

}